\documentclass{elsarticle}
\usepackage{graphicx}      %

\usepackage{amsthm}
\usepackage{times} %
\usepackage{amsmath} %
\usepackage{amssymb}  %
\usepackage{placeins}

\newcommand{\cay}{\mathrm{cay}}

\usepackage{color}
\usepackage{graphics} %
\usepackage{epsfig} %
\usepackage{subfigure}
\usepackage{ulem}
\usepackage{tikz}
\usetikzlibrary{patterns}
\usepackage{algorithmic}
\usepackage{algorithm}
\usepackage{tikz-cd}

\newcommand{\R}{\mathbb{R}}

\newcommand{\J}{\mathcal{J}}

\usepackage{bm}
\usepackage[colorlinks=true]{hyperref}
\hypersetup{urlcolor=blue, citecolor=red}
\newtheorem{theorem}{Theorem}[section]
\newtheorem{corollary}{Corollary}
\newtheorem{lemma}[theorem]{Lemma}
\newtheorem{proposition}{Proposition}

\newtheorem{definition}[theorem]{Definition}
\newtheorem{remark}{Remark}
\newtheorem{example}{Example}

\newcommand{\F}{\mathbb{F}}

\usepackage{amssymb}

\newcommand{\proa}{A^*G \mbox{$\;$}_{\tau^*} \kern-3pt\times_\alpha
G \mbox{$\;$}_\beta \kern-3pt\times_{\tau^*} A^*G}

\newcommand{\Ad}{\mbox{Ad}}

\newtheorem{assumption}[theorem]{Assumption}

\tikzstyle{vertex}=[circle,fill=black!20,minimum size=15pt,inner sep=0pt]
\tikzstyle{selected vertex} = [vertex, fill=red!24]
\tikzstyle{edge} = [draw,thick,-]
\tikzstyle{dedge} = [draw,thick,<->]
\tikzstyle{shadowdedge} = [draw, dotted,->]
\tikzstyle{weight} = [font=\small]
\tikzstyle{selected edge} = [draw,line width=3pt,-,red!50]
\tikzstyle{ignored edge} = [draw,line width=3pt,-,black!20]

\begin{document}

\begin{frontmatter}

\title{Variational integrators for time-dependent Lagrangian systems on Lie groups}%
\author[First]{Leonardo Colombo}
\author[Second]{Mar\'ia Emma Eyrea Iraz\'u}
\author[Third]{David Mart\'in de Diego}
\author[Second]{Mar\'ia Daniela S\'anchez}
\author[Second]{Marcela Zuccalli}
\address[First]{Centre for Automation and Robotics
(UPM-CSIC). Ctra. M300 Campo Real, Km 0,200, Arganda
del Rey, 28500 Madrid, Spain.}%
\address[Second]{Departamento de Matem\'atica, Centro de Matemática (CMaLP), Facultad de Ciencias Exactas, Universidad Nacional de La Plata, Argentina.}%
\address[Third]{Institute of Mathematical Sciences (CSIC-UAM-UCM-UC3M). Calle Nicol\'as Cabrera 13-15, Cantoblanco, 28049, Madrid, Spain.}%

\begin{abstract}
We develop fixed-step variational integrators for time-dependent
Lagrangian systems on Lie groups, derive their reduced dynamics,
and analyze their long-time behavior by backward error analysis. We derive the continuous
Euler-Lagrange equations, energy balance, and Noether
conservation laws, and construct time-indexed discrete
Lagrangians using retractions, with the midpoint-Cayley
scheme treated explicitly. On their regularity domains, the resulting phase-space maps are symplectic. For left-invariant discrete Lagrangians, they reduce
to Lie-Poisson maps that preserve coadjoint orbits and, when $G$
is connected, every Casimir function. A variational error analysis
establishes second-order accuracy of the
midpoint-Cayley method. We also construct an exactly
symplectic, projectable extension on
$T^*(\mathbb R\times G)$ whose formal modified
Hamiltonian is affine in the momentum conjugate to
time and admits symmetry reduction. Under uniform
analyticity, coefficient bounds, and domain-containment
assumptions, optimal truncation yields exponentially
accurate interpolation and near-conservation of the
modified extended Hamiltonian over exponentially long
intervals. For a connected attitude network on $SO(3)^N$ with damping and
uniformly positive, non-increasing coupling weights, we prove
local asymptotic stability of the static consensus set for the
continuous dynamics and exact preservation of total spatial
momentum by the discrete scheme. Numerical experiments on the
network and on a time-dependent rigid body on $SO(3)$ illustrate
second-order convergence and the extended energy balance.
\end{abstract}


\end{frontmatter}

\section{Introduction}

Geometric numerical integration seeks discretizations that reproduce the
qualitative structure of differential equations, rather than only accurate
pointwise approximations on short intervals. This is particularly important
in mechanics, where symplectic forms, momentum maps, first integrals,
and the geometry of the configuration space govern long-time behavior
\cite{hairer,Iserles2000}. Variational integrators are obtained from a
discrete Hamilton principle instead of discretizing the Euler-Lagrange
equations directly \cite{mawest}. Under the usual regularity and symmetry
assumptions, they generate symplectic maps and preserve discrete Noether
momenta; their long-time behavior can be studied by backward error analysis
\cite{hairer,vre}.

Explicit time dependence changes the energy balance. For a Lagrangian
$L:\mathbb R\times TQ\to\mathbb R$, the energy is generally not conserved;
its evolution is governed by $-\partial L/\partial t$. A geometric method
should therefore respect the variational structure and time-dependent
balance laws, rather than impose artificial energy conservation.
Time-dependent variational formulations and variational time integrators
have been studied from several perspectives
\cite{deLeonMartinDiego2002,LewMarsdenOrtizWest2004}. In particular,
\cite{gamonal} develops a fixed-step construction on an arbitrary
configuration manifold $Q$ using time-indexed discrete Lagrangians,
obtaining stepwise symplecticity, discrete Noether laws, and a
backward-error description. Higher-order prolongation-collocation
constructions have recently been considered in \cite{LiEtAl2026}.
Contact variational principles provide an alternative formulation for
dissipative systems
\cite{simoes2026variational,leocontacto,contact1,anahory2021geometry};
here we remain within Hamilton's principle for explicitly time-dependent
Lagrangians.

Many mechanical and engineering systems evolve not merely on manifolds
but on Lie groups: rotations on $SO(3)$, rigid motions on $SE(3)$,
and interacting bodies on product groups. Coordinate discretizations may
leave these spaces and require projection, whereas Lie group integrators
update the configuration intrinsically through group multiplication and
local parametrizations \cite{Iserles2000}. Variational mechanics on Lie
groups developed from the discrete rigid body of Moser and Veselov
\cite{MoserVeselov1991} through discrete Euler-Poincar\'e and
Lie-Poisson reduction \cite{MarsdenPekarskyShkoller2000},
Hamilton-Pontryagin methods \cite{BouRabeeMarsden2009}, and
variational integrators for rigid-body and multibody dynamics
\cite{LeeLeokMcClamroch2007,Arnold2016,BRULS2012121}.
Retractions and discretization maps systematically relate Lie algebra
increments to group configurations while accommodating symmetry and
reduction \cite{BarberoMartinDiego2023,BarberoMarreroMartinDiego2025}.

Time-dependent variational models also arise in formation control and
accelerated optimization. Noether symmetries have been related to
exponential decay in distance-based formation control \cite{leo4}.
Discrete variational principles have produced accelerated optimization
methods in Euclidean spaces \cite{CamposMahilloMartinDiego2023} and
momentum methods on Lie groups \cite{CamposMartinDiegoTorrente2024};
related Bregman formulations and adaptive variational schemes appear in
\cite{LeeTaoLeok2021,DuruisseauxLeok2023}. These contributions
illustrate the relevance of non-autonomous variational mechanics, but
their objectives differ from the combined geometric, reduction, and
long-time analysis pursued here.

The general-manifold setting of \cite{gamonal} is the closest
antecedent to our construction. Rather than reproducing that theory,
we develop the structures specific to a Lie group: adjoint and
coadjoint terms in the continuous and discrete equations, retractions
for relative group increments, and reduced Lie-Poisson dynamics.
We also distinguish pointwise accuracy of a discrete action from the
differentiable control needed to establish the order of its variational
map, and investigate how exact discrete symmetries constrain the
formal modified Hamiltonian. Finally, we formulate long-time
backward error estimates on the extended cotangent bundle, explicitly
stating the uniform analytic and domain assumptions they require.
Some preliminary geometric ingredients appeared in
\cite{MACI2025,Leo25}. The main contributions are as follows.

\begin{enumerate}
\item We formulate time-dependent Lagrangian mechanics on Lie groups in
left-trivialized variables, including the Euler-Lagrange equations,
energy balance, complete and vertical lifts, and Noether conservation
laws.

\item We develop time-indexed discrete Lagrangians on $G\times G$,
derive the discrete Euler-Lagrange equations variationally and through
discrete lifts, and establish discrete Noether conservation.
Retractions yield computable variational schemes, with the Cayley map
treated explicitly. We show that smooth variational contact of order $r$ yields a
phase-space map of order $r$ through a cancellation between the
two endpoint momentum defects (Theorem~\ref{thm:variational-order}).
For the midpoint-Cayley discrete Lagrangian, we verify this criterion
with $r=2$ (Corollary~\ref{cor:midpoint-cayley-order}).

\item The discrete Legendre transformations associated with a regular
discrete Lagrangian define a local symplectic map on $T^*G$.
If the discrete Lagrangian is left-invariant, this map descends
to a Lie-Poisson map on $\mathfrak g^*$. The reduced momentum
update is given by coadjoint transport and therefore preserves
coadjoint orbits. When $G$ is connected, it also preserves every
Casimir function (Theorem~\ref{thm:discrete-Lie-Poisson-reduction}).

\item We construct an exactly symplectic, $\rho$-projectable extension
of the fixed-step map on $T^*(\mathbb R\times G)$
(Proposition~\ref{prop:extended-map-Lie}). Its formal modified
Hamiltonian is affine in the momentum conjugate to time
(Theorem~\ref{thm:affine-modified-Hamiltonian}). Under left invariance and the hypotheses of
Theorem~\ref{thm:reduced-modified-Hamiltonian}, the formal
modified Hamiltonian descends to the reduced variables.
Under uniform analyticity, coefficient bounds, and trajectory
containment, Theorem~\ref{thm:BEA-Lie} gives an exponentially
small one-step interpolation defect and exponentially small
variation of the optimally truncated modified Hamiltonian over
an exponentially long interval. Corollary~\ref{cor:long-time-discrete-energy-balance}
expresses the resulting energy estimates in the physical variables
through the discrete conjugate-momentum update.
\end{enumerate}

The framework is illustrated first by a time-dependent rigid body
on $SO(3)$ and then by an interacting mechanical attitude network
on $SO(3)^N$ with time-dependent weights and a potential depending on
relative attitudes. The latter builds on geometric synchronization
and consensus models
\cite{NairLeonard2007,NairLeonard2008,SarletteSepulchre2009}.
For a connected network with damping and uniformly positive,
non-increasing coupling weights,
Theorem~\ref{thm:time-dependent-local-consensus} establishes local
asymptotic stability of the static consensus set for the continuous
dynamics. For the discrete integrator, diagonal $SO(3)$ symmetry
gives exact preservation of total spatial momentum
(Theorem~\ref{thm:network-noether}) and allows the formal modified
Hamiltonian to be expressed in terms of time, its conjugate
momentum, relative attitudes, and body momenta
(Proposition~\ref{prop:network-diagonal-modified-Hamiltonian}). Numerical experiments assess second-order
convergence of the complete extended map, including the
conjugate-time momentum, and illustrate geometric conservation
and long-time extended-Hamiltonian variation.

The paper is organized as follows. Section~\ref{Sec2-Lie} develops the
continuous theory; Section~\ref{Sec:DiscreteNoether}, the discrete
variational construction; and Section~\ref{sec:ham}, symplecticity, reduction, and
variational error analysis. Section~\ref{sec:BEA} develops the extended
Hamiltonian formulation and backward error analysis;
Section~\ref{sec:app-consensus} treats mechanical attitude consensus.
Section~\ref{sec: con} concludes the paper.

\section{Time-dependent Lagrangian mechanics on Lie groups}\label{Sec2-Lie}

This section recalls the continuous variational framework for
time-dependent Lagrangian systems on Lie groups. After introducing
left-trivialized variables, we derive the corresponding
Euler-Lagrange equations and establish the Lagrangian energy
balance. We then characterize the equations of motion using
complete and vertical lifts and formulate the associated Noether
conservation laws. The constructions are illustrated by a
time-dependent rigid body on $SO(3)$.

\subsection{Left trivialization and Euler–Lagrange equations}\label{subsec:bridge-Q-to-G}

Let $Q$ be the configuration space of a mechanical system, assumed to be a smooth manifold of dimension $n$, with local coordinates $q=(q^1,\dots,q^n)$. Its tangent and cotangent bundles are denoted by $TQ$ and $T^*Q$, respectively; they have local coordinates $(q^i,\dot q^i)$ and $(q^i,p_i)$ and dimension $2n$. For each $q\in Q$, the corresponding fibers are $T_qQ$ and $T_q^*Q$. The canonical projections are denoted by $\tau_Q:TQ\to Q$ and $\pi_Q:T^*Q\to Q$, with local expressions $\tau_Q(q^i,\dot q^i)=q^i$ and $\pi_Q(q^i,p_i)=q^i$.

Consider a time-dependent Lagrangian $L:\R\times TQ\to\mathbb{R}$, and let $\F L:\R\times TQ\to\R\times T^*Q$ be its Legendre transformation, defined locally by $(t,q,\dot q)\mapsto(t,q,p)$, where $p_i=\partial L/\partial\dot q^i$. We assume that $L$ is hyperregular, that is, that $\F L$ is a diffeomorphism. The velocities can then be written as functions $\dot q=\dot q(t,q,p)$, and the associated Hamiltonian $H:\mathbb{R}\times T^*Q\to\mathbb{R}$ is given by $H(t,q,p)=\langle p,\dot q(t,q,p)\rangle-L(t,q,\dot q(t,q,p))$.

The equations of motion follow from Hamilton's variational principle. Given $q_0,q_1\in Q$, let $C^2(q_0,q_1)=\{q:[0,T]\to Q\mid q\text{ is }C^2,\ q(0)=q_0,\ q(T)=q_1\}$ and define the action functional $\J:C^2(q_0,q_1)\to\R$ by $\J(q)=\int_0^T L(t,q,\dot q)\,dt$. Its critical points satisfy the Euler-Lagrange equations $\frac{d}{dt}(\partial L/\partial\dot q^i)-\partial L/\partial q^i=0$, or, equivalently,
\begin{equation}\label{euler-Lagrange}
\frac{\partial^2L}{\partial\dot q^i\partial\dot q^j}\ddot q^j
+\frac{\partial^2L}{\partial\dot q^i\partial q^j}\dot q^j
+\frac{\partial^2L}{\partial\dot q^i\partial t}
-\frac{\partial L}{\partial q^i}=0.
\end{equation}
See, for example, \cite{MarsdenRatiu1999}. Hyperregularity implies that the velocity Hessian $\operatorname{Hess}_{\dot q}(L)=(\partial^2L/\partial\dot q^i\partial\dot q^j)$ is nonsingular. Consequently, \eqref{euler-Lagrange} can be written as an explicit system of time-dependent second-order differential equations (see \cite{leo4,gamonal}). Notice that regularity is a local condition, whereas hyperregularity is the global assumption that makes the Legendre transformation invertible on the whole velocity phase space. Although the subsequent analysis is formulated directly on the Lagrangian side and therefore does not require a Hamiltonian reformulation, this distinction clarifies the relation between the tangent and cotangent bundle descriptions. In general, these equations do not admit closed-form solutions, and numerical integrators are therefore required to approximate their trajectories. 

When the configuration space is a Lie group $G$, its algebraic and geometric structure can be exploited in both the continuous and discrete formulations. Lie group methods provide intrinsic representations of motions such as rotations and translations, avoid the introduction of unnecessary local coordinates, and allow the group constraints and symmetries to be incorporated directly into the numerical method. At the same time, the classical constructions on $TQ$ must be adapted to respect the group multiplication, the Lie algebra representation of velocities, and the corresponding adjoint and coadjoint actions. These are not merely notational changes: they determine the form of the reduced variations, the equations of motion, and the momentum quantities associated with symmetry.

Integrators for Lie group systems and their convergence, particularly in multibody applications, are studied in \cite{Arnold2016,BRULS2012121}. In this paper, we instead construct integrators from a discrete variational principle. The resulting algebraic equations approximate the variational problem itself, rather than being obtained by a direct discretization of the differential equations. This variational construction is especially natural on a Lie group because the discrete trajectory remains on the configuration manifold by design, while the discrete increments can be represented in the Lie algebra by means of a retraction. The continuous framework developed below fixes the conventions that will later be inherited by that discrete formulation.

We now specialize to a finite-dimensional Lie group $G$ with Lie algebra $\mathfrak g=T_eG$. Let $\mathcal L_g:G\to G$ denote left translation, $\mathcal L_g(h)=gh$. We use the left trivialization $\lambda:TG\to G\times\mathfrak g$, defined by $\lambda(v_g)=(g,\xi)$ with $\xi=T_g\mathcal L_{g^{-1}}v_g$. Thus, for a curve $g(t)$, the body velocity is $\xi(t)=T_{g(t)}\mathcal L_{g(t)^{-1}}\dot g(t)$, and the inverse relation is the reconstruction equation
\begin{equation}\label{eq:kinematics-left}
\dot g=T_e\mathcal L_g\,\xi
\qquad
(\text{for matrix groups, }\dot g=g\xi).
\end{equation}
Accordingly, a time-dependent Lagrangian $L:\mathbb{R}\times TG\to\mathbb{R}$ will be regarded as a smooth function $L:\mathbb{R}\times G\times\mathfrak g\to\mathbb{R}$, $(t,g,\xi)\mapsto L(t,g,\xi)$. All the constructions below are expressed in these variables, with \eqref{eq:kinematics-left} linking the group trajectory to its body velocity.

To formulate the variational equations in left-trivialized
variables, we also identify the cotangent bundle $T^*G$
with $G\times\mathfrak{g}^*$ through the corresponding
dual trivialization $\lambda^*:T^*G\to G\times\mathfrak g^*$ defined by $\lambda^*(\alpha_g)=(g,\mu)$, where $\mu=(T_e\mathcal L_g)^*\alpha_g\in\mathfrak g^*$. Conversely, $(g,\mu)$ determines the covector $\alpha_g=(T_g\mathcal L_{g^{-1}})^*\mu\in T_g^*G$.

Let $\partial L/\partial g(t,g,\xi)\in T_g^*G$ and
$\partial L/\partial\xi(t,g,\xi)\in\mathfrak{g}^*$
denote the partial derivatives of the left-trivialized
Lagrangian with respect to $g$ and $\xi$, respectively,
with the other variable held fixed.
We shall use the left-trivialized covector
\begin{equation}\label{eq:DgL}
\mathbf D_gL(t,g,\xi)
:=(T_e\mathcal L_g)^*\!\left(\frac{\partial L}{\partial g}(t,g,\xi)\right)
\in\mathfrak g^*.
\end{equation}

Let $g(\cdot)$ be a smooth curve in $G$ and consider variations $g^\varepsilon(t)=g(t)\exp(\varepsilon\eta(t))$, where $\eta:[0,T]\to\mathfrak g$ satisfies $\eta(0)=\eta(T)=0$. Their infinitesimal variation is $\delta g(t)=T_e\mathcal L_{g(t)}\eta(t)$. The induced variation of the body velocity is the standard identity (see, for instance, \cite{MarsdenRatiu1999})
\begin{equation}\label{eq:delta-xi}
\delta\xi(t)=\dot\eta(t)+[\xi(t),\eta(t)].
\end{equation}
We denote by $\operatorname{ad}_\xi:\mathfrak g\to\mathfrak g$ the map $\operatorname{ad}_\xi\eta=[\xi,\eta]$ and adopt the dual convention
\begin{equation}\label{eq:adstar}
\left\langle\operatorname{ad}_\xi^*\mu,\eta\right\rangle
=\left\langle\mu,[\xi,\eta]\right\rangle,
\qquad
\mu\in\mathfrak g^*,\quad\eta\in\mathfrak g.
\end{equation}
In particular, $\langle\operatorname{ad}_\xi^*\mu,\xi\rangle=\langle\mu,[\xi,\xi]\rangle=0$. Here, $\operatorname{ad}_{\xi}^*$ denotes the algebraic
dual of $\operatorname{ad}_{\xi}$, without an additional
minus sign. We make this convention explicit because
different sign conventions for the infinitesimal
coadjoint action are used in the literature. All continuous and discrete equations in this paper use the convention stated here.

For later use, we also fix the corresponding group-level coadjoint
convention. For $g\in G$, the coadjoint action
$\operatorname{Ad}_g^*:\mathfrak g^*\to\mathfrak g^*$ is defined by
$\langle\operatorname{Ad}_g^*\mu,\eta\rangle
=\langle\mu,\operatorname{Ad}_{g^{-1}}\eta\rangle$ for every
$\mu\in\mathfrak g^*$ and $\eta\in\mathfrak g$. With the convention \eqref{eq:adstar}, its infinitesimal generator
satisfies
$\left.\frac{d}{ds}\right|_{s=0}
\operatorname{Ad}_{\exp(s\xi)}^*\mu=-\operatorname{ad}_\xi^*\mu$.

To obtain the equations of motion, apply \eqref{eq:delta-xi} to the action $\J(g)=\int_0^T L\,dt$, with $L$ evaluated along $(t,g(t),\xi(t))$. Its first variation is
\begin{align*}
\delta\J
={}\int_0^T\left(
\left\langle\mathbf D_gL,\eta\right\rangle
+\left\langle\frac{\partial L}{\partial\xi},\dot\eta+[\xi,\eta]\right\rangle
\right)dt 
={}\int_0^T\left\langle
\mathbf D_gL
-\frac{d}{dt}\frac{\partial L}{\partial\xi}
+\operatorname{ad}_\xi^*\frac{\partial L}{\partial\xi},
\eta\right\rangle dt,
\end{align*}
where the boundary term vanishes because $\eta(0)=\eta(T)=0$. Since $\eta$ is otherwise arbitrary, Hamilton's principle gives the left-trivialized Euler-Lagrange equations (see \cite{Leo25}):
\begin{equation}\label{eq:EL-Lie}
\frac{d}{dt}\frac{\partial L}{\partial\xi}(t,g,\xi)
-\operatorname{ad}_\xi^*\frac{\partial L}{\partial\xi}(t,g,\xi)
-\mathbf D_gL(t,g,\xi)=0,
\qquad
\dot g=T_e\mathcal L_g\,\xi.
\end{equation}

When $L$ is left-invariant, its left-trivialized
representation is independent of $g$. In this case,
the preceding equations reduce to the
time-dependent Euler-Poincar\'e equations
\[
\frac{d}{dt}\frac{\partial\ell}{\partial\xi}
=
\operatorname{ad}_{\xi}^*
\frac{\partial\ell}{\partial\xi},
\qquad
\ell(t,\xi):=L(t,e,\xi),
\]
together with the reconstruction equation $\dot g=T_e\mathcal L_g\,\xi$.

\subsection{Energy balance and variational characterization}

The time-dependent velocity phase space $\mathbb{R}\times TG$ inherits canonical geometric structures from $TG$ that are independent of coordinates. It is the product of time with the tangent bundle, rather than a tangent bundle over an extended configuration space, and we shall therefore refer to it simply as the time-dependent velocity phase space. In particular, the Liouville vector field $\Delta$ extends to
$\R\times TG$ by acting trivially on the time component and provides
an intrinsic definition of the Lagrangian energy.
Under the left trivialization \(TG\simeq G\times\mathfrak g\), the Liouville vector field is the infinitesimal generator of fiber dilations, namely \(\Delta_{(g,\xi)}=\left.\frac{d}{ds}\right|_{s=0}(g,(1+s)\xi)\). If \(\{e_a\}\) is a basis of \(\mathfrak g\) and \(\xi=\xi^a e_a\), then \(\Delta=\xi^a\partial/\partial\xi^a\). %

The energy associated with $L$ is defined by
\begin{equation}\label{eq:energy-def}
E_L:=\Delta(L)-L,
\qquad
E_L(t,g,\xi)
=\left\langle\frac{\partial L}{\partial\xi}(t,g,\xi),\xi\right\rangle
-L(t,g,\xi).
\end{equation}
Unlike in the autonomous case, this energy is generally not conserved. Its evolution is described by the following proposition.

\begin{proposition}\label{prop:energy-balance}
Let $L:\mathbb{R}\times G\times\mathfrak g\to\mathbb{R}$ be a smooth time-dependent Lagrangian, and let $(g(t),\xi(t))$ be a solution of \eqref{eq:EL-Lie}. Then
\[
\frac{d}{dt}E_L(t,g(t),\xi(t))
=-\frac{\partial L}{\partial t}(t,g(t),\xi(t)).
\]
\end{proposition}

\begin{proof}
Differentiating \eqref{eq:energy-def} along a solution gives
\begin{align*}
\frac{d}{dt}E_L
={}&
\left\langle\frac{d}{dt}\left(\frac{\partial L}{\partial\xi}\right),\xi\right\rangle
+\left\langle\frac{\partial L}{\partial\xi},\dot\xi\right\rangle
-\frac{\partial L}{\partial t}
-\left\langle\frac{\partial L}{\partial g},\dot g\right\rangle
-\left\langle\frac{\partial L}{\partial\xi},\dot\xi\right\rangle \\
={}&
\left\langle\frac{d}{dt}\left(\frac{\partial L}{\partial\xi}\right),\xi\right\rangle
-\frac{\partial L}{\partial t}
-\left\langle\frac{\partial L}{\partial g},\dot g\right\rangle.
\end{align*}
By \eqref{eq:kinematics-left} and \eqref{eq:DgL},
\[
\left\langle\frac{\partial L}{\partial g},\dot g\right\rangle
=\left\langle\frac{\partial L}{\partial g},T_e\mathcal L_g\xi\right\rangle
=\left\langle\mathbf D_gL,\xi\right\rangle.
\]
Pairing \eqref{eq:EL-Lie} with $\xi$ yields
\[
\left\langle\frac{d}{dt}\left(\frac{\partial L}{\partial\xi}\right),\xi\right\rangle
=\left\langle\operatorname{ad}_\xi^*\left(\frac{\partial L}{\partial\xi}\right),\xi\right\rangle
+\left\langle\mathbf D_gL,\xi\right\rangle
=\left\langle\mathbf D_gL,\xi\right\rangle,
\]
where the last equality follows from \eqref{eq:adstar}. Substitution into the preceding identity completes the proof.
\end{proof}

The Euler-Lagrange equations also admit an intrinsic
characterization in terms of the standard complete
and vertical lifts of vector fields on $G$.
This characterization provides the link between
the variational equations and the Noether conservation
laws developed below. Let $X\in\mathfrak{X}(G)$
and define its body representation $x:G\to\mathfrak{g}$ by
\begin{equation}\label{eq:body-rep}
x(g):=T_g\mathcal L_{g^{-1}}X(g),
\qquad
\text{equivalently,}\qquad
X(g)=T_e\mathcal L_g\,x(g).
\end{equation}

\begin{definition}\label{def:vertical-lift-Lie}
Let $X\in\mathfrak X(G)$. Its vertical lift
$X^V\in\mathfrak X(TG)$, under the identification
$TG\simeq G\times\mathfrak g$, is given by
$X^V(g,\xi):=(0_g,x(g))$. We extend it to $\mathbb R\times TG$ by setting $X^V(t,g,\xi)=(0_t,X^V(g,\xi))$. 
\end{definition}
Observe that it corresponds to the left trivialization of the standard vertical lift to the tangent bundle.

Consequently, for any $L:\mathbb{R}\times G\times\mathfrak g\to\mathbb{R}$,
\begin{equation}\label{eq:Xv-on-L}
X^V(L)(t,g,\xi)
:=\left\langle\frac{\partial L}{\partial\xi}(t,g,\xi),x(g)\right\rangle.
\end{equation}

The complete lift contains both terms appearing in the variation formula \eqref{eq:delta-xi}. The time component of the extended lift is zero, so it produces no term involving $\partial L/\partial t$.

\begin{proposition}\label{prop:complete-lift-Lie}
Let $X\in\mathfrak X(G)$, and let $x$ be its body representation in \eqref{eq:body-rep}. Under the identification $TG\simeq G\times\mathfrak g$, the complete lift $X^C\in\mathfrak X(TG)$ is
\begin{equation}\label{eq:complete-lift-Lie}
X^C(g,\xi)
:=\Bigl(X(g),\mathbf dx(g)[T_e\mathcal L_g\xi]+[\xi,x(g)]\Bigr),
\end{equation}
where $\mathbf dx(g):T_gG\to\mathfrak g$ is the differential of $x$.
\end{proposition}

\begin{proof}
Let $g(t)$ be a curve with body velocity $\xi(t)=T_{g(t)}\mathcal L_{g(t)^{-1}}\dot g(t)$, and consider the variation induced by $X$, namely $\delta g(t)=X(g(t))$. Its left-trivialized variation is $\eta(t)=T_{g(t)}\mathcal L_{g(t)^{-1}}\delta g(t)=x(g(t))$. Hence
\[
\dot\eta(t)=\mathbf dx(g(t))[\dot g(t)]
=\mathbf dx(g(t))[T_e\mathcal L_{g(t)}\xi(t)].
\]
Using \eqref{eq:delta-xi}, we obtain $\delta\xi=\mathbf dx(g)[T_e\mathcal L_g\xi]+[\xi,x(g)]$. This is the vertical component in \eqref{eq:complete-lift-Lie}, while the base component is $\delta g=X(g)$. Together they give the tangent lift of the flow generated by $X$, which is precisely the complete lift.
\end{proof}

We extend $X^C$ to $\mathbb{R}\times TG$ by $X^C(t,g,\xi)=(0_t,X^C(g,\xi))$. Its action on $L$ is
\begin{equation}\label{eq:Xc-on-L}
\begin{split}
X^C(L)(t,g,\xi)
={}&\left\langle\frac{\partial L}{\partial g}(t,g,\xi),X(g)\right\rangle +\left\langle\frac{\partial L}{\partial\xi}(t,g,\xi),
\mathbf dx(g)[T_e\mathcal L_g\xi]+[\xi,x(g)]\right\rangle.
\end{split}
\end{equation}
Equivalently, using \eqref{eq:DgL} and \eqref{eq:adstar},
\begin{align*}
X^C(L)(t,g,\xi)
={}\langle\mathbf D_gL(t,g,\xi),x(g)\rangle
+\left\langle\frac{\partial L}{\partial\xi}(t,g,\xi),
\mathbf dx(g)[T_e\mathcal L_g\xi]\right\rangle +\left\langle\operatorname{ad}_\xi^*\left(\frac{\partial L}{\partial\xi}(t,g,\xi)\right),x(g)\right\rangle.
\end{align*}

The following result expresses, in left-trivialized
variables, the classical characterization of
Euler-Lagrange solutions through complete and
vertical lifts.

\begin{theorem}\label{thm:EL-via-lifts-Lie}
A curve $g(\cdot)$, or equivalently $(g(\cdot),\xi(\cdot))$, satisfies the left-trivialized Euler-Lagrange equations \eqref{eq:EL-Lie} if and only if
\begin{equation}\label{eq:EL-lifts-Lie}
\frac{d}{dt}\Bigl(X^V(L)(t,g(t),\xi(t))\Bigr)
=X^C(L)(t,g(t),\xi(t)),
\qquad
\text{for every }X\in\mathfrak X(G).
\end{equation}
\end{theorem}

\begin{proof}
Let $X\in\mathfrak X(G)$, and let $x$ be its body representation. From \eqref{eq:Xv-on-L}, $X^V(L)=\langle\partial L/\partial\xi,x(g)\rangle$. Along $(g(t),\xi(t))$,
\[
\frac{d}{dt}X^V(L)
=\left\langle\frac{d}{dt}\left(\frac{\partial L}{\partial\xi}\right),x(g)\right\rangle
+\left\langle\frac{\partial L}{\partial\xi},\mathbf dx(g)[\dot g]\right\rangle,
\qquad
\dot g=T_e\mathcal L_g\xi.
\]
On the other hand, \eqref{eq:Xc-on-L} gives
\[
X^C(L)
=\left\langle\frac{\partial L}{\partial g},X(g)\right\rangle
+\left\langle\frac{\partial L}{\partial\xi},\mathbf dx(g)[\dot g]\right\rangle
+\left\langle\frac{\partial L}{\partial\xi},[\xi,x(g)]\right\rangle.
\]
Subtracting these identities and using \eqref{eq:DgL} and \eqref{eq:adstar}, we obtain
\[
\frac{d}{dt}X^V(L)-X^C(L)
=\left\langle
\frac{d}{dt}\left(\frac{\partial L}{\partial\xi}\right)
-\mathbf D_gL
-\operatorname{ad}_\xi^*\left(\frac{\partial L}{\partial\xi}\right),
x(g)\right\rangle.
\]
Since $X$, and hence $x(g)$, is arbitrary, \eqref{eq:EL-lifts-Lie} holds for every $X$ if and only if the $\mathfrak g^*$-valued expression in the pairing vanishes. This is precisely the dynamical equation in \eqref{eq:EL-Lie}; the reconstruction equation \eqref{eq:kinematics-left} is already built into the definition of $\xi$.
\end{proof}

\subsection{Noether symmetries and momentum maps}\label{sec2.3}

The characterization established in Theorem~\ref{thm:EL-via-lifts-Lie}
allows us to relate symmetries to
conserved quantities. We restrict attention to
symmetries generated by time-independent vector
fields on $G$, acting trivially on the time variable. A symmetry leaves the Lagrangian
invariant, whereas a Noether symmetry may change
it by a total time derivative. Such a boundary term
does not affect the Euler-Lagrange equations and
may depend explicitly on time. We formulate the
corresponding Noether theorem in the left-trivialized
setting, following the standard time-dependent
variational framework \cite{gamonal}.

For $f\in C^\infty(\mathbb{R}\times G)$, its total time derivative along a point $(t,g,\xi)$ of the time-dependent velocity phase space is
\[
d_Tf(t,g,\xi)
:=\frac{\partial f}{\partial t}(t,g)
+\left\langle\frac{\partial f}{\partial g}(t,g),T_e\mathcal L_g\xi\right\rangle
=\frac{\partial f}{\partial t}(t,g)
+\left\langle(T_e\mathcal L_g)^*\frac{\partial f}{\partial g}(t,g),\xi\right\rangle.
\]

\begin{definition}\label{def:symmetry-Lie}
A vector field $X\in\mathfrak X(G)$ is a \textit{symmetry of the Lagrangian} $L$ if $X^C(L)=0$.
\end{definition}

\begin{definition}\label{def:noether-symmetry-Lie}
A vector field $X\in\mathfrak X(G)$ is a \textit{Noether symmetry} of $L$ if there exists a function $f\in C^\infty(\mathbb{R}\times G)$ such that $X^C(L)=d_Tf$.
\end{definition}

The next result relates these symmetries to constants of motion.

\begin{theorem}\label{thm:noether-Lie}
If $X$ is a Noether symmetry of $L$ with associated function $f$, then $J_X:=X^V(L)-f$ is a constant of motion of the Euler-Lagrange equations \eqref{eq:EL-Lie}. Equivalently, $\frac{d}{dt}J_X(t,g(t),\xi(t))=0$ along every solution.
\end{theorem}

\begin{proof}
Along a solution of \eqref{eq:EL-Lie}, Theorem~\ref{thm:EL-via-lifts-Lie} gives $\frac{d}{dt}(X^V(L))=X^C(L)$. If $X$ is a Noether symmetry, then $X^C(L)=d_Tf$, and therefore
\[
\frac{d}{dt}\bigl(X^V(L)-f\bigr)
=X^C(L)-d_Tf=0.
\]
Thus $J_X=X^V(L)-f$ is conserved along solutions.
\end{proof}

We next specialize Theorem~\ref{thm:noether-Lie} to symmetries arising from a Lie group action. This also prepares the rigid body example below. 
Let $\tilde{G}$ be a Lie group with Lie algebra $\mathfrak{\tilde{g}}$, acting smoothly on $G$ from the left by $\Phi:\tilde{G}\times G\to G$. We write $\Phi_a(g)=\Phi(a,g)$, for $a \in \tilde{G}$. The induced tangent action is $T\Phi_a:TG\to TG$, $T\Phi_a(v_g)=T_g\Phi_a(v_g)$. %

\begin{definition}\label{def2.6}
The action of $\tilde{G}$ is a \textit{symmetry of the time-dependent Lagrangian} $L:\mathbb{R}\times TG\to\mathbb{R}$ if
\[
L(t,T\Phi_a(v_g))=L(t,v_g),
\qquad
a\in \tilde{G},\quad v_g\in T_gG,\quad t\in\mathbb{R}.
\]
\end{definition}

For $\zeta\in\widetilde{\mathfrak g}$, differentiating the group action
along the one-parameter subgroup $\exp(s\zeta)$ defines the infinitesimal
generator
\begin{equation}\label{eq:inf-generator-Lie}
\zeta_G(g)
=
\left.\frac{d}{ds}\right|_{s=0}
\Phi_{\exp(s\zeta)}(g).
\end{equation}

Let $\phi_s:=\Phi_{\exp(s\zeta)}$. The corresponding action on $TG$ is the tangent-lifted action
$v_g\mapsto T\phi_s(v_g)$. Its infinitesimal generator is the complete
lift $\zeta_G^C\in\mathfrak X(TG)$ of the vector field
$\zeta_G\in\mathfrak X(G)$ defined in \eqref{eq:inf-generator-Lie};
equivalently,
\[
\zeta_G^C(v_g)
=
\left.\frac{d}{ds}\right|_{s=0}
T\phi_s(v_g).
\]

If the action of $\widetilde G$ is a symmetry of $L$, Definition~\ref{def2.6} gives $L\bigl(t,T\phi_s(v_g)\bigr)=L(t,v_g)$
for every $s$. Hence the map
$s\mapsto L\bigl(t,T\phi_s(v_g)\bigr)$ is constant. Differentiating at $s=0$ and using the chain rule yields
\[
\begin{aligned}
0
&=
\left.\frac{d}{ds}\right|_{s=0}
L\bigl(t,T\phi_s(v_g)\bigr)
=
dL_{(t,v_g)}
\left(
0,
\left.\frac{d}{ds}\right|_{s=0}T\phi_s(v_g)
\right)
=
dL_{(t,v_g)}
\bigl(0,\zeta_G^C(v_g)\bigr)
\\
&=
\zeta_G^C(L)(t,v_g).
\end{aligned}
\]
Therefore,
\begin{equation}\label{eq:inf-invariance-Lie}
\zeta_G^C(L)=0,
\qquad
\zeta\in\widetilde{\mathfrak g}.
\end{equation}

Thus invariance of $L$ under the group action implies the infinitesimal
invariance condition \eqref{eq:inf-invariance-Lie}. Conversely, if
$\widetilde G$ is connected, \eqref{eq:inf-invariance-Lie} for every
$\zeta\in\widetilde{\mathfrak g}$ implies invariance under the full
group action.

For $\zeta\in\widetilde{\mathfrak g}$, let $x_\zeta(g):=T_gL_{g^{-1}}\zeta_G(g)\in\mathfrak g$ be the body representation of the infinitesimal generator. The corresponding momentum map $J:\mathbb R\times G\times\mathfrak g\to\widetilde{\mathfrak g}^{\,*}$ is defined by $\langle J(t,g,\xi),\zeta\rangle:=\langle\partial L/\partial\xi,x_\zeta(g)\rangle$. This is the Lagrangian momentum map associated
with the lifted action, obtained by composing
the canonical momentum map on $T^*G$ with
the time-dependent Legendre transformation. By Theorem~\ref{thm:noether-Lie}, every component $\langle J,\zeta\rangle$ is constant along the Euler-Lagrange flow; hence $J$ is conserved along solutions. In the particular case in which $\widetilde G=G$ acts on itself by left multiplication, $\Phi(a,g)=ag$, one has $\zeta_G(g)=T_eR_g\zeta$ and $x_\zeta(g)=\operatorname{Ad}_{g^{-1}}\zeta$. Therefore $\langle J,\zeta\rangle=\langle\operatorname{Ad}_g^*(\partial L/\partial\xi),\zeta\rangle$, and the spatial momentum $\Pi:=\operatorname{Ad}_g^*(\partial L/\partial\xi)$ is conserved.

We now illustrate these results with a time-dependent
rigid body on $\mathrm{SO}(3)$. The example makes
explicit the distinction between the Lagrangian
energy, the mechanical kinetic energy, and the
spatial momentum preserved by Noether's theorem.

\begin{example}\label{ex:rigid-body-SO3}
 Let $G=SO(3)$. 
We identify $\mathfrak{so}(3)$ with $\mathbb{R}^3$ through the hat map $\widehat{\cdot}:\mathbb{R}^3\rightarrow\mathfrak{so}(3)$, \newline 
$\displaystyle{x=
\begin{pmatrix}
x_1\\x_2\\x_3
\end{pmatrix}
\mapsto
\widehat{x}:=
\begin{pmatrix}
0 & -x_3 & x_2\\
x_3 & 0 & -x_1\\
-x_2 & x_1 & 0
\end{pmatrix}}$, characterized by $\widehat{x}\,y=x\times y$. Its inverse is the vee map
$(\cdot)^\vee:\mathfrak{so}(3)\rightarrow\mathbb{R}^3$, $(\widehat{x})^\vee=x$. Under this identification,
$[\widehat{\xi},\widehat{\eta}]
=
\widehat{\xi\times\eta}$, and we shall identify elements of $\mathfrak{so}(3)$ with their
vector representatives in $\mathbb{R}^3$.

Consider the time-dependent Lagrangian
\[
L(t,R,\xi)=\frac{e^{-\kappa t}}{2}\langle\xi,\mathbb I\xi\rangle,
\qquad
\kappa<0,
\]
where $\mathbb I:\mathfrak{so}(3)\to\mathfrak{so}(3)^*$ is the symmetric positive-definite inertia operator. The vector $\xi$ is the body angular velocity, determined by $\widehat\xi=R^{-1}\dot R$. Since $L$ is independent of $R$, $\partial L/\partial R=0$ and hence $\mathbf D_RL=0$. Its derivative with respect to the body angular velocity is $\partial L/\partial\xi=e^{-\kappa t}\mathbb I\xi$, and therefore
\[
\frac{d}{dt}\frac{\partial L}{\partial\xi}
=e^{-\kappa t}\mathbb I\dot\xi
-\kappa e^{-\kappa t}\mathbb I\xi.
\]
Substitution into \eqref{eq:EL-Lie} gives
\[
e^{-\kappa t}\mathbb I\dot\xi
-\kappa e^{-\kappa t}\mathbb I\xi
-\operatorname{ad}_\xi^*\bigl(e^{-\kappa t}\mathbb I\xi\bigr)=0,
\qquad
\dot R=R\widehat\xi.
\]

With the dual convention \eqref{eq:adstar}, the identification $\mathfrak{so}(3)^*\simeq\mathbb{R}^3$ gives $\operatorname{ad}_\xi^*\mu=\mu\times\xi=-\xi\times\mu$. After cancelling the nonzero factor $e^{-\kappa t}$, the equations of motion take the vector form
\begin{equation}\label{eq:rigid-body-continuous}
\mathbb I\dot\xi
+\xi\times(\mathbb I\xi)
-\kappa\mathbb I\xi=0,
\qquad
\dot R=R\widehat\xi.
\end{equation}
This is the free rigid body equation with an additional linear dissipation term generated by the explicit time dependence of the Lagrangian (see \cite{leocontacto} and the references therein). Thus, although \eqref{eq:rigid-body-continuous} coincides with the equation of a rigid body subject to linear Rayleigh dissipation, with damping proportional to the body angular momentum, no external nonconservative force has been introduced here.

The associated energy is
\[
E_L(t,R,\xi)
=\left\langle\frac{\partial L}{\partial\xi},\xi\right\rangle-L(t,R,\xi)
=\frac{e^{-\kappa t}}{2}\langle\xi,\mathbb I\xi\rangle.
\]
By Proposition~\ref{prop:energy-balance},
\[
\frac{d}{dt}E_L
=-\frac{\partial L}{\partial t}
=\frac{\kappa e^{-\kappa t}}{2}\langle\xi,\mathbb I\xi\rangle
=\kappa E_L.
\]
Hence $E_L(t)=E_L(t_0)e^{\kappa(t-t_0)}$, and in particular $E_L(t)=E_L(0)e^{\kappa t}$ when $t_0=0$.
Thus the Lagrangian energy decays exponentially when $\kappa<0$. The mechanical kinetic energy $K(\xi)=\frac12\langle\xi,\mathbb I\xi\rangle$ is related to $E_L$ by $E_L=e^{-\kappa t}K$. Pairing \eqref{eq:rigid-body-continuous} with $\xi$ also gives $\dot K=2\kappa K$, because $\langle\xi\times(\mathbb I\xi),\xi\rangle=0$. Thus the physical kinetic energy decays as $K(t)=K(t_0)e^{2\kappa(t-t_0)}$, consistently with the dissipative interpretation, and in particular $K(t)=K(0)e^{2\kappa t}$ when $t_0=0$. Although the prefactor $e^{-\kappa t}$ grows for $\kappa<0$, the kinetic energy $K$ decays at twice the exponential rate, so that $E_L=e^{-\kappa t}K$ still decays as $e^{\kappa t}$.

Finally, $L$ is invariant under the left action of $SO(3)$ on itself, $\Phi(Q,R)=QR$. Indeed, the tangent action sends $(R,\dot R)$ to $(QR,Q\dot R)$ and leaves the body velocity unchanged, since $(QR)^{-1}(Q\dot R)=R^{-1}\dot R$. For $\zeta\in\mathfrak{so}(3)$, the infinitesimal generator is $\zeta_G(R)=\widehat\zeta R$, whose body representation is $x(R)=\operatorname{Ad}_{R^{-1}}\zeta$. Noether's theorem therefore gives the conserved quantity
\[
J_\zeta
=\left\langle e^{-\kappa t}\mathbb I\xi,
\operatorname{Ad}_{R^{-1}}\zeta\right\rangle
=\left\langle\operatorname{Ad}_R^*\bigl(e^{-\kappa t}\mathbb I\xi\bigr),\zeta\right\rangle.
\]
Since this holds for every $\zeta$, the exponentially weighted spatial angular momentum
\[
\Pi:=\operatorname{Ad}_R^*\bigl(e^{-\kappa t}\mathbb I\xi\bigr)
\qquad
(\text{equivalently, }\Pi=e^{-\kappa t}R\mathbb I\xi\text{ in vector form})
\]
is constant. Notice that the corresponding body quantity $e^{-\kappa t}\mathbb I\xi$ is generally not constant; it evolves by the coadjoint term in \eqref{eq:EL-Lie}.

The continuous energy balance and Noether conservation
laws established in this section provide the geometric
reference for the discrete variational framework
developed next. In particular, the preservation of
symmetry-induced momenta and the approximation of the
time-dependent energy balance will play distinct roles
in the numerical analysis.
\end{example}

\section{Time-dependent discrete variational mechanics on Lie groups}
\label{Sec:DiscreteNoether}

The continuous variational framework established in
Section~\ref{Sec2-Lie} provides the geometric reference for its
fixed-step discretization. We first recall 
 discrete variational principles for time-dependent Lagrangians on a
general configuration manifold and then specialize
it to Lie groups, introducing the corresponding
discrete Legendre transformations. We subsequently
characterize the discrete time-dependent Euler-Lagrange equations
through discrete lifts and establish the associated
Noether conservation laws. Finally, we develop
retraction-based Lie group variational integrators,
with particular emphasis on the Cayley map and its
application to the time-dependent rigid body.

\subsection{Discrete variational principle and Legendre transformations}\label{subsec:discrete-variational-manifolds}

Let $Q$ be a smooth configuration manifold and consider a uniform temporal
grid $t_k=t_0+kh$, for $k=0,\dots,N$. A time-dependent discrete Lagrangian is
a family of smooth functions
$L_d^k:Q\times Q\to\mathbb{R}$, $k=0,\dots,N-1$, where $L_d^k(q_k,q_{k+1})$
approximates the action of a continuous regular time-dependent Lagrangian
$L:\mathbb{R}\times TQ\to\mathbb{R}$ over the interval
$[t_k,t_{k+1}]$ \cite{CamposMahilloMartinDiego2023}. The dependence on $k$ accounts for the explicit time
dependence of the continuous system; the dependence on the fixed step size
$h$ is suppressed from the notation.

A discrete path is a sequence
$q_d=(q_0,\dots,q_N)\in Q^{N+1}$. For fixed endpoints $q_0$ and $q_N$,
the discrete action functional is
\[
\J_d(q_d)
:=
\sum_{k=0}^{N-1}L_d^k(q_k,q_{k+1}).
\]
Variations of $q_d$ are sequences
$\delta q_d=(\delta q_0,\dots,\delta q_N)$, with
$\delta q_k\in T_{q_k}Q$ and $\delta q_0=\delta q_N=0$. The first
variation of the discrete action is
\[
\delta\J_d(q_d)
=
\sum_{k=1}^{N-1}
\left\langle
D_2L_d^{k-1}(q_{k-1},q_k)
+
D_1L_d^k(q_k,q_{k+1}),
\delta q_k
\right\rangle.
\]
Therefore, the stationarity condition $\delta\J_d(q_d)=0$ for every
admissible variation is equivalent to the time-dependent discrete
Euler-Lagrange equations
\[
D_2L_d^{k-1}(q_{k-1},q_k)
+
D_1L_d^k(q_k,q_{k+1})
=
0,
\qquad
k=1,\dots,N-1.
\]

A discrete Lagrangian is said to be regular if the mixed derivative
$D_{12}L_d^k$ is nondegenerate. Under this assumption, the discrete
Euler-Lagrange equations determine locally $q_{k+1}$ from
$(q_{k-1},q_k)$ and define a family of local discrete flows
\[
\Psi_d^k:Q\times Q\rightarrow Q\times Q,
\qquad
(q_{k-1},q_k)\longmapsto(q_k,q_{k+1}).
\]

The exact discrete Lagrangian provides the fundamental connection between the
continuous and discrete variational theories \cite{hairer,mawest}. Assume that
$h$ is sufficiently small and that the endpoints $q_k$ and $q_{k+1}$ are
sufficiently close for the Euler-Lagrange boundary-value problem to admit a
unique solution. The exact discrete Lagrangian associated with the interval
$[t_k,t_{k+1}]$ is then defined by
\[
L_{d,E}^k(q_k,q_{k+1};h)
:=
\int_{t_k}^{t_{k+1}}
L\bigl(t,q(t),\dot q(t)\bigr)\,dt,
\]
where $q(t)$ is the solution of the continuous Euler-Lagrange equations
satisfying $q(t_k)=q_k$ and $q(t_{k+1})=q_{k+1}$. In contrast with the
autonomous case, the exact discrete Lagrangian depends not only on $h$ and the
endpoints, but also on the initial time $t_k$, here encoded by the index $k$.

The discrete Euler-Lagrange equations associated with the family
$L_{d,E}^k$ reproduce exactly the continuous evolution sampled at the grid
points. Since the exact discrete Lagrangian is generally unavailable in
closed form, practical variational integrators are obtained by constructing
approximations
$L_d^k\simeq L_{d,E}^k$ using quadrature rules and suitable approximations of
the continuous trajectory between $q_k$ and $q_{k+1}$ (see \cite{mawest}).

\medskip

We now specialize the preceding construction to a finite-dimensional Lie
group $G$ with Lie algebra $\mathfrak g$. The product $G\times G$ is the
discrete analogue of the velocity phase space $TG$: a pair
$(g_k,g_{k+1})$ represents a discrete trajectory segment. From the Lie groupoid viewpoint, $G\times G
\rightrightarrows G$ is the pair groupoid,
whose discrete Lagrangian mechanics is included
in the general framework developed in
\cite{MarreroMartinDeDiegoMartinez2006}. The corresponding
Lie algebra representation of this segment will be introduced later by means
of a retraction. This construction is consistent with the continuous
left-trivialized formulation introduced in
Subsection~\ref{subsec:bridge-Q-to-G}
(see \cite{Leo25}).

A time-dependent discrete Lagrangian on $G$ is a family of smooth functions
$L_d^k:G\times G\to\mathbb{R}$, for $k=0,\dots,N-1$. For a discrete
trajectory $g_d=(g_0,\dots,g_N)\in G^{N+1}$, its discrete action is
\begin{equation}
\label{eq:disc-action-G}
\J_d(g_d)
=
\sum_{k=0}^{N-1}L_d^k(g_k,g_{k+1}).
\end{equation}

Variations compatible with the group structure can be written as
$g_k^\varepsilon=g_k\exp(\varepsilon\eta_k)$, where
$\eta_k\in\mathfrak g$ and $\eta_0=\eta_N=0$. Their infinitesimal form is
$\delta g_k=T_e\mathcal L_{g_k}\eta_k\in T_{g_k}G$. Since these variations
span $T_{g_k}G$ at every interior point, stationarity of
\eqref{eq:disc-action-G} gives the discrete Euler-Lagrange equations
\begin{equation}
\label{eq:DEL-G}
D_2L_d^{k-1}(g_{k-1},g_k)
+
D_1L_d^k(g_k,g_{k+1})
=
0,
\qquad
k=1,\dots,N-1.
\end{equation}
If the discrete Lagrangians are regular, these equations define local
discrete flows
\[
\Psi_d^k:G\times G\rightarrow G\times G,
\qquad
(g_{k-1},g_k)\longmapsto(g_k,g_{k+1}).
\]

The discrete Euler-Lagrange equations admit a phase-space interpretation
through the discrete Legendre transformations. At an interior node $g_k$, we
define the pre- and post-momenta by
\begin{equation}
\label{eq:prepost-momenta-G}
\begin{aligned}
p_k^-
&:=
-D_1L_d^k(g_k,g_{k+1})
\in T_{g_k}^*G,\\
p_k^+
&:=
D_2L_d^{k-1}(g_{k-1},g_k)
\in T_{g_k}^*G.
\end{aligned}
\end{equation}
Equivalently, these covectors are obtained from the discrete Legendre
transformations
\[
\begin{aligned}
\mathbb{F}^-L_d^k(g_k,g_{k+1})
&:=
\bigl(g_k,-D_1L_d^k(g_k,g_{k+1})\bigr),\\
\mathbb{F}^+L_d^{k-1}(g_{k-1},g_k)
&:=
\bigl(g_k,D_2L_d^{k-1}(g_{k-1},g_k)\bigr).
\end{aligned}
\]
Equation~\eqref{eq:DEL-G} is precisely the momentum-matching condition
$p_k^-=p_k^+$, for $k=1,\dots,N-1$. Under the regularity assumption, this
condition determines the discrete evolution in phase space.

To connect this formulation with the Lie algebra description used in the
continuous theory, we left-trivialize the covectors in
\eqref{eq:prepost-momenta-G}. Define
\begin{equation}
\label{eq:left-triv-covector}
\mu_k^-
:=
(T_e\mathcal L_{g_k})^*p_k^-
\in\mathfrak g^*,
\qquad
\mu_k^+
:=
(T_e\mathcal L_{g_k})^*p_k^+
\in\mathfrak g^*.
\end{equation}
The discrete Euler-Lagrange equations are therefore equivalent to
$\mu_k^-=\mu_k^+$. Along a solution, we write their common value as
$\mu_k$. Thus, the natural left-trivialized phase variables are
$(g_k,\mu_k)\in G\times\mathfrak g^*$. A discrete body velocity
$\xi_k\in\mathfrak g$ will be associated with the pair
$(g_k,g_{k+1})$ later, after choosing a retraction from $\mathfrak g$ to $G$.

\subsection{Discrete lifts, symmetries, and Noether theorem}\label{sec2.2}

To formulate symmetries and conservation laws intrinsically,
we first introduce discrete analogues of the complete and vertical lifts
used in the continuous theory. These objects provide an intrinsic
characterization of the discrete Euler-Lagrange equations and will be used
below to formulate the discrete Noether theorem. Although the following
constructions are valid on any smooth configuration manifold, we state them on
$G$ in order to remain within the Lie group setting of the preceding
subsection. 
These constructions are closely related to the
discrete variational calculus on $Q\times Q$
developed in \cite{DiazMartinDeDiego2018}.

For $(g_0,g_1)\in G\times G$, the canonical product decomposition gives
\[
T_{(g_0,g_1)}(G\times G)
\simeq T_{g_0}G\oplus T_{g_1}G,
\qquad
T^*_{(g_0,g_1)}(G\times G)
\simeq T^*_{g_0}G\oplus T^*_{g_1}G.
\]
Accordingly, the differential of a discrete Lagrangian
$L_d^k:G\times G\to\mathbb{R}$ is represented by
$dL_d^k(g_0,g_1)=(D_1L_d^k(g_0,g_1),D_2L_d^k(g_0,g_1))$, where
$D_1L_d^k(g_0,g_1)\in T_{g_0}^*G$ and
$D_2L_d^k(g_0,g_1)\in T_{g_1}^*G$.

\begin{definition}
\label{def:discrete-lifts}
Let $X\in\mathfrak X(G)$. Its \textit{discrete complete lift} and its first
and second \textit{discrete vertical lifts} are the vector fields on
$G\times G$ defined by
\begin{equation*}
\label{eq:lifts-discrete}
\begin{aligned}
X^{C,d}(g_0,g_1)
&:=
\bigl(X(g_0),X(g_1)\bigr),\\
X^{V_1,d}(g_0,g_1)
&:=
\bigl(X(g_0),0_{g_1}\bigr),\\
X^{V_2,d}(g_0,g_1)
&:=
\bigl(0_{g_0},X(g_1)\bigr).
\end{aligned}
\end{equation*}
In particular,
$X^{C,d}=X^{V_1,d}+X^{V_2,d}$ (see \cite{MarreroMartinDeDiegoMartinez2006}).
\end{definition}

If $\varphi_s$ denotes the local flow of $X$, then $X^{C,d}$ is the
infinitesimal generator of the diagonal action
$(g_0,g_1)\mapsto(\varphi_s(g_0),\varphi_s(g_1))$, whereas
$X^{V_1,d}$ and $X^{V_2,d}$ generate the corresponding actions on the
first and second factors, respectively. Their action on $L_d^k$ is therefore
given by
\begin{align*}
X^{C,d}(L_d^k)(g_0,g_1)
&=
\left\langle D_1L_d^k(g_0,g_1),X(g_0)\right\rangle
+
\left\langle D_2L_d^k(g_0,g_1),X(g_1)\right\rangle,\\
X^{V_1,d}(L_d^k)(g_0,g_1)
&=
\left\langle D_1L_d^k(g_0,g_1),X(g_0)\right\rangle,\\
X^{V_2,d}(L_d^k)(g_0,g_1)
&=
\left\langle D_2L_d^k(g_0,g_1),X(g_1)\right\rangle.
\end{align*}

The following result is the discrete counterpart of the continuous
lift-based characterization established in
Theorem~\ref{thm:EL-via-lifts-Lie}.

\begin{theorem}
\label{thm:DEL-via-discrete-lifts}
Let $L_d^k:G\times G\to\mathbb{R}$, $k=0,\dots,N-1$, be a family of
time-dependent discrete Lagrangians. A discrete curve
$g_d=(g_0,\dots,g_N)$ satisfies the discrete Euler-Lagrange equations
\eqref{eq:DEL-G} if and only if, for every $X\in\mathfrak X(G)$,
\begin{equation}
\label{eq:DEL-lifts-identity}
X^{V_2,d}(L_d^{k-1})(g_{k-1},g_k)
+
X^{V_1,d}(L_d^k)(g_k,g_{k+1})
=
0,
\qquad
k=1,\dots,N-1.
\end{equation}
Equivalently,
\begin{equation}
\label{eq:DEL-lifts-complete}
\begin{split}
X^{C,d}(L_d^{k-1})(g_{k-1},g_k)
={}&
X^{V_1,d}(L_d^{k-1})(g_{k-1},g_k)\\
&-
X^{V_1,d}(L_d^k)(g_k,g_{k+1}),
\qquad
k=1,\dots,N-1.
\end{split}
\end{equation}
\end{theorem}

\begin{proof}
By the definitions of the two discrete vertical lifts, the left-hand side of
\eqref{eq:DEL-lifts-identity} is
\[
\left\langle
D_2L_d^{k-1}(g_{k-1},g_k)
+
D_1L_d^k(g_k,g_{k+1}),
X(g_k)
\right\rangle.
\]
If $g_d$ satisfies the discrete Euler-Lagrange equations
\eqref{eq:DEL-G}, this expression vanishes for every
$X\in\mathfrak X(G)$.

Conversely, suppose that \eqref{eq:DEL-lifts-identity} holds for every
$X\in\mathfrak X(G)$. Since the value $X(g_k)\in T_{g_k}G$ can be chosen
arbitrarily, it follows that
\[
D_2L_d^{k-1}(g_{k-1},g_k)
+
D_1L_d^k(g_k,g_{k+1})
=
0.
\]
These are precisely the discrete Euler-Lagrange equations
\eqref{eq:DEL-G}.

Finally, since
$X^{C,d}=X^{V_1,d}+X^{V_2,d}$, equation
\eqref{eq:DEL-lifts-identity} can be rewritten as
\eqref{eq:DEL-lifts-complete}.
\end{proof}

With this intrinsic characterization in place, we turn to
symmetries of time-dependent discrete Lagrangian systems and their
associated conservation laws. A symmetry of the configuration manifold induces
a diagonal action on $G\times G$, and invariance of every member of the family
$L_d^k$ gives rise to a discrete momentum that is preserved by the discrete
Euler-Lagrange evolution. This provides the discrete counterpart of the
Noether theorem established in Section~\ref{Sec2-Lie}. 
Unlike the more general Noether symmetries considered
in Subsection~\ref{sec2.3}, we restrict the discrete formulation
to Lie group actions that leave each time-indexed
discrete Lagrangian invariant. This assumption
is sufficient for the discrete momentum conservation
laws established below.

Let $\tilde G$ be a Lie group with Lie algebra $\mathfrak {\tilde g}$, acting smoothly on $G$
from the left by
$\Phi:\tilde G\times G\to G$. We write $\Phi_a(g)=\Phi(a,g),$ for $a\in \tilde G$. For
$\zeta\in\mathfrak {\tilde g}$, let $\zeta_G\in\mathfrak X(G)$ be the infinitesimal
generator defined in \eqref{eq:inf-generator-Lie}. The induced diagonal action
on $G\times G$ is
$\Phi_a^{G\times G}(g_0,g_1)
=
\bigl(\Phi_a(g_0),\Phi_a(g_1)\bigr)$,
and its infinitesimal generator is the discrete complete lift
\[
\zeta_{G\times G}(g_0,g_1)
=
\bigl(\zeta_G(g_0),\zeta_G(g_1)\bigr)
=
\zeta_G^{C,d}(g_0,g_1).
\]

\begin{definition}
\label{def:discrete-symmetry}
The action of $\tilde G$ is a \textit{symmetry of the time-dependent discrete
Lagrangian} if
\[
L_d^k\bigl(\Phi_a(g_0),\Phi_a(g_1)\bigr)
=
L_d^k(g_0,g_1)
\]
for every $a\in \tilde G$, $(g_0,g_1)\in G\times G$, and
$k=0,\dots,N-1$.
\end{definition}

Differentiating the invariance condition along the one-parameter subgroup
$\exp(s\zeta)$ gives the infinitesimal condition
\begin{equation}
\label{eq:discrete-inf-invariance}
\zeta_G^{C,d}(L_d^k)=0,
\qquad
\zeta\in\mathfrak {\tilde g},\quad k=0,\dots,N-1.
\end{equation}
Invariance always implies \eqref{eq:discrete-inf-invariance}; if $\tilde G$ is
connected, the two conditions are equivalent.

For each $k$, define the pre- and post-discrete momentum maps
$J_{d,-}^k,J_{d,+}^k:G\times G\to\mathfrak {\tilde g}^*$ by
\begin{equation*}
\label{eq:discrete-momentum-maps}
\begin{aligned}
\left\langle J_{d,-}^k(g_k,g_{k+1}),\zeta\right\rangle
&:=
-\left\langle
D_1L_d^k(g_k,g_{k+1}),
\zeta_G(g_k)
\right\rangle,\\
\left\langle J_{d,+}^k(g_k,g_{k+1}),\zeta\right\rangle
&:=
\left\langle
D_2L_d^k(g_k,g_{k+1}),
\zeta_G(g_{k+1})
\right\rangle,
\end{aligned}
\end{equation*}
for every $\zeta\in\mathfrak {\tilde g}$. By
\eqref{eq:prepost-momenta-G}, these expressions can also be written as
\[
\begin{aligned}
\left\langle J_{d,-}^k(g_k,g_{k+1}),\zeta\right\rangle
&=
\left\langle p_k^-,\zeta_G(g_k)\right\rangle,\\
\left\langle J_{d,+}^k(g_k,g_{k+1}),\zeta\right\rangle
&=
\left\langle p_{k+1}^+,\zeta_G(g_{k+1})\right\rangle.
\end{aligned}
\]

If $L_d^k$ is invariant, then
\[
0
=
\zeta_G^{C,d}(L_d^k)
=
\left\langle
D_1L_d^k,\zeta_G(g_k)
\right\rangle
+
\left\langle
D_2L_d^k,\zeta_G(g_{k+1})
\right\rangle.
\]
Consequently, $J_{d,-}^k=J_{d,+}^k$. We denote their common value by
$J_d^k$ and, for each $\zeta\in\mathfrak {\tilde g}$, define its scalar component by
\begin{equation}
\label{eq:disc-noether-quantity}
J_\zeta^k
:=
\left\langle J_d^k(g_k,g_{k+1}),\zeta\right\rangle
=
-\left\langle
D_1L_d^k(g_k,g_{k+1}),
\zeta_G(g_k)
\right\rangle.
\end{equation}

\begin{theorem}
\label{thm:discrete-noether}
Suppose that $\tilde G$ is a symmetry group of the family of discrete Lagrangians
$L_d^k$. If $g_d=(g_0,\dots,g_N)$ satisfies the discrete Euler-Lagrange
equations \eqref{eq:DEL-G}, then
\[
J_d^{k+1}(g_{k+1},g_{k+2})
=
J_d^k(g_k,g_{k+1}),
\qquad
k=0,\dots,N-2.
\]
Equivalently, $J_\zeta^{k+1}=J_\zeta^k$ for every
$\zeta\in\mathfrak {\tilde g}$.
\end{theorem}

\begin{proof}
Let $\zeta\in\mathfrak {\tilde g}$. Applying
\eqref{eq:DEL-lifts-complete} to the infinitesimal generator $\zeta_G$ gives
\[
\begin{split}
\zeta_G^{C,d}(L_d^k)(g_k,g_{k+1})
={}&
\zeta_G^{V_1,d}(L_d^k)(g_k,g_{k+1})-
\zeta_G^{V_1,d}(L_d^{k+1})(g_{k+1},g_{k+2}).
\end{split}
\]
By the infinitesimal invariance condition
\eqref{eq:discrete-inf-invariance}, the left-hand side vanishes. Therefore,
\[
-\zeta_G^{V_1,d}(L_d^k)(g_k,g_{k+1})
=
-\zeta_G^{V_1,d}(L_d^{k+1})(g_{k+1},g_{k+2}),
\]
which, by \eqref{eq:disc-noether-quantity}, is precisely
$J_\zeta^k=J_\zeta^{k+1}$.

Equivalently, the same conclusion follows directly from momentum matching:
\[
\begin{aligned}
\left\langle J_d^k,\zeta\right\rangle=
\left\langle p_{k+1}^+,\zeta_G(g_{k+1})\right\rangle=
\left\langle p_{k+1}^-,\zeta_G(g_{k+1})\right\rangle=
\left\langle J_d^{k+1},\zeta\right\rangle,
\end{aligned}
\]
where the second equality follows from the discrete Euler-Lagrange equations.
Since this holds for every $\zeta\in\mathfrak {\tilde g}$, the result follows.
\end{proof}

To express the discrete Noether quantity in left-trivialized variables, define
the body representation of the infinitesimal generator by
$x_\zeta(g)
:=
T_g\mathcal L_{g^{-1}}\zeta_G(g)
\in\mathfrak g$. Using the left-trivialized momenta introduced in
\eqref{eq:left-triv-covector}, we obtain
\[
\left\langle J_d^k,\zeta\right\rangle
=
\left\langle\mu_k^-,x_\zeta(g_k)\right\rangle
=
\left\langle\mu_{k+1}^+,x_\zeta(g_{k+1})\right\rangle.
\]
Thus, for a general action, the conserved quantity depends both on the
left-trivialized momentum and on the body representation of the corresponding
infinitesimal generator.

We now consider the particular case in which $\tilde G=G$ acts on itself by left
multiplication, $\Phi_a(g)=ag$. For $\zeta\in\mathfrak g$, the corresponding
infinitesimal generator is
$\displaystyle{\zeta_G(g)
=
\left.\frac{d}{ds}\right|_{s=0}\exp(s\zeta)g}$,
and its body representation is
$x_\zeta(g)=\operatorname{Ad}_{g^{-1}}\zeta$. Therefore,
\[
\left\langle J_d^k,\zeta\right\rangle
=
\left\langle
\mu_k^-,
\operatorname{Ad}_{g_k^{-1}}\zeta
\right\rangle
=
\left\langle
\operatorname{Ad}_{g_k}^*\mu_k^-,
\zeta
\right\rangle.
\]

Along a solution, $\mu_k^-=\mu_k^+=:\mu_k$, and the discrete Noether theorem
implies that the spatial momentum
$\Pi_k
:=
\operatorname{Ad}_{g_k}^*\mu_k$ is preserved:
$\Pi_{k+1}=\Pi_k$.

If $F_k=g_k^{-1}g_{k+1}$, this conservation law is equivalently expressed as
the coadjoint transport equation
\[
\mu_k
=
\operatorname{Ad}_{F_k}^*\mu_{k+1}.
\]
This is the discrete counterpart of the continuous coadjoint evolution and
will be used below when the group increment $F_k$ is parametrized through a
retraction.

\subsection{Retractions and Lie group variational integrators}
\label{subsec:retractions}

We now express the discrete dynamics in Lie algebra variables by means of a
retraction (see e.g. \cite{Absil, BarberoMartinDiego2023}). This provides a local representation of the group increment
$g_k^{-1}g_{k+1}$ and leads to a computable formulation of the discrete
Euler-Lagrange equations on the left-trivialized phase space
$G\times\mathfrak g^*$. We use throughout the coadjoint convention fixed in
Section~\ref{Sec2-Lie}. We shall also use the standard notions of coadjoint orbit and Casimir
function on $\mathfrak g^*$; their precise definitions and their relation
with the Lie-Poisson structure are recalled in
Section~\ref{sec:ham}.

\begin{definition}
A local retraction on $G$ is a smooth map
$\mathcal R:U\subset\mathfrak g\to G$, where $U$
is an open neighborhood of $0$, satisfying
$\mathcal R(0)=e$, $T_0\mathcal R=\mathrm{id}_{\mathfrak g}$.
By the inverse function theorem, after restricting
$U$ if necessary, $\mathcal R$ is a diffeomorphism
onto a neighborhood of $e$. We denote its local
inverse by $\mathcal R^{-1}$. For each $g\in G$,
the corresponding left-translated retraction is
$\mathcal R_g(\xi):=g\mathcal R(\xi)$.
\end{definition}

Under the left trivialization $TG\simeq G\times\mathfrak g$, determined by
$\dot g=T_e\mathcal L_g\xi$, the continuous reconstruction equation is
locally approximated by
\begin{equation}
\label{eq:disc-retraction-step}
g_{k+1}=g_k\mathcal R(h\xi_k),
\qquad
\xi_k=\frac1h\mathcal R^{-1}\bigl(g_k^{-1}g_{k+1}\bigr).
\end{equation}
We denote the group increment by
$F_k:=g_k^{-1}g_{k+1}=\mathcal R(h\xi_k)$.

Let $L:\mathbb R\times G\times\mathfrak g\to\mathbb R$ be the
left-trivialized representation of a time-dependent Lagrangian. A
retraction-based approximation of the midpoint configuration is
$g_{k+\frac12}:=g_k\mathcal R(\frac h2\xi_k)$, with
$t_{k+\frac12}:=t_k+\frac h2$. The midpoint discrete Lagrangian is defined by
\begin{equation}
\label{eq:Ld-midpoint-retraction}
L_d^k(g_k,g_{k+1})
=
hL\left(t_{k+\frac12},g_{k+\frac12},\xi_k\right),
\end{equation}
where $\xi_k$ is determined by \eqref{eq:disc-retraction-step}. The midpoint construction does not automatically
yield a second-order method for an arbitrary
retraction. Its order depends on the approximation
properties of the chosen retraction and the
resulting discrete Lagrangian. The second-order
accuracy of the midpoint-Cayley scheme will
be established in Section~\ref{sec:ham}.

It is convenient to introduce the mixed and Lie algebra representations
\begin{equation}
\label{eq:Ld-mixed-variables}
\widehat L_d^k(g_k,F_k)
:=
L_d^k(g_k,g_kF_k),
\qquad
\ell_d^k(g_k,\xi_k)
:=
\widehat L_d^k\bigl(g_k,\mathcal R(h\xi_k)\bigr).
\end{equation}

For an arbitrary retraction $\mathcal R$, let
$\mathcal R^{-1}$ denote its local inverse. We define the
left-trivialized differential of $\mathcal R^{-1}$ at $F$ by
$\mathrm d^L\mathcal R^{-1}_{F}
:=
T_F\mathcal R^{-1}\circ T_e\mathcal L_F
:
\mathfrak g\rightarrow\mathfrak g$.
Thus,
\[
\mathrm d^L\mathcal R^{-1}_{F}(\eta)
=
T_F\mathcal R^{-1}\bigl(T_e\mathcal L_F\,\eta\bigr),
\qquad
\eta\in\mathfrak g.
\]
Equivalently, it is the differential of $\mathcal R^{-1}$ evaluated
on the tangent vector at $F$ represented by $\eta$ under left
trivialization.

\begin{proposition}
\label{prop:general-retraction-momentum-update}
Let $\mathcal R:U\subset\mathfrak g\to G$ be a retraction, and let
$\widehat L_d^k$ and $\ell_d^k$ be defined by
\eqref{eq:Ld-mixed-variables}, with
$F_k=\mathcal R(h\xi_k)$. Assume that $F_k$ belongs to the domain of
the local inverse $\mathcal R^{-1}$. Define the left-trivialized derivative
$\mathbf D_g\widehat L_d^k(g_k,F_k)
:=
\bigl(T_e\mathcal L_{g_k}\bigr)^*
D_1\widehat L_d^k(g_k,F_k)$\footnote{
Here, differentiation with respect to $g_k$ in
the mixed representation is performed while
holding the group increment $F_k$ fixed.
} and the post-momentum by
\[
\mu_{k+1}
:=
\frac{1}{h}
\left(
\mathrm d^L\mathcal R^{-1}_{F_k}
\right)^*
\frac{\partial\ell_d^k}{\partial\xi}(g_k,\xi_k).
\]
Then the discrete Euler-Lagrange equations are equivalent to
\[
\mu_k
=
\operatorname{Ad}_{F_k}^*\mu_{k+1}
-
\mathbf D_g\widehat L_d^k(g_k,F_k),
\qquad
g_{k+1}=g_kF_k,
\qquad
F_k=\mathcal R(h\xi_k).
\]
If $L_d^k$ is invariant under the diagonal left action, equivalently,
if $\widehat L_d^k$ is independent of $g_k$, then
$\mu_k=\operatorname{Ad}_{F_k}^*\mu_{k+1}$, or equivalently
$\mu_{k+1}=\operatorname{Ad}_{F_k^{-1}}^*\mu_k$.
Consequently, the momentum evolves by coadjoint transport and remains
on the same coadjoint orbit. If $G$ is connected, every Casimir
function is therefore preserved.
\end{proposition}

\begin{proof}
Consider variations
$\delta g_k=T_e\mathcal L_{g_k}\eta_k$, where
$\eta_0=\eta_N=0$. Since $F_k=g_k^{-1}g_{k+1}$, the induced variation is
\[
\delta F_k
=
F_k\left(
\eta_{k+1}-\operatorname{Ad}_{F_k^{-1}}\eta_k
\right).
\]
Using $F_k=\mathcal R(h\xi_k)$ and differentiating
$h\xi_k=\mathcal R^{-1}(F_k)$ gives
\[
\delta\xi_k
=
\frac{1}{h}
\mathrm d^L\mathcal R^{-1}_{F_k}
\left(
\eta_{k+1}
-
\operatorname{Ad}_{F_k^{-1}}\eta_k
\right).
\]
Therefore, the variation of one term of the discrete action is
\[
\delta\widehat L_d^k
=
\left\langle
\mathbf D_g\widehat L_d^k
-
\operatorname{Ad}_{F_k}^*\mu_{k+1},
\eta_k
\right\rangle
+
\left\langle
\mu_{k+1},\eta_{k+1}
\right\rangle,
\]
where we used the coadjoint convention fixed in
Section~\ref{Sec2-Lie} to write
$\langle\mu_{k+1},\operatorname{Ad}_{F_k^{-1}}\eta_k\rangle
=
\langle\operatorname{Ad}_{F_k}^*\mu_{k+1},\eta_k\rangle$.

Summing over $k=0,\dots,N-1$, reindexing the last term, and using
$\eta_0=\eta_N=0$ gives
\[
\delta\mathcal J_d
=
\sum_{k=1}^{N-1}
\left\langle
\mu_k
+
\mathbf D_g\widehat L_d^k(g_k,F_k)
-
\operatorname{Ad}_{F_k}^*\mu_{k+1},
\eta_k
\right\rangle.
\]
Since the interior variations are arbitrary, stationarity is equivalent
to the stated momentum equation. If $L_d^k$ is invariant under the diagonal left action, then its mixed
representation $\widehat L_d^k$ is independent of $g_k$. Hence
$\mathbf D_g\widehat L_d^k=0$, and the coadjoint transport relation
follows.
\end{proof}

The choice of retraction determines the relation
between discrete velocity and momentum, but
does not alter the coadjoint form of the
momentum update in the left-invariant case.
The corresponding Lie-Poisson structure
will be analyzed in Section~\ref{sec:ham}.

If the exponential map itself is chosen as the retraction, then
$\mathcal R^{-1}=\log$ locally and
$\mathrm d^L\mathcal R^{-1}_{F}$ is the corresponding
left-trivialized differential of the logarithm. In the remainder of
this subsection, we instead specialize the general construction to
the Cayley retraction. This construction applies to matrix Lie groups
for which the Cayley transform maps the Lie algebra locally into the
group, in particular to $SO(n)$. The Cayley map and its local inverse
are
\begin{equation*}
\label{eq:cayley-and-inverse}
\operatorname{cay}(\Xi)
=
\left(I-\frac12\Xi\right)^{-1}
\left(I+\frac12\Xi\right),
\qquad
\operatorname{cay}^{-1}(F)
=
2(F-I)(F+I)^{-1}.
\end{equation*}
assuming that both $I-\frac12\Xi$ and $I+\frac12\Xi$ are
invertible.
The map is defined whenever the indicated inverses exist and satisfies
$\operatorname{cay}(-\Xi)=\operatorname{cay}(\Xi)^{-1}$. Hence
$F_k=\operatorname{cay}(h\xi_k)$ and
$\xi_k=h^{-1}\operatorname{cay}^{-1}(F_k)$.

Differentiating the inverse Cayley map gives
\begin{equation*}
\label{eq:dcay-inv}
d(\operatorname{cay}^{-1})_F(\delta F)
=
4(F+I)^{-1}\delta F(F+I)^{-1}.
\end{equation*}

We denote the corresponding left-trivialized inverse differential by
\begin{equation}
\label{eq:dlog-cayley}
\mathrm d^L\operatorname{cay}^{-1}_{F}(\eta)
:=
T_F\operatorname{cay}^{-1}
\bigl(T_e\mathcal L_F\,\eta\bigr)
=
4(F+I)^{-1}F\eta(F+I)^{-1}.
\end{equation}

For the midpoint discrete Lagrangian \eqref{eq:Ld-midpoint-retraction}, its
Lie algebra representation becomes
\begin{equation*}
\label{eq:ell-d-cayley}
\ell_d^k(g_k,\xi_k)
:=
hL\left(
t_{k+\frac12},
g_k\operatorname{cay}\left(\frac h2\xi_k\right),
\xi_k
\right).
\end{equation*}

The Cayley-based equations are now an immediate specialization of
Proposition~\ref{prop:general-retraction-momentum-update}.

\begin{corollary}
\label{prop:cayley-momentum-update}
Let $\widehat L_d^k$ and $\ell_d^k$ be defined by
\eqref{eq:Ld-mixed-variables}, with
$F_k=\operatorname{cay}(h\xi_k)$. Then
\begin{equation}
\label{eq:force-and-postmomentum}
\mu_{k+1}
:=
\frac{1}{h}
\left(
\mathrm d^L\operatorname{cay}^{-1}_{F_k}
\right)^*
\frac{\partial\ell_d^k}{\partial\xi}(g_k,\xi_k),
\end{equation}
and the discrete Euler-Lagrange equations are equivalent to
\begin{equation}
\label{eq:LGVI-momentum}
\mu_k
=
\operatorname{Ad}_{F_k}^*\mu_{k+1}
-
\mathbf D_g\widehat L_d^k(g_k,F_k),
\qquad
g_{k+1}=g_kF_k,
\qquad
F_k=\operatorname{cay}(h\xi_k).
\end{equation}
If $L_d^k$ is invariant under the diagonal left action, then
\begin{equation}
\label{eq:discrete-Euler-Poincare}
\mu_k=\operatorname{Ad}_{F_k}^*\mu_{k+1},
\qquad
\mu_{k+1}=\operatorname{Ad}_{F_k^{-1}}^*\mu_k.
\end{equation}
\end{corollary}

\begin{proof}
This follows directly from
Proposition~\ref{prop:general-retraction-momentum-update} by taking
$\mathcal R=\operatorname{cay}$ and using
\eqref{eq:dlog-cayley}.
\end{proof}

The index in the definition of $\mu_{k+1}$ is essential: it is the
left-trivialized post-momentum associated with the segment
$(g_k,g_{k+1})$. More precisely,
\begin{equation}
\label{eq:mixed-pre-post-momenta}
\begin{aligned}
\mu_{k+1}
&=
\bigl(T_e\mathcal L_{g_{k+1}}\bigr)^*
D_2L_d^k(g_k,g_{k+1}),\\
\operatorname{Ad}_{F_k}^*\mu_{k+1}
-
\mathbf D_g\widehat L_d^k(g_k,F_k)
&=
-
\bigl(T_e\mathcal L_{g_k}\bigr)^*
D_1L_d^k(g_k,g_{k+1}).
\end{aligned}
\end{equation}
Thus \eqref{eq:LGVI-momentum} is precisely the left-trivialized
momentum-matching condition introduced in
\eqref{eq:left-triv-covector}.

\begin{example}
\label{ejemplo-discretizado}

Consider the time-dependent rigid body Lagrangian introduced in
Example~\ref{ex:rigid-body-SO3},
\begin{equation*}
\label{eq:rigid-continuous-system}
L(t,R,\xi)
=
\frac12e^{-\kappa t}
\left\langle\xi,\mathbb I\xi\right\rangle,
\qquad
\frac{d}{dt}
\left(e^{-\kappa t}\mathbb I\xi\right)
=
\operatorname{ad}_\xi^*
\left(e^{-\kappa t}\mathbb I\xi\right),
\end{equation*}
where $R\in SO(3)$,
$\xi\in\mathfrak{so}(3)\simeq\mathbb R^3$, and
$\mathbb I$ is a symmetric positive-definite inertia operator.

The corresponding energy
$E(t)=\frac12e^{-\kappa t}\langle\xi(t),\mathbb I\xi(t)\rangle$
satisfies $\dot E=\kappa E$. Therefore,
\begin{equation}
\label{eq:continuous-energy-law}
E(t)=E(0)e^{\kappa t},
\end{equation}
and the energy decays when $\kappa<0$.

For the discrete system, let
$F_k=R_k^{-1}R_{k+1}=\operatorname{cay}(h\widehat\xi_k)$. The reconstruction
equation is
\begin{equation*}
\label{eq:R-update-explicit}
R_{k+1}
=
R_kF_k
=
R_k
\left(I-\frac h2\widehat\xi_k\right)^{-1}
\left(I+\frac h2\widehat\xi_k\right).
\end{equation*}
The midpoint discrete Lagrangian is
\begin{equation}
\label{eq:Ld-rigid}
L_d^k(R_k,R_{k+1})
=
\frac h2e^{-\kappa t_{k+\frac12}}
\left\langle\xi_k,\mathbb I\xi_k\right\rangle.
\end{equation}
It depends on $(R_k,R_{k+1})$ only through
$F_k=R_k^{-1}R_{k+1}$. It is therefore invariant under the diagonal left
action $(R_k,R_{k+1})\mapsto(QR_k,QR_{k+1})$, and its mixed representation
is independent of $R_k$.

Since
$\partial\ell_d^k/\partial\xi
=h e^{-\kappa t_{k+\frac12}}\mathbb I\xi_k$, the post-momentum associated
with $(R_k,R_{k+1})$ and its evolution are
\begin{equation*}
\label{eq:rigid-momentum-equations}
\mu_{k+1}
=
e^{-\kappa t_{k+\frac12}}
\left(
\mathrm d^L\operatorname{cay}^{-1}_{F_k}
\right)^*
(I\xi_k),
\qquad
\mu_k=\operatorname{Ad}_{F_k}^*\mu_{k+1}.
\end{equation*}

Under the Euclidean identification
$\mathfrak{so}(3)^*\simeq\mathbb R^3$, the coadjoint action satisfies
$\operatorname{Ad}_F^*\mu=F\mu$. Hence,
\begin{equation*}
\label{eq:mu-update-explicit}
\mu_k=F_k\mu_{k+1},
\qquad
\mu_{k+1}=F_k^{-1}\mu_k.
\end{equation*}

For $F_k=\operatorname{cay}(h\widehat\xi_k)$, the identity
$F_k+I=2(I-\frac{h}{2}\widehat\xi_k)^{-1}$ and
\eqref{eq:dlog-cayley} give
\[
\mathrm d^L\operatorname{cay}^{-1}_{F_k}(\widehat\eta)
=
\left(I+\frac{h}{2}\widehat\xi_k\right)
\widehat\eta
\left(I-\frac{h}{2}\widehat\xi_k\right).
\]
Using
$[\widehat\xi,\widehat\eta]=\widehat{\xi\times\eta}$ and
$\widehat\xi\,\widehat\eta\,\widehat\xi
=-(\xi^\top\eta)\widehat\xi$, its vector representation is
$\mathrm d^L\operatorname{cay}^{-1}_{F_k}(\eta)
=
\eta
+
\frac{h}{2}\xi_k\times\eta
+
\frac{h^2}{4}\xi_k(\xi_k^\top\eta)$.
Consequently, with respect to the Euclidean inner product,
\begin{equation*}
\label{eq:inverse-cayley-differential-vector}
\mathrm d^L\operatorname{cay}^{-1}_{F_k}
=
I+\frac{h}{2}\widehat\xi_k
+\frac{h^2}{4}\xi_k\xi_k^\top,
\qquad
\left(
\mathrm d^L\operatorname{cay}^{-1}_{F_k}
\right)^*
=
I-\frac{h}{2}\widehat\xi_k
+\frac{h^2}{4}\xi_k\xi_k^\top.
\end{equation*}

Given $(R_k,\mu_k)$, the discrete velocity must be found from
\begin{equation}
\label{eq:rigid-implicit-xi}
\mu_k
=
F_k e^{-\kappa t_{k+\frac12}}
\left(
I-\frac h2\widehat\xi_k
+\frac{h^2}{4}\xi_k\xi_k^\top
\right)
\mathbb I\xi_k,
\qquad
F_k=\operatorname{cay}(h\widehat\xi_k).
\end{equation}

For numerical implementation, define the residual
$\mathcal{R}_k:\mathbb{R}^3\to\mathbb{R}^3$ by
\[
\mathcal{R}_k(\xi)
:=
e^{-\kappa t_{k+\frac12}}
\operatorname{cay}(h\widehat{\xi})
\left(
I_3-\frac{h}{2}\widehat{\xi}
+\frac{h^2}{4}\xi\xi^\top
\right)\mathbb{I}\xi-\mu_k.
\]
At each time step, the current momentum $\mu_k$
and time $t_k$ are known. The only unknown is
the discrete body velocity $\xi_k\in\mathbb{R}^3$,
which is determined by solving
$\mathcal{R}_k(\xi_k)=0$.
Once $\xi_k$ has been obtained, the Cayley
increment $F_k$ is computed and the configuration
and momentum are updated explicitly using (35).
No independent nonlinear equation for
$R_{k+1}$ or $\mu_{k+1}$ is required.

Local invertibility of the discrete Legendre
transformation guarantees a unique smooth
local solution of equation~\eqref{eq:rigid-implicit-xi} for sufficiently
small $h$. The resulting configuration and
momentum updates are
\begin{equation}
\label{eq:rigid-complete-update}
R_{k+1}=R_kF_k,
\qquad
\mu_{k+1}=F_k^{-1}\mu_k.
\end{equation}

Since $F_k\in SO(3)$, the quadratic Casimir is preserved exactly:
$
\|\mu_{k+1}\|^2=\|\mu_k\|^2$.

The same invariance gives preservation of the discrete Noether momentum. For
the left action $\Phi_Q(R)=QR$, the infinitesimal generator associated with
$\zeta\in\mathfrak{so}(3)$ is
$\zeta_{SO(3)}(R)=\widehat\zeta R$, whose body representation is
$\operatorname{Ad}_{R^{-1}}\zeta$. Therefore,
\begin{equation*}
\label{eq:rigid-noether-quantity}
J_\zeta^k
=
\left\langle
\mu_k,\operatorname{Ad}_{R_k^{-1}}\zeta
\right\rangle
=
\left\langle
\operatorname{Ad}_{R_k}^*\mu_k,\zeta
\right\rangle.
\end{equation*}
Thus the spatial momentum
$\Pi_k:=\operatorname{Ad}_{R_k}^*\mu_k=R_k\mu_k$ is preserved because
\begin{equation*}
\label{eq:Pi-preserved}
\Pi_{k+1}
=
R_{k+1}\mu_{k+1}
=
R_kF_kF_k^{-1}\mu_k
=
\Pi_k.
\end{equation*}

Unlike the Casimir and the Noether momentum, the continuous energy balance
\eqref{eq:continuous-energy-law} is not expected to be reproduced exactly by
a fixed-step variational integrator for a non-autonomous Lagrangian. We
therefore define the midpoint energy and its reference value by
\begin{equation*}
\label{eq:discrete-energy-reference}
E_d^{k+\frac12}
:=
\frac12e^{-\kappa t_{k+\frac12}}
\left\langle\xi_k,\mathbb I\xi_k\right\rangle,
\qquad
E_{\mathrm{ref}}^{k+\frac12}
:=
E_0e^{\kappa t_{k+\frac12}},
\end{equation*}
where for $t_0=0$, $E_0=\frac12\langle\xi(0),\mathbb I\xi(0)\rangle$. The energy approximation is assessed by comparison with the continuous
reference balance, while the convergence order of the complete phase-space
method is verified separately using final-time configuration and momentum
errors. The implementation is summarized in Table~\ref{tab:so3-lgvi}.

\begin{table}[H]
\centering
\begin{tabular}{p{0.94\linewidth}}
\hline
\textbf{Input:}
$t_0\in\mathbb R$, $R_0\in\mathrm{SO}(3)$,
$\mu_0\in\mathbb R^3$, $h>0$, $N\in\mathbb N$,
$\kappa\in\mathbb R$, and
$\mathbb I=\mathbb I^\top>0$.\\
\hline
\textbf{Initialization:}
$\xi_{\mathrm{guess}}
=e^{\kappa t_0}\mathbb I^{-1}\mu_0$.\\
\hline
\textbf{For} $k=0,\dots,N-1$: \\
\quad 1. Set $t_{k+\frac12}=t_0+(k+\frac12)h$.\\
\quad 2. Solve $\mathcal R_k(\xi_k)=0$ for $\xi_k$,
using $\xi_{\mathrm{guess}}$ as the initial iterate,
with $\|\mathcal R_k(\xi_k)\|\leq 10^{-10}$.\\
\quad 3. Compute
$F_k=\operatorname{cay}(h\widehat{\xi}_k)$.\\
\quad 4. Update $R_{k+1}=R_kF_k$.\\
\quad 5. Update $\mu_{k+1}=F_k^{-1}\mu_k$.\\
\quad 6. Compute $E_d^{k+\frac12}$,
$\|\mu_{k+1}\|^2$ and
$\Pi_{k+1}=R_{k+1}\mu_{k+1}$.\\
\quad 7. Set $\xi_{\mathrm{guess}}=\xi_k$.\\
\hline
\end{tabular}
\caption{Cayley-based variational integrator for the
time-dependent rigid body.}
\label{tab:so3-lgvi}
\end{table}

The long-time numerical metrics are
\begin{equation*}
\label{eq:rigid-diagnostics}
\begin{aligned}
\varepsilon_E^{k+\frac12}
&:=
\frac{
E_d^{k+\frac12}
-
E_0e^{\kappa t_{k+\frac12}}
}{
E_0e^{\kappa t_{k+\frac12}}
},
&
\varepsilon_C^k
&:=
\frac{
\left|
\|\mu_k\|^2-\|\mu_0\|^2
\right|
}{
\|\mu_0\|^2
},
\\
\varepsilon_\Pi^k
&:=
\frac{
\left\|
R_k\mu_k-R_0\mu_0
\right\|
}{
\|R_0\mu_0\|
},
&
\varepsilon_R^k
&:=
\left\|
R_k^\top R_k-I
\right\|_{\mathrm F}.
\end{aligned}
\end{equation*}
These metrics assess distinct properties of
the numerical method. The signed relative
energy-balance error $\varepsilon_E^{k+\frac12}$
measures the deviation of the discrete midpoint
energy from the exact continuous balance law \eqref{eq:continuous-energy-law}.
Since the reference energy decays exponentially,
normalization by its value at the same time prevents
the decreasing energy scale from concealing
relative discrepancies. A zero value indicates
exact agreement with the reference balance, while
negative and positive values indicate
underestimation and overestimation, respectively.
This metric tests the approximation of the
prescribed energy evolution, not energy conservation.

In contrast, the quadratic Casimir
$C(\mu)=\|\mu\|^2$ and the spatial momentum
$\Pi=R\mu$ are preserved exactly by the discrete
equations in exact arithmetic. Their relative
errors $\varepsilon_C^k$ and $\varepsilon_\Pi^k$
measure deviations from the corresponding initial
values, normalized to obtain dimensionless
quantities independent of their initial magnitudes.
The former shows coadjoint orbit preservation,
whereas the latter shows the discrete Noether
conservation law.

For a single rigid body, these two invariants
are not algebraically independent: since
$R_k\in\mathrm{SO}(3)$, conservation of
$\Pi_k=R_k\mu_k$ also implies conservation
of $\|\mu_k\|^2$. Nevertheless, the two
metrics separately illustrate the
coadjoint and Noether interpretations of
the discrete momentum update.

Finally, $\varepsilon_R^k$ measures the
orthogonality constraint through the Frobenius norm.
Unlike the preceding metrics, this is an
absolute constraint residual rather than a
relative error. Since the Cayley update preserves
$SO(3)$ in exact arithmetic, nonzero
values reflect floating-point effects rather
than a failure of the geometric construction.

For the long-time experiment, we take
$\kappa=-10^{-3}$,
$\mathbb I=\operatorname{diag}(2,1,0.5)$,
$R_0=I$, and
$\xi(0)=(0.8,1.2,-0.6)^\top$.
The initial left-trivialized momentum
$\mu_0\in\mathfrak{so}(3)^*$ is obtained from the continuous Legendre
transformation evaluated at $(t_0,R_0,\xi(0))$:
\begin{equation}
\label{eq:rigid-initial-momentum}
\mu_0
=
\frac{\partial L}{\partial\xi}
\bigl(t_0,R_0,\xi(0)\bigr)
=
e^{-\kappa t_0}\mathbb I\xi(0).
\end{equation}

We use $h=5\times10^{-2}$ and integrate up to
$T=10^4$, corresponding to $2\times10^5$ time
steps. At each step, the three-dimensional
nonlinear system \eqref{eq:rigid-implicit-xi} is solved using Powell's hybrid method. The prescribed continuous initial
velocity $\xi(0)$ is used as the initial iterate
for the first step, and each subsequent step
uses the previously computed discrete velocity.
The Euclidean norm of the implicit residual is
required to be at most $10^{-10}$. We also check
consistency between the updated momentum and
the post-momentum obtained from the discrete
Legendre transformation, with tolerance
$2\times10^{-10}$. The results are shown in Figure~\ref{fig:so3_diagnostics}.

\begin{figure}[h!]
\centering
\begin{minipage}{0.32\textwidth}
    \centering
    \includegraphics[width=\textwidth]{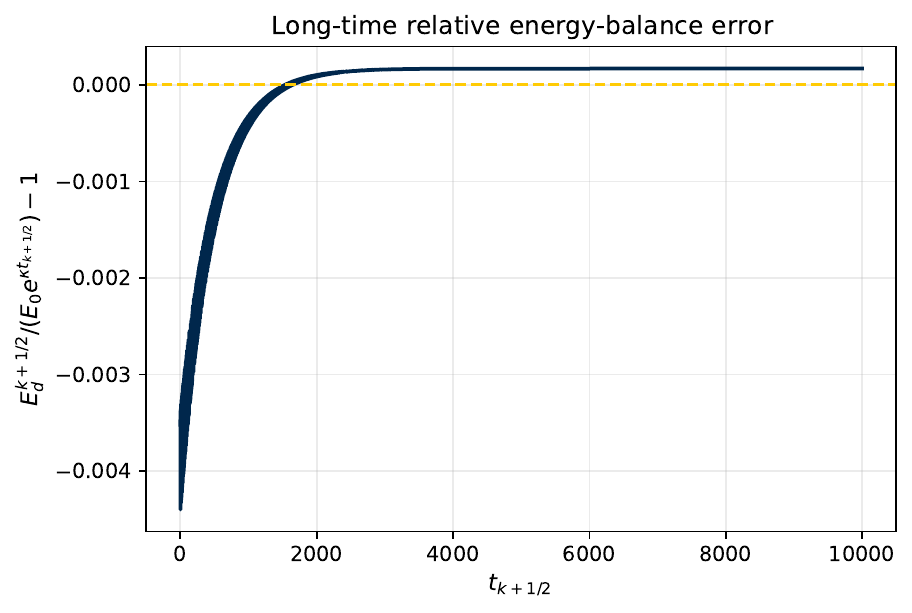}
    \small (a) Normalized energy-balance error
\end{minipage}\hfill
\begin{minipage}{0.32\textwidth}
    \centering
    \includegraphics[width=\textwidth]{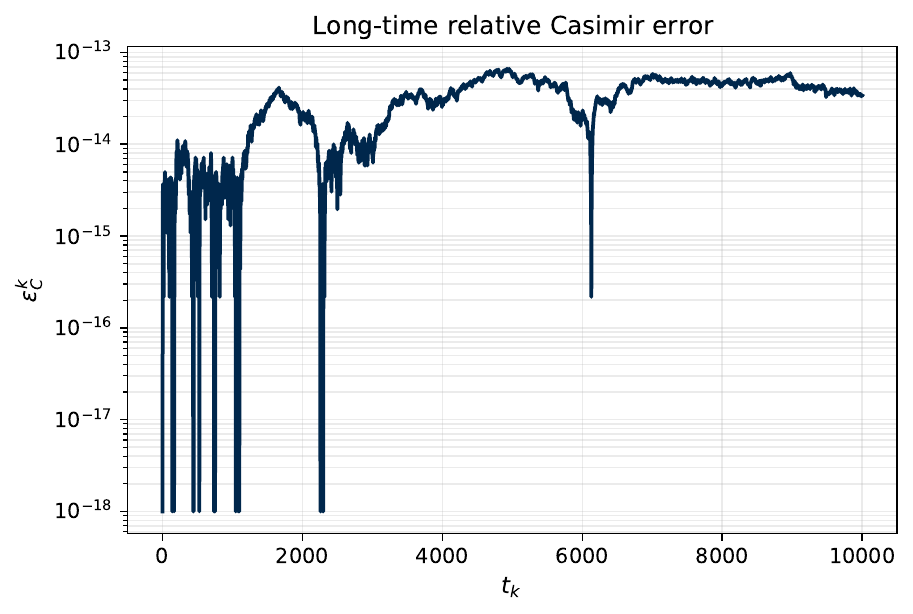}
    \small (b) Relative Casimir error
\end{minipage}\hfill
\begin{minipage}{0.32\textwidth}
    \centering
    \includegraphics[width=\textwidth]{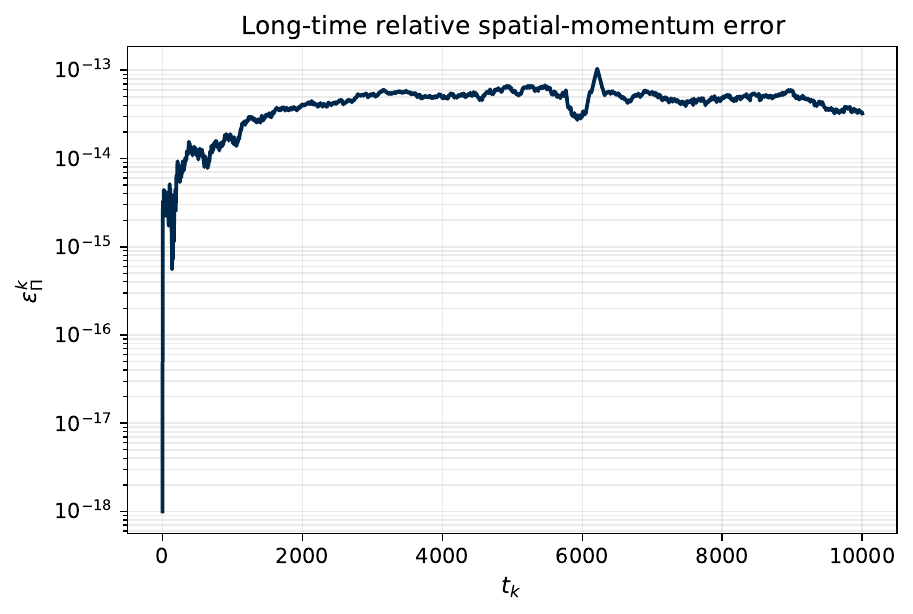}
    \small (c) Relative spatial-momentum error
\end{minipage}
\caption{Long-time numerical results for the Cayley-based
variational integrator on $SO(3)$ with
$\kappa=-10^{-3}$, $h=5\times10^{-2}$, and $T=10^4$.
(a) Normalized energy-balance error
$\varepsilon_E^{k+\frac12}$; the dashed yellow line represents exact
agreement with the continuous balance law.
(b) Relative quadratic Casimir error $\varepsilon_C^k$.
(c) Relative spatial-momentum error $\varepsilon_\Pi^k$.}
\label{fig:so3_diagnostics}
\end{figure}

Figure~\ref{fig:so3_diagnostics} shows that the signed relative
energy-balance error undergoes an initial transient and subsequently
remains close to zero from below. Over the complete integration
interval, it remains bounded, without visible
secular growth, and satisfies
$\max_k|\varepsilon_E^{k+\frac12}|
=4.419\times10^{-3}$.
The dashed horizontal line represents exact
agreement with the continuous energy balance,
not conservation of the physical energy. The quadratic Casimir and spatial momentum remain preserved, with
maximum relative errors of $6.710\times10^{-14}$ and
$1.047\times10^{-13}$, respectively. The maximum residual in the
implicit momentum equation is $1.821\times10^{-15}$. 

Although $\varepsilon_R^k$ is not displayed in
Figure~\ref{fig:so3_diagnostics}, it is evaluated at every time step
and satisfies
$\max_k\varepsilon_R^k=5.399\times10^{-13}$.
Since $\widehat\xi_k\in\mathfrak{so}(3)$, each Cayley increment
$F_k=\operatorname{cay}(h\widehat\xi_k)$ belongs to $SO(3)$, and the
exact-arithmetic update $R_{k+1}=R_kF_k$ preserves the group
constraint exactly. The reported nonzero value of
$\varepsilon_R^k$ therefore measures floating-point accumulation
rather than a truncation-induced drift away from $SO(3)$.%
\end{example}

\section{Symplecticity, Lie-Poisson reduction,
and variational error analysis}\label{sec:ham}

\subsection{Symplectic phase-space formulation
and Lie-Poisson reduction}\label{sec4.1}

We now pass from the discrete Lagrangian description on $G\times G$
to the corresponding evolution on the left-trivialized cotangent bundle
$G\times\mathfrak g^*$. This transition is not merely a change of variables.
The discrete Legendre transformations transfer the
discrete variational structure to the canonical symplectic structure of
$T^*G$, and we shall show below that a regular discrete Lagrangian
induces a symplectic phase-space map.
This provides the appropriate phase-space map for the backward error
analysis developed in Section~\ref{sec:BEA}.

We emphasize that $G\times\mathfrak g^*$ is the left-trivialized full
phase space, not the reduced phase space. Reduction by left translations
will be performed below and leads to a map on $\mathfrak g^*$.

Recall the dual left trivialization
$\lambda^*:T^*G\to G\times\mathfrak g^*$
introduced in Subsection~\ref{subsec:bridge-Q-to-G}, defined by
$\lambda^*(\alpha_g)
=(g,(T_eL_g)^*\alpha_g)$.
We denote by
$\theta_G$ and $\omega_G=-d\theta_G$ the canonical one- and two-forms
transported to $G\times\mathfrak g^*$. If a tangent vector at $(g,\mu)$
is represented as $(\eta,\nu)\in\mathfrak g\times\mathfrak g^*$ through
left trivialization, these forms are
\begin{equation*}
\label{eq:left-trivialized-canonical-forms}
\theta_G(g,\mu)(\eta,\nu)
=
\langle\mu,\eta\rangle,
\qquad
\omega_G(g,\mu)
\bigl((\eta_1,\nu_1),(\eta_2,\nu_2)\bigr)
=
\langle\nu_2,\eta_1\rangle
-
\langle\nu_1,\eta_2\rangle
+
\langle\mu,[\eta_1,\eta_2]\rangle.
\end{equation*} 
A local diffeomorphism
$\Phi:G\times\mathfrak g^*\to G\times\mathfrak g^*$
is called symplectic if $\Phi^*\omega_G=\omega_G$.

For each time interval $[t_k,t_{k+1}]$, the left-trivialized discrete
Legendre transformations associated with
$L_d^k:G\times G\to\mathbb R$ (\cite{Bates}) are defined by
\begin{equation}
\label{eq:left-trivialized-discrete-Legendre}
\begin{aligned}
\mathbb F_k^-L_d(g_k,g_{k+1})
&:=
\left(
g_k,
-
(T_e\mathcal L_{g_k})^*
D_1L_d^k(g_k,g_{k+1})
\right),\\
\mathbb F_k^+L_d(g_k,g_{k+1})
&:=
\left(
g_{k+1},
(T_e\mathcal L_{g_{k+1}})^*
D_2L_d^k(g_k,g_{k+1})
\right).
\end{aligned}
\end{equation}
According to the notation of Subsection~\ref{subsec:retractions}, these
maps satisfy
$\mathbb F_k^-L_d(g_k,g_{k+1})=(g_k,\mu_k)$ and
$\mathbb F_k^+L_d(g_k,g_{k+1})=(g_{k+1},\mu_{k+1})$ along a solution.

The discrete Euler-Lagrange equations are precisely the
momentum-matching condition
\begin{equation*}
\label{eq:Legendre-momentum-matching}
\mathbb F_{k-1}^+L_d(g_{k-1},g_k)
=
\mathbb F_k^-L_d(g_k,g_{k+1}).
\end{equation*}
Thus, a discrete trajectory can equivalently be represented by the
configuration pairs $(g_k,g_{k+1})$ or by the phase-space variables
$(g_k,\mu_k)$.

\begin{definition}
\label{def:regular-time-dependent-Ld}
The family $L_d^k$ is \textit{regular} if both discrete Legendre
transformations in \eqref{eq:left-trivialized-discrete-Legendre} are local
diffeomorphisms for every $k$.
\end{definition}

Under this regularity assumption, the discrete Euler-Lagrange evolution
on consecutive configuration pairs and the \textit{induced phase-space maps} are
defined by
\begin{equation*}
\label{eq:edge-and-phase-flows}
\Psi_d^k
:=
(\mathbb F_k^-L_d)^{-1}
\circ\mathbb F_{k-1}^+L_d,
\qquad
\widetilde\Psi_d^k
:=
\mathbb F_k^+L_d
\circ
(\mathbb F_k^-L_d)^{-1}.
\end{equation*}
Here $\Psi_d^k(g_{k-1},g_k)=(g_k,g_{k+1})$, whereas
$\widetilde\Psi_d^k(g_k,\mu_k)=(g_{k+1},\mu_{k+1})$. Since the
Lagrangian depends explicitly on the time index, the evolution is
generated by a sequence of maps rather than by iterates of a single
autonomous map.

The relation between the two representations is summarized by the
following commutative diagram.

\begin{figure}[h!]
\centering
\begin{tikzcd}[column sep=large,row sep=large]
G\times G
\arrow[r,"\Psi_d^k"]
\arrow[d,"\mathbb F_{k-1}^+L_d"']
&
G\times G
\arrow[d,"\mathbb F_k^+L_d"]
\\
G\times\mathfrak g^*
\arrow[r,"\widetilde\Psi_d^k"']
&
G\times\mathfrak g^*
\end{tikzcd}
\caption{Commutative diagram relating the discrete evolution on $G\times G$ to the left-trivialized phase-space evolution
on $G\times\mathfrak g^*$. The momentum-matching condition identifies
$\mathbb F_{k-1}^+L_d(g_{k-1},g_k)$ with
$\mathbb F_k^-L_d(g_k,g_{k+1})$, while both paths in the diagram yield
$(g_{k+1},\mu_{k+1})$.}
\label{fig:discrete-Legendre-diagram}
\end{figure}
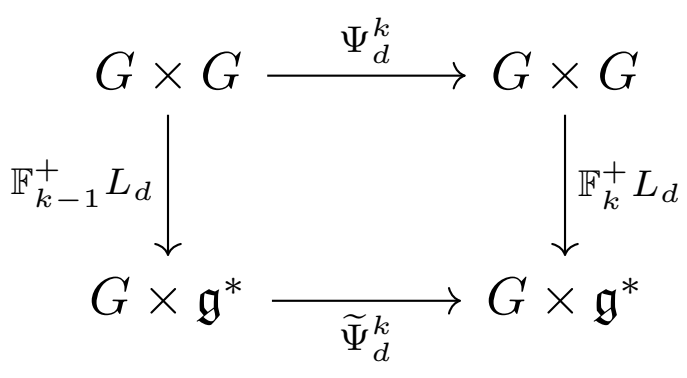

 The identity encoded
by Figure~\ref{fig:discrete-Legendre-diagram} is
\begin{equation}
\label{eq:Legendre-conjugacy}
\mathbb F_k^+L_d\circ\Psi_d^k
=
\widetilde\Psi_d^k\circ\mathbb F_{k-1}^+L_d.
\end{equation}

We next introduce the discrete Poincaré-Cartan one-forms
$\Theta_k^\pm\in\Omega^1(G\times G)$ and the associated
closed two-form
$\Omega_k\in\Omega^2(G\times G)$ by
\begin{equation}
\Theta_k^\pm := (\mathbb F_k^\pm L_d)^*\theta_G,
\qquad
\Omega_k := -d\Theta_k^+ = -d\Theta_k^- .
\label{eq:discrete-poincare-cartan}
\end{equation}
The equality between the two expressions for $\Omega_k$ follows from
$\Theta_k^+-\Theta_k^-=dL_d^k$.
Since $\omega_G=-d\theta_G$, one also has
\[
\Omega_k
=(\mathbb F_k^+L_d)^*\omega_G
=(\mathbb F_k^-L_d)^*\omega_G.
\]
When $L_d^k$ is regular, $\Omega_k$ is nondegenerate and hence defines
a symplectic form on $G\times G$. We refer to it as the discrete
symplectic form associated with the interval $[t_k,t_{k+1}]$.

We use the following standard notion. A smooth function
$S:G\times G\to\mathbb R$ is a generating function of the first kind
for a local map $\Phi:T^*G\to T^*G$ if the graph of $\Phi$ is locally
parametrized by $(g_0,g_1)\in G\times G$ through
\[
p_0=-D_1S(g_0,g_1),
\qquad
p_1=D_2S(g_0,g_1),
\] that is,
$\Phi(g_0,p_0)=(g_1,p_1)$.
Under the left trivialization
$T^*G\simeq G\times\mathfrak g^*$, the same relations define the
corresponding local phase-space map on $G\times\mathfrak g^*$.

\begin{theorem}
\label{thm:discrete-phase-symplecticity}
Let $L_d^k:G\times G\to\mathbb R$ be regular. Then $L_d^k$ is a
generating function of the first kind for the phase-space map
$\widetilde\Psi_d^k$. In particular,
\begin{equation}
\label{eq:phase-and-edge-symplecticity}
(\widetilde\Psi_d^k)^*\omega_G=\omega_G,
\qquad
(\Psi_d^k)^*\Omega_k=\Omega_{k-1}.
\end{equation}
Therefore, although the family of maps is time-dependent,
every one-step phase-space map
$\widetilde\Psi_d^k$ is symplectic on $G\times\mathfrak g^*$,
while the configuration-pair evolution transports the discrete
symplectic form associated with one time interval to that of the
preceding interval.
\end{theorem}

\begin{proof}
Let $(g_k,\mu_k)\in G\times\mathfrak g^*$ and determine $g_{k+1}$
implicitly from
$\mathbb F_k^-L_d(g_k,g_{k+1})=(g_k,\mu_k)$. By
\eqref{eq:left-trivialized-discrete-Legendre}, the corresponding output
satisfies
$\mathbb F_k^+L_d(g_k,g_{k+1})=(g_{k+1},\mu_{k+1})$. Equivalently, the
initial and final covectors are determined by
$-D_1L_d^k(g_k,g_{k+1})$ and $D_2L_d^k(g_k,g_{k+1})$, respectively.
Hence $L_d^k$ generates $\widetilde\Psi_d^k$ as a generating function of
the first kind.

Using \eqref{eq:discrete-poincare-cartan} and
$\widetilde\Psi_d^k=\mathbb F_k^+L_d\circ(\mathbb F_k^-L_d)^{-1}$, we
obtain
\begin{equation*}
\label{eq:phase-symplecticity-proof}
\begin{aligned}
(\widetilde\Psi_d^k)^*\omega_G
&=
\bigl((\mathbb F_k^-L_d)^{-1}\bigr)^*
(\mathbb F_k^+L_d)^*\omega_G=
\bigl((\mathbb F_k^-L_d)^{-1}\bigr)^*\Omega_k
=
\omega_G.
\end{aligned}
\end{equation*}
Similarly, the conjugacy relation \eqref{eq:Legendre-conjugacy} gives
\begin{equation*}
\label{eq:edge-symplecticity-proof}
\begin{aligned}
(\Psi_d^k)^*\Omega_k
&=
(\Psi_d^k)^*(\mathbb F_k^+L_d)^*\omega_G=
(\mathbb F_{k-1}^+L_d)^*
(\widetilde\Psi_d^k)^*\omega_G
=
\Omega_{k-1}.
\end{aligned}
\end{equation*}
\end{proof}

The second identity in \eqref{eq:phase-and-edge-symplecticity} should
not be replaced by $(\Psi_d^k)^*\Omega_k=\Omega_k$ unless the discrete
Lagrangian is autonomous. For a time-dependent family,
$\Omega_{k-1}$ and $\Omega_k$ need not coincide. By contrast, after
passing to $G\times\mathfrak g^*$ through the discrete Legendre
transformations, every $\widetilde\Psi_d^k$ preserves the same canonical
form $\omega_G$.

This distinction is essential for Section~\ref{sec:BEA}: the
time-dependent composition
\begin{equation*}
\widetilde\Psi_d^{N-1}
\circ\cdots\circ
\widetilde\Psi_d^1
\circ
\widetilde\Psi_d^0
\end{equation*}
is symplectic on $G\times\mathfrak g^*$
because each one-step map
$\widetilde\Psi_d^k$ is symplectic with respect to the same form
$\omega_G$.
The extended phase-space formulation used in the next section
incorporates time as a dynamical variable and yields
a single exactly symplectic map on the extended phase space.

\medskip

We next use this symplectic phase-space formulation to
describe the reduced dynamics in the presence of left invariance. Assume that every $L_d^k$ is invariant under the diagonal left action,
\[
L_d^k(hg_0,hg_1)=L_d^k(g_0,g_1),
\qquad h\in G.
\]
The cotangent lift of the left action of $G$ on itself induces an action
on $T^*G$. Under the left trivialization
$T^*G\simeq G\times\mathfrak g^*$ introduced above, the corresponding
quotient is naturally identified with $\mathfrak g^*$:
$T^*G/G \simeq \mathfrak g^*$. Equivalently, the quotient projection is represented in left-trivialized
variables by $\pi:G\times\mathfrak g^*\rightarrow\mathfrak g^*$, $\pi(g,\mu)=\mu$.

The quotient $\mathfrak g^*$ carries the Lie-Poisson structure
associated with the coadjoint convention fixed in Section~\ref{Sec2-Lie}.
For $f,k\in C^\infty(\mathfrak g^*)$, define their functional
derivatives by
\[
\left\langle \delta\mu,
\frac{\delta f}{\delta\mu}(\mu)\right\rangle
=
\left.\frac{d}{ds}\right|_{s=0}
f(\mu+s\delta\mu),
\qquad
\delta\mu\in\mathfrak g^*.
\]
With our convention, the Lie-Poisson bracket is
\begin{equation*}
\label{eq:Lie-Poisson-bracket}
\{f,k\}_{\mathrm{LP}}(\mu)
=
-\left\langle
\mu,
\left[
\frac{\delta f}{\delta\mu},
\frac{\delta k}{\delta\mu}
\right]
\right\rangle .
\end{equation*}

A smooth map
$\psi:\mathfrak g^*\to\mathfrak g^*$ is called a
\textit{Poisson map} if
\[
\{f\circ\psi,k\circ\psi\}_{\mathrm{LP}}
=
\{f,k\}_{\mathrm{LP}}\circ\psi
\]
for all $f,k\in C^\infty(\mathfrak g^*)$. %

For $\mu\in\mathfrak g^*$, its \textit{coadjoint orbit} is
\[
\mathcal O_\mu
:=
\left\{
\operatorname{Ad}_g^*\mu
\;:\;
g\in G
\right\}.
\]
The connected components of the coadjoint orbits are the symplectic
leaves of the Lie-Poisson structure. In particular, when $G$ is
connected, each coadjoint orbit is connected and hence is itself a
symplectic leaf.

A function $C\in C^\infty(\mathfrak g^*)$ is a \textit{Casimir} if
\[
\{C,f\}_{\mathrm{LP}}=0
\qquad
\text{for every }f\in C^\infty(\mathfrak g^*).
\]
Equivalently, a Casimir is constant on every symplectic leaf and,
when $G$ is connected, on every coadjoint orbit. By \textit{coadjoint transport} we mean a discrete momentum
update of the form
\[
\mu_{k+1}
=
\operatorname{Ad}_{a_k}^*\mu_k,
\qquad a_k\in G.
\]
Such an update moves the momentum within its coadjoint orbit.
If $G$ is connected, each coadjoint orbit is a symplectic leaf and
the update therefore preserves every Casimir. Notice, however, that orbit
preservation alone does not imply that the resulting map
$\mu_k\mapsto\mu_{k+1}$ is Poisson when $a_k$ depends on $\mu_k$.
The Poisson property must follow from the reduced geometric origin
of the map.

\begin{theorem}
\label{thm:discrete-Lie-Poisson-reduction}

Suppose that $L_d^k$ is regular and left-invariant
for every $k$. On any $G$-invariant regularity
domain where the phase-space map
$\widetilde\Psi_d^k$ is well defined,
this map descends through the quotient projection
$\pi:G\times\mathfrak g^*\to\mathfrak g^*$
to a local Poisson map
\[
\psi_d^k:U_k\subset\mathfrak g^*
\longrightarrow\mathfrak g^*,
\qquad
\pi\circ\widetilde\Psi_d^k
=\psi_d^k\circ\pi,
\]
where $U_k$ is an open subset of $\mathfrak g^*$
and the identity holds on the corresponding
regularity domain.

For the Cayley parametrization introduced in Subsection~\ref{subsec:retractions}, this reduced
map is determined implicitly by
\[
\mu_{k+1}
=
\psi_d^k(\mu_k)
=
\Ad^*_{F_k^{-1}}\mu_k,
\qquad
F_k=\cay(h\xi_k).
\]
Here, $\xi_k$ is determined implicitly by
combining the discrete post-Legendre
relation~\eqref{eq:force-and-postmomentum} with the momentum-matching
condition~\eqref{eq:discrete-Euler-Poincare}.
Consequently, $\psi_d^k$ preserves coadjoint orbits. If, in addition,
$G$ is connected, it preserves every Casimir of the Lie-Poisson structure.
\end{theorem}

\begin{proof}
Left invariance of $L_d^k$ implies that both
discrete Legendre transformations are equivariant
with respect to the diagonal left action on
$G\times G$ and the corresponding cotangent-lifted
action. Consequently, the phase-space map
$\widetilde\Psi_d^k$ is equivariant on every
$G$-invariant regularity domain. Under left
trivialization, this action takes the form
$a\cdot(g,\mu)=(ag,\mu)$.
Therefore, $\widetilde\Psi_d^k$ induces a
well-defined local map $\psi_d^k$ on the quotient
$G\times\mathfrak g^*/G\simeq\mathfrak g^*$,
satisfying
\[
\pi\circ\widetilde\Psi_d^k
=\psi_d^k\circ\pi.
\]

We now establish the Poisson property.
Let $\{\cdot,\cdot\}_G$ denote the canonical
Poisson bracket on $G\times\mathfrak g^*$.
The quotient projection $\pi$ is a Poisson map
with respect to the Lie-Poisson bracket
introduced above, namely
\[
\{f_1\circ\pi,f_2\circ\pi\}_G
=
\{f_1,f_2\}_{LP}\circ\pi.
\]
This identity follows directly from the
canonical symplectic form $\omega_G$ and
the Lie-Poisson bracket defined above.

Since $\widetilde\Psi_d^k$ is symplectic by
Theorem~\ref{thm:discrete-phase-symplecticity}, for any smooth functions
$f_1,f_2$ defined on a neighborhood of the
image of $\psi_d^k$, we obtain
\begin{align*}
\{f_1\circ\psi_d^k,
   f_2\circ\psi_d^k\}_{LP}\circ\pi
&=
\{f_1\circ\pi\circ\widetilde\Psi_d^k,
  f_2\circ\pi\circ\widetilde\Psi_d^k\}_G=
\{f_1\circ\pi,f_2\circ\pi\}_G
\circ\widetilde\Psi_d^k=
\{f_1,f_2\}_{LP}
\circ\pi\circ\widetilde\Psi_d^k\\
&\quad=
\{f_1,f_2\}_{LP}
\circ\psi_d^k\circ\pi.
\end{align*}
Since $\pi$ is surjective onto the reduced
regularity domain $U_k$, it follows that
\[
\{f_1\circ\psi_d^k,
  f_2\circ\psi_d^k\}_{LP}
=
\{f_1,f_2\}_{LP}\circ\psi_d^k.
\]
Thus $\psi_d^k$ is a Poisson map.

For the Cayley parametrization, the momentum
update follows directly from \eqref{eq:discrete-Euler-Poincare}:
\[
\mu_{k+1}
=\operatorname{Ad}_{F_k^{-1}}^*\mu_k,
\qquad
F_k=\operatorname{cay}(h\xi_k).
\]
Hence $\mu_{k+1}$ belongs to the coadjoint
orbit through $\mu_k$. If $G$ is connected,
these orbits are the symplectic leaves of the
Lie-Poisson structure. Therefore, every
Casimir function is preserved along the
discrete evolution.
\end{proof}

\subsection{Variational error analysis and second-order accuracy}\label{sec4.2}

We now study the approximation order of the discrete phase-space map. A pointwise estimate for the exact discrete Lagrangian does not necessarily give the same order for the induced variational integrator. Indeed, the discrete Legendre transformations involve derivatives with respect to both endpoints. These derivatives may have a lower order of accuracy than the discrete action itself. This difficulty is at the heart of the variational error analysis of Patrick and Cuell \cite{PatrickCuell2009}.

To illustrate the problem, consider a free particle on $Q=\mathbb R$ with $L(q,\dot q)=\frac12\dot q^2$. Its exact discrete Lagrangian is $L_{d,E}(q_0,q_1;h)=(q_1-q_0)^2/(2h)$. Setting $\Delta=q_1-q_0$, consider
\[
L_d(q_0,q_1;h)=\frac{\Delta^2}{2h}
+\varepsilon h^3\sin\!\left(\frac{\Delta}{h^2}\right),
\qquad 0<|\varepsilon|<1.
\]
The action defect satisfies $|L_d-L_{d,E}|\leq |\varepsilon|h^3$. Nevertheless, the discrete Legendre transformations give
\[
p_0=-D_1L_d=\frac{\Delta}{h}
+\varepsilon h\cos\!\left(\frac{\Delta}{h^2}\right),
\qquad p_1=D_2L_d=p_0,
\]
and hence
\[
q_1-(q_0+hp_0)=-\varepsilon h^2
\cos\!\left(\frac{\Delta}{h^2}\right).
\]
This error is of order $h^2$. For example, if $p_0=0$ and $\Delta=h^2s$, then $s+\varepsilon\cos s=0$. This equation has a unique nonzero solution because its derivative is strictly positive and its value at $s=0$ is nonzero. Moreover,
\[
\frac{\partial p_0}{\partial\Delta}
=\frac{1-\varepsilon\sin(\Delta/h^2)}{h}
\geq\frac{1-|\varepsilon|}{h}>0.
\]
Thus the discrete Lagrangian is regular, but its $O(h^3)$
pointwise action defect alone does not suffice to establish
second-order accuracy of the associated variational integrator.

The problem is the lack of differentiable control of the normalized defect. Along an exact solution $q_1=q_0+hv$ we have
\[
h^{-3}(L_d-L_{d,E})=\varepsilon\sin(v/h),
\qquad
\frac{\partial}{\partial v}\bigl[h^{-3}(L_d-L_{d,E})\bigr]
=\frac{\varepsilon}{h}\cos(v/h).
\]
The second expression is not uniformly bounded as $h\to0$. Therefore, an estimate of the action values alone is insufficient. We need an approximation property that also controls derivatives with respect to the initial data.

Patrick and Cuell addressed this issue through a desingularized variational principle and a past-future symmetry that accounts for the cancellation of lower-order terms \cite{PatrickCuell2009}. We use the same cancellation principle but work directly with the exact endpoint map and the discrete Legendre transformations. This gives a proof in the left-trivialized phase space and accommodates explicit time dependence. We then apply it to the midpoint-Cayley integrator.

Recall the left trivialization $TG\simeq G\times\mathfrak g$ from Subsection~2.1. Let $\varphi^L_{t,h}$ be the local Euler-Lagrange evolution from time $t$ to $t+h$. In these variables write
\[
\varphi^L_{t,h}(g,\xi)
=\bigl(\Gamma_{t,h}(g,\xi),\Xi_{t,h}(g,\xi)\bigr),
\]
where $\Gamma_{t,h}$ and $\Xi_{t,h}$ denote the final configuration and body velocity. The exact endpoint map is $\mathcal E_{t,h}(g,\xi):=(g,\Gamma_{t,h}(g,\xi))$. Let $L_{d,E}(t,\cdot,\cdot;h)$ be the exact discrete Lagrangian over $[t,t+h]$ and let $L_d(t,\cdot,\cdot;h)$ be its approximation, with $L_d^k=L_d(t_k,\cdot,\cdot;h)$. We define the action defect and the \textit{pulled-back action defect} by
\[
\begin{aligned}
\delta L_h(t,g_0,g_1)
&:=L_d(t,g_0,g_1;h)-L_{d,E}(t,g_0,g_1;h),\\
\Delta_d(t,h,g,\xi)
&:=\delta L_h\bigl(t,\mathcal E_{t,h}(g,\xi)\bigr).
\end{aligned}
\]
The first function is defined on pairs of configurations. The second is its pullback to $G\times\mathfrak g$ along exact solutions. This distinction is essential because the regularity assumption is imposed on $\Delta_d$ while the discrete dynamics depends on the endpoint derivatives of $\delta L_h$.

A pointwise estimate $\Delta_d=O(h^{r+1})$ controls only the values of the action defect. We use the following stronger condition.

\begin{definition}
We say that $L_d$ has smooth variational contact of order $r$ with $L_{d,E}$ if
\begin{equation}
\Delta_d(t,h,g,\xi)=h^{r+1}\rho(t,h,g,\xi),
\label{eq:smooth-variational-contact}
\end{equation}
where $\rho$ and its first derivatives with respect to $(g,\xi)$ extend continuously to $h=0$ and remain uniformly bounded for $t$ in compact intervals and $(g,\xi)$ in compact subsets of $G\times\mathfrak g$.
\end{definition}

Throughout this subsection, group derivatives
of scalar functions are left-trivialized
as in \eqref{eq:DgL}. Derivatives with respect to $\xi$ are taken
in the vector space $\mathfrak g$.
In particular, smooth variational contact
implies
$\mathbf{D}_g\Delta_d=O(h^{r+1})$
and $D_\xi\Delta_d=O(h^{r+1})$
uniformly on compact subsets.
These bounds will allow us to recover the
correct order of the phase-space map even
though each endpoint derivative of
$\delta L_h$ may only be $O(h^r)$.

\begin{theorem}\label{thm:variational-order}
Let $L:\mathbb R\times G\times\mathfrak g\to\mathbb R$ be a sufficiently smooth regular time-dependent Lagrangian. Let $L_d(t,\cdot,\cdot;h)$ be a sufficiently smooth family of regular discrete Lagrangians for all sufficiently small $h>0$. Assume that the exact Euler-Lagrange flow, the exact discrete Lagrangian and the exact and approximate discrete Legendre transformations are defined on common local domains containing the trajectory segments issued from a compact subset of $G\times\mathfrak g$.

If $L_d$ has smooth variational contact of order $r\geq1$ in the sense of \eqref{eq:smooth-variational-contact}, then its left-trivialized phase-space map $\widetilde\Psi_{t,h}$ approximates the exact Hamiltonian evolution $\Phi_{t,h}$ on $G\times\mathfrak g^*$ with local error $O(h^{r+1})$, uniformly on compact subsets of the corresponding regularity domains. Hence $\widetilde\Psi_{t,h}$ is a one-step method of order $r$.
\end{theorem}

\begin{proof}
We work in $G\times\mathfrak g$ and use left trivialization for all tangent and cotangent variables. The argument is local and all estimates are uniform on the compact sets under consideration.

Fix $(t,g,\xi)$ and set $g_1=\Gamma_{t,h}(g,\xi)$. Smooth dependence on the initial data and the reconstruction equation (2) give $\log(g^{-1}\Gamma_{t,h}(g,\xi))=h\xi+O(h^2)$. Introduce the left-trivialized derivatives of the final configuration,
\[
\begin{aligned}
A_{t,h}&:=T_{g_1}L_{g_1^{-1}}\circ D_\xi\Gamma_{t,h},\\
B_{t,h}&:=T_{g_1}L_{g_1^{-1}}\circ D_g\Gamma_{t,h}\circ T_eL_g.
\end{aligned}
\]
They are linear maps from $\mathfrak g$ to itself. Taylor expansion gives $A_{t,h}=hI+O(h^2)$ and $B_{t,h}=I+O(h)$. Thus $A_{t,h}$ is invertible for sufficiently small $h>0$ and $A_{t,h}^{-1}=h^{-1}(I+O(h))$.

Let $\delta L_h$ be the action defect defined above. At the exact endpoint pair $(g,g_1)$ set
\[
a_1:=(T_eL_g)^*D_1\delta L_h,
\qquad
a_2:=(T_eL_{g_1})^*D_2\delta L_h.
\]
Both belong to $\mathfrak g^*$. Applying the chain rule to the pulled-back
action defect gives
\[
D_\xi\Delta_d=A_{t,h}^*a_2,
\qquad
\mathbf{D}_g\Delta_d=a_1+B_{t,h}^*a_2.
\]
By \eqref{eq:smooth-variational-contact}, both derivatives on the left are $O(h^{r+1})$. Invertibility of $A_{t,h}$ gives $a_2=O(h^r)$. Consequently,
\[
\begin{aligned}
a_1+a_2
&=\mathbf{D}_g\Delta_d+(I-B_{t,h}^*)a_2\\
&=O(h^{r+1})+O(h)O(h^r)
=O(h^{r+1}).
\end{aligned}
\]
This is the cancellation of the leading endpoint defects. It is the step that restores the order lost when the endpoint derivatives are estimated separately.

We next compare the exact and discrete phase-space evolutions. Define the initial and final continuous left-trivialized momenta along the exact solution by
\[
\begin{aligned}
M_0(t,g,\xi)&:=\frac{\partial L}{\partial\xi}(t,g,\xi),\\
M_1(t,h,g,\xi)&:=\frac{\partial L}{\partial\xi}
\bigl(t+h,\Gamma_{t,h}(g,\xi),\Xi_{t,h}(g,\xi)\bigr).
\end{aligned}
\]
Smoothness of the flow implies $M_1-M_0=O(h)$ and $D_\xi(M_1-M_0)=O(h)$. The exact endpoint identities, expressed in left-trivialized momenta, are
\[
\begin{aligned}
-(T_eL_g)^*D_1L_{d,E}(t,g,g_1;h)&=M_0(t,g,\xi),\\
(T_eL_{g_1})^*D_2L_{d,E}(t,g,g_1;h)&=M_1(t,h,g,\xi).
\end{aligned}
\]

Fix an initial phase-space point $(g,\mu)$ and let $\xi$ be determined by the continuous Legendre relation $\mu=M_0(t,g,\xi)$. Let $\bar\xi$ denote the initial body velocity whose exact endpoint agrees with the endpoint selected by the discrete negative Legendre transformation. Thus
\[
g_1^d=\Gamma_{t,h}(g,\bar\xi),\qquad
\mu=-(T_eL_g)^*D_1L_d(t,g,g_1^d;h).
\]
The exact endpoint identity gives $M_0(t,g,\bar\xi)-M_0(t,g,\xi)=a_1(t,h,g,\bar\xi)=O(h^r)$. Since $L$ is regular, the continuous velocity Legendre transformation is locally invertible. It follows that $\bar\xi-\xi=O(h^r)$. As $D_\xi\Gamma_{t,h}=O(h)$, the final configurations satisfy
\[
\chi(g_1^d)-\chi(g_1^E)=O(h^{r+1}),
\qquad g_1^E:=\Gamma_{t,h}(g,\xi),
\]
in any local chart $\chi$ containing both endpoints.

The final discrete and exact body momenta are
\[
\mu_1^d=M_1(t,h,g,\bar\xi)+a_2(t,h,g,\bar\xi),
\qquad \mu_1^E=M_1(t,h,g,\xi).
\]
Using $D_\xi(M_1-M_0)=O(h)$ together with $\bar\xi-\xi=O(h^r)$, we obtain
\[
\begin{aligned}
\mu_1^d-\mu_1^E
&=M_0(t,g,\bar\xi)-M_0(t,g,\xi)+a_2+O(h^{r+1})\\
&=a_1+a_2+O(h^{r+1})=O(h^{r+1}),
\end{aligned}
\]
where $a_1$ and $a_2$ are evaluated at the endpoint pair determined by $\bar\xi$. Both phase-space components therefore have local error $O(h^{r+1})$. This estimate is invariant under smooth changes of local coordinates and proves the result on $G\times\mathfrak g^*$.
\end{proof}

\begin{corollary}\label{cor:midpoint-cayley-order}
Let $G$ be a matrix Lie group on which the Cayley transform defines a local retraction, and let $L:\mathbb R\times G\times\mathfrak g\to\mathbb R$ be a sufficiently smooth regular time-dependent Lagrangian. Assume that the exact trajectories issued from a compact subset of $G\times\mathfrak g$ remain in common local domains where the group logarithm, the Cayley map and its inverse, and the exact discrete Lagrangian are well defined. Assume also that the midpoint-Cayley discrete Legendre transformations are locally invertible for sufficiently small $h>0$.

Then the midpoint-Cayley discrete Lagrangian \eqref{eq:Ld-midpoint-retraction} has smooth variational contact of order two with the exact discrete Lagrangian. Its phase-space map $\widetilde\Psi_d^k$ is therefore symplectic and has order two on $G\times\mathfrak g^*$. If $L$ is left-invariant, the induced reduced map $\psi_d^k$ is a second-order Poisson map on $\mathfrak g^*$.
\end{corollary}

\begin{proof}
Fix $(g,\xi)\in G\times\mathfrak g$ at time $t$ and let $(g(s),\xi(s))$ be the corresponding exact solution, with $g(t)=g$ and $\xi(t)=\xi$. Write $t_m=t+h/2$ and $F_h=g^{-1}g(t+h)$. For sufficiently small $h$, the increment $F_h$ lies in the common domain of the logarithm and the inverse Cayley map.

The centered Magnus expansion gives
\[
\log F_h=h\xi(t_m)+h^3A(t,h,g,\xi),
\]
where $A$ extends smoothly to $h=0$ and depends smoothly on the initial data. Since $\operatorname{cay}^{-1}(\exp\eta)=\eta+O(\|\eta\|^3)$ with a smooth normalized remainder near the origin, the discrete velocity evaluated on exact endpoints satisfies
\[
\xi_h:=\frac1h\operatorname{cay}^{-1}(F_h)
=\xi(t_m)+h^2B(t,h,g,\xi),
\]
for a smooth function $B$.

Let $g_m=g(t_m)$ be the exact midpoint configuration and let $g_m^d:=g\operatorname{cay}(\frac h2\xi_h)$ be the midpoint configuration used by the discrete Lagrangian. The half-step reconstruction expansion and the preceding estimate give
\[
\log\bigl(g_m^{-1}g_m^d\bigr)
=h^2C(t,h,g,\xi),
\]
where $C$ is smooth on the compact family under consideration. Hence smoothness of $L$ implies
\[
L(t_m,g_m^d,\xi_h)
=L(t_m,g_m,\xi(t_m))+h^2D(t,h,g,\xi),
\]
with $D$ smooth. Evaluating \eqref{eq:Ld-midpoint-retraction} on the exact endpoints therefore gives
\[
L_d(t,g,g(t+h);h)
=hL(t_m,g_m,\xi(t_m))+h^3D(t,h,g,\xi).
\]
On the other hand, midpoint quadrature applied to the exact action yields
\[
\begin{aligned}
L_{d,E}(t,g,g(t+h);h)
&=\int_t^{t+h}L(s,g(s),\xi(s))\,ds\\
&=hL(t_m,g_m,\xi(t_m))+h^3E(t,h,g,\xi),
\end{aligned}
\]
where $E$ has the same smooth dependence. Subtracting the two expansions gives
\[
\Delta_d(t,h,g,\xi)
=h^3\bigl(D(t,h,g,\xi)-E(t,h,g,\xi)\bigr).
\]
The normalized defect and its first derivatives with respect to $(g,\xi)$ are therefore uniformly bounded on compact subsets. This is smooth variational contact of order two. Theorem~\ref{thm:variational-order} gives local error $O(h^3)$ for the phase-space map, and Theorem~\ref{thm:discrete-phase-symplecticity} gives symplecticity.

If $L$ is left-invariant, the midpoint-Cayley construction is left-invariant because the relative increment $g_0^{-1}g_1$ is unchanged by simultaneous left translation of both endpoints. Theorem~\ref{thm:discrete-Lie-Poisson-reduction} then implies that the reduced map is Poisson. Its local error remains $O(h^3)$ after reduction, so the reduced method is a second-order Lie-Poisson integrator.
\end{proof}

\begin{example}
\label{ex:numerical-variational-order}

We now verify numerically the second-order
convergence of the midpoint-Cayley variational
integrator established in Corollary~\ref{cor:midpoint-cayley-order}, using
the time-dependent rigid body of Example~\ref{ejemplo-discretizado}.
The corollary gives a local error of order
$O(h^3)$, which yields a global error of
order $O(h^2)$ on fixed time intervals under
the standard smoothness and stability
assumptions for one-step methods. Accordingly,
we fix the final time and successively refine
the step size, measuring the errors in both
components of the left-trivialized phase space.

For this experiment, we use
$\kappa=-0.4$,
$\mathbb I=\operatorname{diag}(2,1,0.5)$,
$R_0=I$, and
$\xi(0)=(0.8,1.2,-0.6)^\top$.
The initial left-trivialized momentum is determined by
\eqref{eq:rigid-initial-momentum}. The equations are integrated up to
$T=10$ with step sizes
$h_j\in\{0.08,0.04,0.02,0.01,0.005\}$.

The values of $\kappa$ are selected for
different numerical purposes, rather than
for a comparison of stability regimes. In the convergence experiment, $\kappa=-0.4$ makes the
explicit time dependence appreciable over the fixed interval $T=10$.
Using the same value in the long-time experiment with $T=10^4$ would
give the reference factor
$e^{\kappa T}=e^{-4000}$, which is below floating-point range and would
make the normalized energy-balance experiment numerically meaningless.
The value $\kappa=-10^{-3}$ is therefore used in the long-time
experiment so that the prescribed energy balance remains numerically
resolvable throughout the integration interval. The exact discrete
preservation of the Casimir and spatial momentum does not depend on
this choice of $\kappa$; it follows from the coadjoint update and the
discrete Noether theorem, respectively.

A reference solution of the continuous equations is computed with the
eighth-order adaptive DOP853 method, using relative and absolute
tolerances $2\times10^{-13}$ and $2\times10^{-14}$, respectively.
We denote its final configuration and momentum by
$R_{\mathrm{ref}}(T)$ and $\mu_{\mathrm{ref}}(T)$.

The final configuration error is measured by the geodesic distance on
$SO(3)$, while the momentum error is measured with the Euclidean norm.
The corresponding errors and observed orders are
\begin{equation*}
\begin{aligned}
e_R(h_j)
&:=
\left\|
\left(
\log\left(
R_{\mathrm{ref}}(T)^\top R_{h_j}(T)
\right)
\right)^\vee
\right\|,
&
e_\mu(h_j)
&:=
\left\|
\mu_{h_j}(T)-\mu_{\mathrm{ref}}(T)
\right\|,
\\
p_R(h_j)
&:=
\log_2
\left(
\frac{e_R(2h_j)}{e_R(h_j)}
\right),
&
p_\mu(h_j)
&:=
\log_2
\left(
\frac{e_\mu(2h_j)}{e_\mu(h_j)}
\right).
\end{aligned}
\end{equation*}
Here $\log:SO(3)\to\mathfrak{so}(3)$ denotes the principal matrix
logarithm in a neighborhood of the identity. The observed orders are
reported for every step size $h_j$ for which the computation with the
coarser step $2h_j$ is available.

For this left-invariant rigid body, the
configuration and momentum errors are related
through conservation of the spatial momentum,
$R_h(T)\mu_h(T)=R_{\mathrm{ref}}(T)
\mu_{\mathrm{ref}}(T)=\Pi_0$.
Nevertheless, we report both errors to
assess convergence in the two components
of the left-trivialized phase space.

\begin{figure}[h!]
\centering
\includegraphics[width=0.58\textwidth]{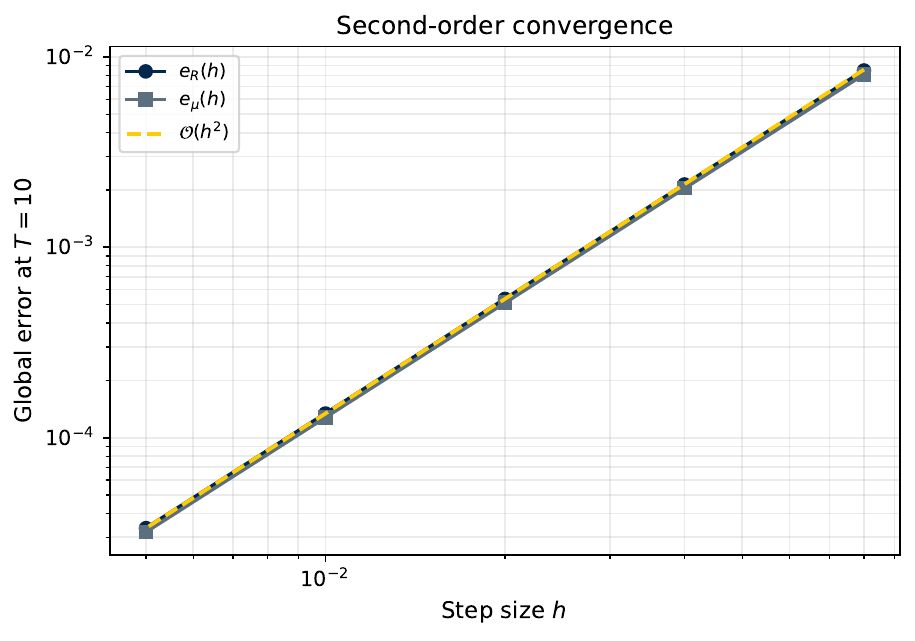}
\caption{Global configuration and momentum errors at $T=10$ as
functions of the step size $h$, together with a reference line
proportional to $h^2$.}
\label{fig:so3-convergence}
\end{figure}

\begin{table}[H]
\centering
\begin{tabular}{ccccc}
\hline
$h$ & $e_R(h)$ & $p_R(h)$ & $e_\mu(h)$ & $p_\mu(h)$\\
\hline
$0.080$ & $8.488\times10^{-3}$ & --- &
$8.107\times10^{-3}$ & ---\\
$0.040$ & $2.139\times10^{-3}$ & $1.9886$ &
$2.044\times10^{-3}$ & $1.9876$\\
$0.020$ & $5.358\times10^{-4}$ & $1.9971$ &
$5.122\times10^{-4}$ & $1.9969$\\
$0.010$ & $1.340\times10^{-4}$ & $1.9993$ &
$1.281\times10^{-4}$ & $1.9992$\\
$0.005$ & $3.351\times10^{-5}$ & $1.9999$ &
$3.203\times10^{-5}$ & $1.9998$\\
\hline
\end{tabular}
\caption{Final-time configuration and momentum errors and their
observed convergence orders at $T=10$. The symbol `---' indicates that
no coarser step size is available for computing the corresponding
observed order.}
\label{tab:so3-convergence}
\end{table}

\FloatBarrier

Figure~\ref{fig:so3-convergence} shows that both error curves are
parallel to the reference line proportional to $h^2$ over the complete
range of tested step sizes. Quantitatively,
Table~\ref{tab:so3-convergence} shows that the observed configuration
order increases from $1.9886$ to $1.9999$, while the observed momentum
order increases from $1.9876$ to $1.9998$. These results provide
numerical evidence for the second-order convergence of the
midpoint-Cayley variational integrator in both components of the
left-trivialized phase space, in agreement with
Corollary~\ref{cor:midpoint-cayley-order}.

To assess the robustness of the observed convergence with respect to the
time-dependence parameter, we repeated the same experiment for
$\kappa\in\{-10^{-3},-10^{-2},-5\times10^{-2},-10^{-1},
-2\times10^{-1},-4\times10^{-1},-8\times10^{-1}\}$,
keeping the initial data, final time, reference solver, and step sizes
unchanged. In all cases, both the configuration and momentum errors
exhibit second-order convergence. Using the two finest step sizes, the
observed orders remain within $1.9994\leq p_R\leq1.9999$ and
$1.9992\leq p_\mu\leq2.0001$. As shown in
Figure~\ref{fig:kappa-sweep}, the observed orders remain essentially
equal to two throughout the tested range. These results show that the
second-order behavior of the midpoint-Cayley variational integrator is
robust with respect to $\kappa$ and is not tied to the particular value
$\kappa=-0.4$ used in Figure~\ref{fig:so3-convergence} and Table~\ref{tab:so3-convergence}.

\begin{figure}[ht]
\centering
\includegraphics[width=0.72\textwidth]{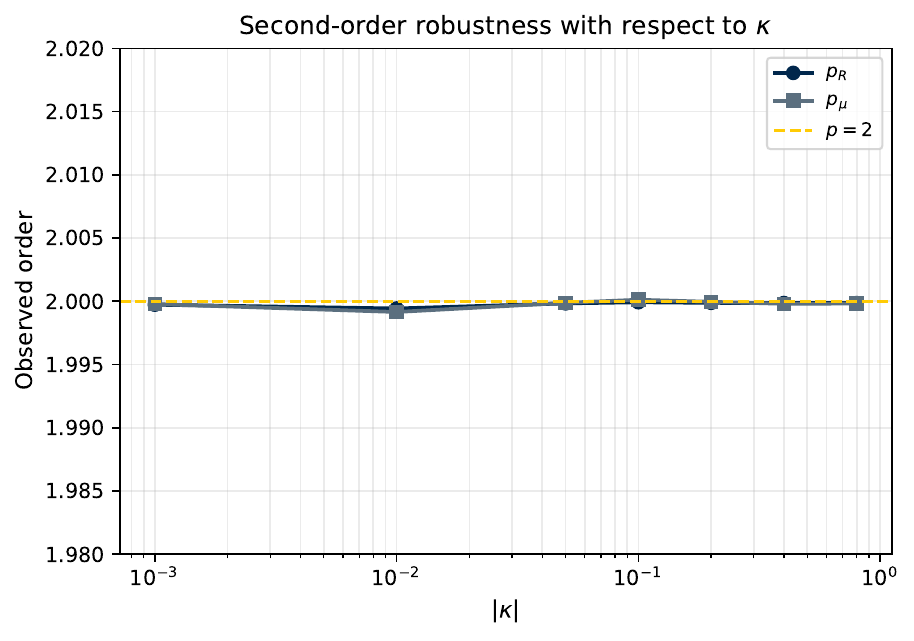}
\caption{Observed convergence orders for different values of $\kappa$.}
\label{fig:kappa-sweep}
\end{figure}

\end{example}

The phase-space map $\widetilde\Psi_d^k$ is the map that will be lifted
to the extended phase space in Section~\ref{sec:BEA}. The preceding
analysis shows that it is a
second-order symplectic map generated by the discrete
Lagrangian. In the left-invariant case, its reduction is a second-order
Poisson map on $\mathfrak g^*$ and, by
Theorem~\ref{thm:discrete-Lie-Poisson-reduction}, preserves coadjoint
orbits and Casimir functions exactly.

\section{Extended Hamiltonian formulation and backward error analysis}
\label{sec:BEA}

The time-dependent variational construction
developed in Sections~\ref{Sec:DiscreteNoether} and~\ref{sec:ham} yields a time-dependent
family of symplectic maps on
$T^*G\simeq G\times\mathfrak g^*$.
Since this family does not define a single
autonomous map on $T^*G$, we introduce
the extended cotangent bundle
$T^*(\mathbb R\times G)$ by adjoining
time and its conjugate momentum.
In this formulation, the extended Hamiltonian
is conserved along the exact autonomous flow,
whereas the original time-dependent
Hamiltonian is generally not conserved.

We first recall the continuous extended
Hamiltonian formulation and its Poisson
projection onto $\mathbb R\times T^*G$.
We then construct an exactly symplectic
and $\rho$-projectable extension of the
time-dependent discrete variational flow \cite{Asorey}.
This construction provides the basis for
establishing the affine dependence of the
formal modified Hamiltonian on the momentum
conjugate to time. Under suitable uniform analyticity and
coefficient bounds, and provided that
the numerical trajectories and the
corresponding modified flows remain
in the prescribed analytic domains,
we establish conditional long-time
near-conservation estimates.

\subsection{Extended Hamiltonian systems and Poisson projection}
\label{subsec:extended-Hamiltonian}

We retain the left trivialization
$T^*G\simeq G\times\mathfrak g^*$ and the
canonical forms $\theta_G$ and
$\omega_G=-d\theta_G$ introduced in
Subsection~\ref{sec4.1}. Let
$H:\mathbb R\times G\times\mathfrak g^*
\to\mathbb R$ be a smooth time-dependent
Hamiltonian.
Its fiber derivative
$\delta H/\delta\mu\in\mathfrak g$
and left-trivialized group derivative
$\mathbf{D}_gH\in\mathfrak g^*$
are defined by
\[
\left\langle\delta\mu,
\frac{\delta H}{\delta\mu}\right\rangle
=
\left.\frac{d}{ds}\right|_{s=0}
H(t,g,\mu+s\delta\mu),
\]
and
\[
\langle\mathbf{D}_gH,\eta\rangle
=
\left.\frac{d}{ds}\right|_{s=0}
H(t,g\exp(s\eta),\mu),
\]
for $\delta\mu\in\mathfrak g^*$ and
$\eta\in\mathfrak g$, respectively.

With the coadjoint convention~(5),
Hamilton's equations are
\[
\dot g=T_eL_g\frac{\delta H}{\delta\mu},
\qquad
\dot\mu=
\operatorname{ad}^*_{\delta H/\delta\mu}\mu
-\mathbf{D}_gH.
\]

Consider the extended configuration manifold
$\overline G:=\mathbb R\times G$ and its cotangent bundle
$\overline M:=T^*\overline G$. Using the preceding left
trivialization, we identify
$\overline M\simeq
\mathbb R\times\mathbb R\times G\times\mathfrak g^*$
and write its elements as $x=(t,\nu,g,\mu)$, where $\nu$ is the
momentum conjugate to time.

The canonical one-form and symplectic form on $\overline M$ are
\begin{equation*}
\label{eq:extended-canonical-forms}
\theta_{\mathrm{ext}}
=
\nu\,dt+\theta_G,
\qquad
\omega_{\mathrm{ext}}
=
-d\theta_{\mathrm{ext}}
=
dt\wedge d\nu+\omega_G.
\end{equation*}
For $\mathcal H\in C^\infty(\overline M)$, its Hamiltonian vector
field $X_{\mathcal H}$ is defined by
$\iota_{X_{\mathcal H}}\omega_{\mathrm{ext}}=d\mathcal H$.
We define the associated Poisson bracket by
$\{\mathcal F,\mathcal K\}_{\mathrm{ext}}
:=
\omega_{\mathrm{ext}}(X_{\mathcal F},X_{\mathcal K})$.
With these conventions,
$X_{\mathcal K}[\mathcal F]
=
\{\mathcal F,\mathcal K\}_{\mathrm{ext}}$.

The autonomous extended Hamiltonian associated with $H$ is
\begin{equation*}
\label{eq:extended-Hamiltonian}
H^{\mathrm{ext}}(t,\nu,g,\mu)
=
\nu+H(t,g,\mu).
\end{equation*}
Its Hamiltonian vector field satisfies
\begin{equation*}
\label{eq:extended-Hamilton-equations}
\dot t=1,
\qquad
\dot\nu=-\frac{\partial H}{\partial t},
\qquad
\dot g
=
T_e\mathcal L_g\frac{\delta H}{\delta\mu},
\qquad
\dot\mu
=
\operatorname{ad}_{\delta H/\delta\mu}^*\mu-\mathbf{D}_gH.
\end{equation*}
Thus the projection of an integral curve of
$X_{H^{\mathrm{ext}}}$ onto $(t,g,\mu)$ satisfies the original
non-autonomous Hamilton equations. Moreover,
$H^{\mathrm{ext}}$ is constant along the extended flow, since
$dH^{\mathrm{ext}}(X_{H^{\mathrm{ext}}})
=
\omega_{\mathrm{ext}}
(X_{H^{\mathrm{ext}}},X_{H^{\mathrm{ext}}})=0$.
By contrast, the physical Hamiltonian $H(t,g,\mu)$ is generally not
conserved.

Let
$\rho:\overline M\to\mathbb R\times T^*G$ be the surjective
submersion defined by
\begin{equation*}
\label{eq:extended-projection}
\rho(t,\nu,g,\mu)=(t,g,\mu).
\end{equation*}
The manifold $\mathbb R\times T^*G$ carries the product Poisson
structure for which time is a Casimir and whose restriction to every
constant-time fiber is the canonical Poisson structure of $T^*G$.
In the left-trivialized variables introduced above, this bracket is
\begin{equation}
\label{eq:time-dependent-Poisson-bracket}
\begin{split}
\{F,K\}_0(t,g,\mu)
={}&
\left\langle
\mathbf{D}_gF,
\frac{\delta K}{\delta\mu}
\right\rangle
-
\left\langle
\mathbf{D}_gK,
\frac{\delta F}{\delta\mu}
\right\rangle
-
\left\langle
\mu,
\left[
\frac{\delta F}{\delta\mu},
\frac{\delta K}{\delta\mu}
\right]
\right\rangle
\end{split}
\end{equation}
for $F,K\in C^\infty(\mathbb R\times T^*G)$. In particular,
$\{t,F\}_0=0$ for every
$F\in C^\infty(\mathbb R\times T^*G)$.

Since $\rho^*F$ and $\rho^*K$ are independent of $\nu$, direct
substitution in the extended canonical bracket gives
\begin{equation*}
\label{eq:rho-Poisson}
\{\rho^*F,\rho^*K\}_{\mathrm{ext}}
=
\rho^*\{F,K\}_0.
\end{equation*}
Hence $\rho$ is a Poisson map.

Let $\mathcal U,\mathcal V\subset\overline M$ be open subsets. A
local diffeomorphism $\Psi:\mathcal U\to\mathcal V$ is called
\textit{symplectic} if
$\Psi^*(\omega_{\mathrm{ext}}|_{\mathcal V})
=
\omega_{\mathrm{ext}}|_{\mathcal U}$.
It is called \textit{exactly symplectic} if there exists
$S\in C^\infty(\mathcal U)$ such that
\begin{equation}
\label{eq:exact-symplectic-map}
\Psi^*(\theta_{\mathrm{ext}}|_{\mathcal V})
-
\theta_{\mathrm{ext}}|_{\mathcal U}
=
dS.
\end{equation}
Every exactly symplectic map is symplectic, as follows by applying the
exterior derivative to
\eqref{eq:exact-symplectic-map}.

\begin{definition}
\label{def:projectable-extended-map}
Let $\mathcal U,\mathcal V\subset\overline M$ be open subsets. A smooth
map $\Psi:\mathcal U\to\mathcal V$ is called
\textit{$\rho$-projectable} if there exists a smooth map
$\Phi:\rho(\mathcal U)\to\rho(\mathcal V)$ such that
\begin{equation*}
\label{eq:projectable-map}
\rho\circ\Psi
=
\Phi\circ\rho|_{\mathcal U}.
\end{equation*}
In this case, $\Phi$ is called the projected map associated with
$\Psi$.
\end{definition}

\begin{proposition}
\label{prop:projected-structures-Lie}
Let $\mathcal U,\mathcal V\subset\overline M$ be open subsets, and let
$\Psi:\mathcal U\to\mathcal V$ be a symplectic,
$\rho$-projectable local diffeomorphism. Then its projected map
$\Phi:\rho(\mathcal U)\to\rho(\mathcal V)$ is a Poisson map with
respect to the bracket
\eqref{eq:time-dependent-Poisson-bracket}.
\end{proposition}

\begin{proof}
Let $F,K\in C^\infty(\rho(\mathcal V))$. Since $\Psi$ is symplectic,
its pullback preserves the extended Poisson bracket. Since $\rho$ is
Poisson and
$\rho\circ\Psi=\Phi\circ\rho|_{\mathcal U}$, we obtain
$\rho^*\Phi^*\{F,K\}_0
=
\rho^*\{\Phi^*F,\Phi^*K\}_0$
on $\mathcal U$.

The restriction
$\rho|_{\mathcal U}:\mathcal U\to\rho(\mathcal U)$ is surjective.
Consequently, its pullback on smooth functions is injective, and
therefore
$\Phi^*\{F,K\}_0
=
\{\Phi^*F,\Phi^*K\}_0$.
Thus $\Phi$ is a Poisson map.
\end{proof}

Suppose that the projected map advances time by a fixed step $h$, so
that $\Phi_h(t,z)=(t+h,\phi_{t,h}(z))$. For each $t$, define
$\mathcal D_t:=\{z\in T^*G:(t,z)\in\rho(\mathcal U)\}$ and
$\mathcal D'_{t+h}:=
\{z\in T^*G:(t+h,z)\in\rho(\mathcal V)\}$.

\begin{corollary}
\label{cor:fiberwise-symplectic}
Under the hypotheses of
Proposition~\ref{prop:projected-structures-Lie}, each map
$\phi_{t,h}:\mathcal D_t\to\mathcal D'_{t+h}$ is a symplectic local
diffeomorphism onto its image.
\end{corollary}

\begin{proof}
Let $f$ and $k$ be smooth functions defined on a neighborhood of
$\phi_{t,h}(\mathcal D_t)$ in $T^*G$, and regard them locally as
functions on $\mathbb R\times T^*G$ that are independent of time.
Since $\Phi_h$ is Poisson,
\begin{equation}
\label{eq:fiberwise-Poisson}
\{f\circ\phi_{t,h},k\circ\phi_{t,h}\}_{T^*G}
=
\{f,k\}_{T^*G}\circ\phi_{t,h}.
\end{equation}
Because the Poisson tensor of $T^*G$ is nondegenerate,
\eqref{eq:fiberwise-Poisson} implies that
$T_z\phi_{t,h}$ has full rank for every $z\in\mathcal D_t$.
Since the source and target have the same dimension,
$\phi_{t,h}$ is a local diffeomorphism. A Poisson local
diffeomorphism between symplectic manifolds is symplectic.
Consequently,
$\phi_{t,h}^*\omega_G
=
\omega_G|_{\mathcal D_t}$.
\end{proof}

\subsection{Exact symplectic extension of the discrete variational flow}
\label{subsec:extended-discrete-map}

Let
$L_d:\mathbb R\times G\times G\times(0,h_0)\to\mathbb R$
be a smooth time-dependent discrete Lagrangian. For fixed $h$, the
value $L_d(t,g_0,g_1;h)$ approximates the exact action over
$[t,t+h]$. The family introduced in
Section~\ref{Sec:DiscreteNoether} is recovered by setting
$L_d^k(g_k,g_{k+1})=L_d(t_k,g_k,g_{k+1};h)$.

For fixed $(t,h)$, define the negative and positive discrete Legendre
transformations by
\begin{equation*}
\label{eq:time-dependent-discrete-Legendre}
\begin{aligned}
\mathbb F^-L_d(t,g_0,g_1;h)
&=
\bigl(g_0,-D_1L_d(t,g_0,g_1;h)\bigr),\\
\mathbb F^+L_d(t,g_0,g_1;h)
&=
\bigl(g_1,D_2L_d(t,g_0,g_1;h)\bigr).
\end{aligned}
\end{equation*}
Here $D_1L_d(t,g_0,g_1;h)\in T_{g_0}^*G$ and
$D_2L_d(t,g_0,g_1;h)\in T_{g_1}^*G$ denote differentiation with
respect to the first and second Lie group arguments, respectively.

The discrete Lagrangian is called \textit{regular} at
$(t,g_0,g_1;h)$ if
$\mathbb F^-L_d(t,\cdot,\cdot;h)$ is a local diffeomorphism at
$(g_0,g_1)$. Equivalently, the mixed derivative
$D_{12}L_d(t,g_0,g_1;h)$, regarded as a linear map from
$T_{g_1}G$ to $T_{g_0}^*G$, is an isomorphism. Since
$D_{21}L_d$ is the transpose of $D_{12}L_d$, this condition is also
equivalent to the local invertibility of
$\mathbb F^+L_d(t,\cdot,\cdot;h)$.

Given $(t,\nu_0,g_0,\mu_0)$, regularity determines $g_1$ locally
through the first momentum relation in
\begin{equation}
\label{eq:extended-discrete-Legendre-relations}
\begin{aligned}
t_1
&=
t+h,\\
\mu_0
&=
-\bigl(T_e\mathcal L_{g_0}\bigr)^*
D_1L_d(t,g_0,g_1;h),\\
\mu_1
&=
\bigl(T_e\mathcal L_{g_1}\bigr)^*
D_2L_d(t,g_0,g_1;h),\\
\nu_1
&=
\nu_0+
\frac{\partial L_d}{\partial t}(t,g_0,g_1;h).
\end{aligned}
\end{equation}
These relations define a local map
$\Psi_h^{\mathrm{ext}}(t,\nu_0,g_0,\mu_0)
=(t_1,\nu_1,g_1,\mu_1)$ on every open set on which the first momentum
relation admits a unique smooth local solution for $g_1$.

\begin{proposition}
\label{prop:extended-map-Lie}
On every such regularity domain,
$\Psi_h^{\mathrm{ext}}$ is exactly symplectic and
$\rho$-projectable. Its projected map is
\begin{equation}
\label{eq:projected-discrete-map}
\Phi_h(t,g_0,\mu_0)
=
\bigl(t+h,\Psi_d^{t,h}(g_0,\mu_0)\bigr),
\end{equation}
where $\Psi_d^{t,h}$ is the discrete Hamiltonian map generated by
$L_d(t,\cdot,\cdot;h)$.
\end{proposition}

\begin{proof}
Let $p_0\in T_{g_0}^*G$ and $p_1\in T_{g_1}^*G$ be the
nontrivialized momenta defined by
$p_0=-D_1L_d(t,g_0,g_1;h)$ and
$p_1=D_2L_d(t,g_0,g_1;h)$. Their left-trivialized representations are
$\mu_i=(T_e\mathcal L_{g_i})^*p_i$, for $i=0,1$.

Let $\Gamma_h\subset\overline M\times\overline M$ be the graph of
$\Psi_h^{\mathrm{ext}}$, let
$\iota_h:\Gamma_h\hookrightarrow
\overline M\times\overline M$ denote the inclusion, and let
$\operatorname{pr}_1$ and $\operatorname{pr}_2$ be the projections
onto the first and second factors. Locally, regularity allows
$\Gamma_h$ to be parametrized by $(t,\nu_0,g_0,g_1)$.

Since $t_1=t+h$ and $h$ is fixed, one has $dt_1=dt$. Therefore, the
relations
\eqref{eq:extended-discrete-Legendre-relations} give
\begin{equation}
\label{eq:generating-function-identity}
\iota_h^*
\left(
\operatorname{pr}_2^*\theta_{\mathrm{ext}}
-
\operatorname{pr}_1^*\theta_{\mathrm{ext}}
\right)
=
dL_d(t,g_0,g_1;h),
\end{equation}
where the differential on the right-hand side is taken with respect
to $(t,g_0,g_1)$ while $h$ is held fixed.

Indeed, when evaluated on a tangent variation
$(\delta t,\delta\nu_0,\delta g_0,\delta g_1)$ of the graph, the
left-hand side of
\eqref{eq:generating-function-identity} is
$(\nu_1-\nu_0)\delta t
-\langle p_0,\delta g_0\rangle
+\langle p_1,\delta g_1\rangle$.
Using
\eqref{eq:extended-discrete-Legendre-relations}, this becomes
$(\partial L_d/\partial t)\delta t
+\langle D_1L_d,\delta g_0\rangle
+\langle D_2L_d,\delta g_1\rangle$, which is precisely the
differential of $L_d$.

By regularity, $g_1$ is a smooth function of
$(t,g_0,\mu_0)$ on the domain of the map. Define
$S_h(t,\nu_0,g_0,\mu_0)
:=
L_d(t,g_0,g_1(t,g_0,\mu_0);h)$.
Pulling
\eqref{eq:generating-function-identity} back through the
parametrization of $\Gamma_h$ by its first component yields
\begin{equation}
\label{eq:exact-symplectic-extended-map}
(\Psi_h^{\mathrm{ext}})^*\theta_{\mathrm{ext}}
-
\theta_{\mathrm{ext}}
=
dS_h.
\end{equation}
Thus $\Psi_h^{\mathrm{ext}}$ is exactly symplectic in the sense of
\eqref{eq:exact-symplectic-map}. Applying the exterior derivative to
\eqref{eq:exact-symplectic-extended-map} gives
$(\Psi_h^{\mathrm{ext}})^*\omega_{\mathrm{ext}}
=\omega_{\mathrm{ext}}$.

The equation that determines $g_1$ from $(t,g_0,\mu_0)$ is independent
of $\nu_0$. Consequently, $g_1$ and $\mu_1$ depend only on
$(t,g_0,\mu_0)$, whereas $\nu_1$ depends affinely on $\nu_0$.
Therefore,
$\rho\circ\Psi_h^{\mathrm{ext}}
=
\Phi_h\circ\rho$,
with $\Phi_h$ given by
\eqref{eq:projected-discrete-map}. This proves
$\rho$-projectability.
\end{proof}

The fixed-step discrete Lagrangian
$L_d(t,g_0,g_1;h)$ is defined on
$\mathbb R\times G\times G$ for each fixed $h$. It should not be
identified directly with an ordinary discrete Lagrangian on an open
subset of $(\mathbb R\times G)^2$, because the relation
$t_1-t_0=h$ restricts the two time variables to a codimension-one
submanifold. Consequently, the fixed-step function used here does not,
by itself, define Legendre transformations with respect to two
independent time variables.

One may instead introduce a discrete Lagrangian
$\overline L_d(t_0,g_0,t_1,g_1)$ with $t_0$ and $t_1$ treated as
independent  \cite{Kane}. Stationarity with respect to the discrete time
variables then produces an additional discrete energy equation and,
in general, a variable-time variational integrator. This is a
different construction from the fixed-step method studied here.

The relations
\eqref{eq:extended-discrete-Legendre-relations} do not claim that such
a variable-time formulation is impossible or necessarily degenerate.
They provide a direct generating-function extension of the prescribed
fixed-step map, retaining $t_1=t_0+h$ while adjoining the conjugate
momentum $\nu$. The resulting extension is exactly symplectic and
projects to the original non-autonomous discrete Hamiltonian map.

Subsections~\ref{subsec:extended-Hamiltonian} and
\ref{subsec:extended-discrete-map} establish the geometric
infrastructure required for backward error analysis: the
fixed-step non-autonomous variational map has been lifted to an
exactly symplectic, $\rho$-projectable map on the extended phase
space. No formal interpolation or exponentially small estimate has
yet been asserted. The formal modified vector field and its
Hamiltonian are constructed in the next subsection through the
formal logarithm of the pullback operator. The analytic domains,
coefficient estimates, optimal truncation index
$N(h)=\lfloor c_0/h\rfloor$, and exponentially small remainder are
introduced separately in Subsection~\ref{subsec:exponential-BEA}.

\subsection{Modified Hamiltonian and its affine structure}
\label{subsec:formal-modified-Hamiltonian}

We construct the modified vector field associated
with the extended variational map as a formal
power series in the step size and establish the
existence of a corresponding modified Hamiltonian
with affine dependence on the momentum conjugate
to time. No convergence of the series is asserted.
Analytic truncation and exponentially small
estimates are considered separately in
Subsection~\ref{subsec:exponential-BEA}.

\begin{definition}
Let $\mathcal U\subset\overline M$ be an open set.
A function $f\in C^\infty(\mathcal U)$ is called
\textit{$\rho$-basic} if there exists
$F\in C^\infty(\rho(\mathcal U))$ such that
$f=\rho^*F$ on $\mathcal U$.
\end{definition}

Every $\rho$-basic function satisfies
$\partial f/\partial\nu=0$. Conversely, if every
nonempty fiber of the restricted projection
$\rho|_{\mathcal U}:\mathcal U\to\rho(\mathcal U)$
is connected, then $f\in C^\infty(\mathcal U)$
is $\rho$-basic if and only if
$\partial f/\partial\nu=0$ on $\mathcal U$.

\begin{definition}
\label{def:formal-modified-fields}
Let $\mathcal U\subset\overline M$ be an open set.
We denote by $C^\infty(\mathcal U)[[h]]$
the algebra of formal power series in $h$
with coefficients in $C^\infty(\mathcal U)$,
where equality is understood coefficientwise.

A \textit{formal vector field} on $\mathcal U$
is a series
\[
X_h=\sum_{r\geq0}h^rX_r,
\qquad X_r\in\mathfrak X(\mathcal U).
\]
A \textit{formal Hamiltonian} is a series
\[
\mathcal H_h=\sum_{r\geq0}h^r\mathcal H_r,
\qquad
\mathcal H_r\in C^\infty(\mathcal U),
\]
whose associated formal Hamiltonian vector
field is determined coefficientwise by
\[
\iota_{X_{\mathcal H_r}}\omega_{\mathrm{ext}}
=d\mathcal H_r.
\]
\end{definition}
\begin{definition}\label{def5.4}
Let $\Psi_h:\mathcal U\to\mathcal U$ be a
smooth family of local diffeomorphisms with
$\Psi_0=\operatorname{id}_{\mathcal U}$.
A formal vector field $X_h$ is called a
\textit{formal modified vector field}
of $\Psi_h$ if
$\Psi_h^*=\exp(hX_h)$ as formal operators on
$C^\infty(\mathcal U)[[h]]$.

If $X_h$ is Hamiltonian, a formal
Hamiltonian $\mathcal H_h$ satisfying
\[
\iota_{X_h}\omega_{\mathrm{ext}}
=d\mathcal H_h
\]
is called a \textit{formal modified
Hamiltonian} of $\Psi_h$.
\end{definition}

For a family $\Psi_h$ as in
Definition~\ref{def5.4},
the formal logarithm of its pullback
operator is defined by
\begin{equation}
\label{eq:formal-logarithm}
\mathcal D_h
:=
\frac{1}{h}\log(\Psi_h^*)
=
\frac{1}{h}
\sum_{m=1}^{\infty}
\frac{(-1)^{m+1}}{m}
\bigl(\Psi_h^*-\operatorname{id}\bigr)^m.
\end{equation}
This expression is well defined as a formal series because
$\Psi_h^*-\operatorname{id}=O(h)$; consequently, only finitely many
terms contribute to the coefficient of any prescribed power of $h$.

\begin{lemma}
\label{lem:formal-generator-projectable}
Let $\Psi_h:\mathcal U\to\mathcal U$ be a smooth family of exactly
symplectic local diffeomorphisms such that
$\Psi_0=\operatorname{id}_{\mathcal U}$. Assume that the primitives in
the exact-symplectic identities can be chosen smoothly in $h$, so that
\[
\Psi_h^*\theta_{\mathrm{ext}}-\theta_{\mathrm{ext}}=dS_h,
\qquad
S_0=0.
\]
If $\Psi_h$ is $\rho$-projectable and satisfies
$t\circ\Psi_h=t+h$, then the following statements hold.

\begin{enumerate}
\item The formal operator $\mathcal D_h$ is a derivation of
$C^\infty(\mathcal U)[[h]]$ and hence determines a unique formal vector
field $X_h$.

\item The formal vector field $X_h$ is Hamiltonian on $\mathcal U$:
there exists
$\overline H^{\mathrm{ext}}\in C^\infty(\mathcal U)[[h]]$ such that
$\iota_{X_h}\omega_{\mathrm{ext}}
=d\overline H^{\mathrm{ext}}$.

\item The time component of the formal vector field is exact:
$X_h[t]=1$.

\item The operator $X_h$ maps $\rho$-basic formal functions to
$\rho$-basic formal functions. Consequently, $X_h$ projects to a formal
vector field on $\rho(\mathcal U)$.
\end{enumerate}
\end{lemma}

\begin{proof}
Set $P_h:=\Psi_h^*$. The pullback $P_h$
is a unital algebra homomorphism satisfying
$P_h=\mathrm{id}+O(h)$, and is therefore
formally invertible. 
To show that its formal logarithm is a
derivation, introduce the formal
one-parameter family
$A_{s,h}:=\exp(s\log P_h)$.
For every nonnegative integer $n$, one has
$A_{n,h}=P_h^n$, and therefore
\[
A_{n,h}(fk)=A_{n,h}(f)A_{n,h}(k).
\]
At every fixed order in $h$, both sides
of the corresponding identity for
$A_{s,h}$ are polynomials in $s$.
Since they agree for all nonnegative
integers, they agree identically in $s$.
Differentiating at $s=0$ shows that
$\log P_h(fk)=(\log P_hf)k+f(\log P_hk)$.
Hence $\mathcal D_h=h^{-1}\log P_h$
is a derivation.

Every derivation of $C^\infty(\mathcal U)$ is induced by a unique smooth
vector field. Applying this fact coefficientwise yields a unique formal
vector field $X_h$ satisfying $\mathcal D_hf=X_h[f]$ for every
$f\in C^\infty(\mathcal U)[[h]]$.

It remains to prove that $X_h$ is Hamiltonian. This conclusion does not
follow from symplecticity alone, because a symplectic vector field need
not be globally Hamiltonian on a domain with nontrivial first de Rham
cohomology. The exact-symplectic hypothesis is therefore essential.

The logarithm of the pullback on differential forms is the Lie
derivative of the formal vector field:
$\log(\Psi_h^*)=h\mathcal L_{X_h}$. Indeed, the logarithm is a
degree-zero derivation of the algebra of differential forms, it commutes
with the exterior derivative, and its restriction to functions is
$hX_h$. Therefore,
$\Psi_h^*=\exp(h\mathcal L_{X_h})$ as a formal operator on differential
forms.

Define the formally invertible operator
\[
\mathcal Q_h
:=
\frac{\exp(h\mathcal L_{X_h})-\operatorname{id}}
     {h\mathcal L_{X_h}}
=
\sum_{m=0}^{\infty}
\frac{h^m}{(m+1)!}\mathcal L_{X_h}^{\,m}.
\]
Its constant term is the identity. Using exact symplecticity, we obtain
\[
dS_h
=
\bigl(\Psi_h^*-\operatorname{id}\bigr)\theta_{\mathrm{ext}}
=
h\mathcal Q_h\mathcal L_{X_h}\theta_{\mathrm{ext}}.
\]
Since $S_0=0$, the quotient $S_h/h$ is a well-defined formal function.
Moreover, $\mathcal L_{X_h}$ commutes with $d$, and hence so do
$\mathcal Q_h$ and its formal inverse. It follows that
$\mathcal L_{X_h}\theta_{\mathrm{ext}}=dB_h$, where
$B_h:=h^{-1}\mathcal Q_h^{-1}S_h$.

Cartan's identity and the convention
$d\theta_{\mathrm{ext}}=-\omega_{\mathrm{ext}}$ give
\[
\mathcal L_{X_h}\theta_{\mathrm{ext}}
=
\iota_{X_h}d\theta_{\mathrm{ext}}
+
d\bigl(\iota_{X_h}\theta_{\mathrm{ext}}\bigr)
=
-\iota_{X_h}\omega_{\mathrm{ext}}
+
d\bigl(\iota_{X_h}\theta_{\mathrm{ext}}\bigr).
\]
Combining this identity with
$\mathcal L_{X_h}\theta_{\mathrm{ext}}=dB_h$ yields
$\iota_{X_h}\omega_{\mathrm{ext}}
=d\overline H^{\mathrm{ext}}$, with
$\overline H^{\mathrm{ext}}
:=\iota_{X_h}\theta_{\mathrm{ext}}-B_h$.
Thus $X_h$ is a formal Hamiltonian vector field on $\mathcal U$. Notice
that no application of the Poincaré lemma and no contractibility
assumption on $\mathcal U$ are required.

The exact time update gives
$\Psi_h^*t=t+h$. Consequently,
$(\Psi_h^*-\operatorname{id})t=h$, while
$(\Psi_h^*-\operatorname{id})^mt=0$ for every $m\geq2$, because pullback
preserves constants. Substitution into
\eqref{eq:formal-logarithm} gives $\mathcal D_ht=1$, and therefore
$X_h[t]=1$.

Finally, $\rho$-projectability gives
$\Psi_h^*\rho^*F=\rho^*\Phi_h^*F$ for every
$F\in C^\infty(\rho(\mathcal U))$. Thus $\Psi_h^*$ preserves the
subalgebra of $\rho$-basic functions. The same is true of
$\Psi_h^*-\operatorname{id}$, all its powers, and its formal logarithm.
Hence $X_h$ maps $\rho$-basic formal functions to $\rho$-basic formal
functions. This is precisely the infinitesimal criterion for
projectability of the formal vector field.
\end{proof}

The distinction between symplectic and exactly symplectic maps in the
preceding lemma cannot be omitted. For example, on
$T^*S^1$ with canonical coordinates $(q,p)$, the translation
$(q,p)\mapsto(q,p+c)$ is symplectic, but
$\Psi^*\theta-\theta=c\,dq$ is not exact when $c\neq0$. Its infinitesimal
generator is symplectic but not globally Hamiltonian. Exact
symplecticity rules out this obstruction and is what permits the global
construction of $\overline H^{\mathrm{ext}}$ on $\mathcal U$.

\begin{definition}
\label{def:extended-order}
Let $\varphi_h^{H^{\mathrm{ext}}}$ denote
the time-$h$ flow of the exact extended
Hamiltonian vector field.
A family of maps $\Psi_h^{\mathrm{ext}}$
is said to have \textit{order $p$}
with respect to this flow if, for every
coordinate chart
$\chi:\mathcal W\to W\subset\mathbb R^d$
and every compact subset
$\mathcal K\subset\mathcal W$ on which
both maps take values in $\mathcal W$,
\begin{equation*}
\sup_{x\in\mathcal K}
\left\|
\chi\bigl(\Psi_h^{\mathrm{ext}}(x)\bigr)
-
\chi\bigl(
\varphi_h^{H^{\mathrm{ext}}}(x)
\bigr)
\right\|
=O(h^{p+1})
\end{equation*}
as $h\to0$.
This definition is independent of the
chosen chart because the derivatives
of smooth transition maps are bounded
on compact subsets.
\end{definition}

\begin{theorem}
\label{thm:affine-modified-Hamiltonian}
Let $\Psi_h^{\mathrm{ext}}$ be the extended variational map of
Proposition~\ref{prop:extended-map-Lie}. Assume that it extends
smoothly to $h=0$, with
$\Psi_0^{\mathrm{ext}}=\operatorname{id}$, and that it has order $p$ for
the flow of
$H^{\mathrm{ext}}=\nu+H(t,g,\mu)$. 
On every open set
$\mathcal U\subset\overline M$
satisfying the hypotheses of
Lemma~\ref{lem:formal-generator-projectable}
and such that every nonempty fiber of
$\rho|_{\mathcal U}$ is connected,
its formal modified Hamiltonian can
be chosen in the form

\begin{equation}
\label{eq:modified-H-reduced}
\overline H^{\mathrm{ext}}(t,\nu,g,\mu;h)
=
\nu
+
H(t,g,\mu)
+
\sum_{r=p}^{\infty}h^rH_r(t,g,\mu).
\end{equation}
In particular, every correction coefficient $H_r$ is independent of
$\nu$.
\end{theorem}

\begin{proof}
By Proposition~\ref{prop:extended-map-Lie}, the map
$\Psi_h^{\mathrm{ext}}$ is exactly symplectic,
$\rho$-projectable, and satisfies
$t\circ\Psi_h^{\mathrm{ext}}=t+h$. Therefore,
Lemma~\ref{lem:formal-generator-projectable} provides a formal vector
field $X_h$ and a formal Hamiltonian
$\overline H^{\mathrm{ext}}$ satisfying
$\iota_{X_h}\omega_{\mathrm{ext}}
=d\overline H^{\mathrm{ext}}$ and $X_h[t]=1$.

For the convention
$\omega_{\mathrm{ext}}=dt\wedge d\nu+\omega_G$, the coefficient of
$d\nu$ in
$\iota_{X_{\mathcal H}}\omega_{\mathrm{ext}}=d\mathcal H$ gives
$X_{\mathcal H}[t]=\partial\mathcal H/\partial\nu$. Hence
$X_h[t]=1$ implies
$\partial\overline H^{\mathrm{ext}}/\partial\nu=1$. It follows that
$\overline H^{\mathrm{ext}}-\nu$ is independent of $\nu$. Since the
fibers of $\rho|_{\mathcal U}$ are connected, there exists a formal
function $\overline H(t,g,\mu;h)$ on $\rho(\mathcal U)$ such that
$\overline H^{\mathrm{ext}}=\nu+\overline H(t,g,\mu;h)$.

The Hamiltonian associated with a given Hamiltonian vector field is
determined up to a formal constant on each connected component of
$\mathcal U$. We fix this freedom by prescribing the value of each
coefficient at one reference point in every connected component.

The exact flow satisfies
$(\varphi_h^{H^{\mathrm{ext}}})^*
=\exp(h\mathcal L_{X_{H^{\mathrm{ext}}}})$. Since the method has order
$p$, its pullback and the pullback of the exact flow agree through order
$h^p$. Taking their formal logarithms and dividing by $h$ shows that
$X_h-X_{H^{\mathrm{ext}}}=O(h^p)$. Nondegeneracy of
$\omega_{\mathrm{ext}}$ then implies
$d(\overline H^{\mathrm{ext}}-H^{\mathrm{ext}})=O(h^p)$. After fixing
the formal additive constants as above, the coefficients of orders
$1,\dots,p-1$ vanish. Therefore,
$\overline H(t,g,\mu;h)
=H(t,g,\mu)+\sum_{r=p}^{\infty}h^rH_r(t,g,\mu)$, which proves
\eqref{eq:modified-H-reduced}.
\end{proof}

The affine dependence on $\nu$ is a consequence of the exact discrete
time translation $t\mapsto t+h$. Exact symplecticity guarantees the
global existence of the formal Hamiltonian on the chosen domain,
whereas $\rho$-projectability ensures that the associated formal vector
field descends to a formal non-autonomous vector field on
$\mathbb R\times T^*G$. These properties play distinct roles in the
construction. On domains with connected
$\rho$-fibers, however, symplecticity together
with the exact time translation
$t\circ\Psi_h=t+h$ already implies
$\rho$-projectability.

\subsection{Analytic domains and exponentially accurate interpolation}
\label{subsec:exponential-BEA}

The formal series constructed in
Subsection~\ref{subsec:formal-modified-Hamiltonian} need not converge.
Backward error analysis instead truncates this series at an order that
depends on the step size. Under uniform analyticity assumptions, the
optimally truncated modified vector field interpolates the discrete map
with an exponentially small one-step error. We now state precisely the
domains and bounds required for this conclusion.

Let $d=2(\dim G+1)$ and let
$\chi:\mathcal V\to V\subset\mathbb R^d$ be a real-analytic coordinate
chart on the extended phase space $\overline M$. Let
$\mathcal K\subset\mathcal V$ be compact and write
$K:=\chi(\mathcal K)$. We assume that $K$ has nonempty interior and is
contained in a bounded subset of $\mathbb R^d$.

For $\rho>0$, define the complex $\rho$-neighborhood of $K$ by
\(K_\rho:=\{z\in\mathbb C^d:
\inf_{y\in K}\|z-y\|_\infty<\rho\}\), where
\(\|z\|_\infty:=\max_{1\leq j\leq d}|z_j|\).
If $0<\rho_2<\rho_1<\rho$, then
$K_{\rho_2}\subset K_{\rho_1}\subset K_\rho$. We write
$K_\rho^{\mathbb R}:=K_\rho\cap\mathbb R^d$.

For a bounded holomorphic map
$f:K_\rho\to\mathbb C^m$, its analytic supremum norm is
$\|f\|_\rho:=\sup_{z\in K_\rho}\|f(z)\|_\infty$. The same notation is
used for the coordinate representation of a holomorphic vector field.
Distances between real points contained in the chart are measured by
$d_\chi(x,y):=\|\chi(x)-\chi(y)\|_\infty$.

The use of the two smaller radii $\rho_1$ and $\rho_2$ is essential.
The formal logarithm is estimated on $K_{\rho_1}$ after losing analytic
width from $K_\rho$, whereas derivatives of the truncated modified
Hamiltonian are estimated on $K_{\rho_2}$ by applying Cauchy's estimate
between $K_{\rho_2}$ and $K_{\rho_1}$.

A single chart need not contain an entire long-time trajectory. We
therefore use the following uniformity assumption.

\begin{definition}
\label{def:uniform-analytic-family}
A sequence of extended numerical states
$x_k=(t_k,\nu_k,g_k,\mu_k)$ is said to admit a
\textit{uniformly controlled analytic family of charts} if there exist
numbers $0<\rho_2<\rho_1<\rho,\, h_0>0$,
a compact set $K\subset\mathbb R^d$, and constants $M,L>0$ such that,
for every $k$, there is a real-analytic chart
$\chi_k:\mathcal V_k\to V_k\subset\mathbb R^d$ with the following
properties.

\begin{enumerate}
\item The states $x_k$ and
$x_{k+1}=\Psi_h^{\mathrm{ext}}(x_k)$
belong to
$\chi_k^{-1}(K_{\rho_2}^{\mathbb R})$.

\item The chart $\chi_k$, its inverse, the symplectic form
$\omega_{\mathrm{ext}}$, and the associated Poisson tensor admit
holomorphic coordinate representations on $K_\rho$. Their first
derivatives and matrix norms are bounded there by $L$, independently of
$k$.

\item The coordinate representation
$\widehat\Psi_{h,k}:=
\chi_k\circ\Psi_h^{\mathrm{ext}}\circ\chi_k^{-1}$
is jointly holomorphic in $(z,h)$ for
$z\in K_\rho$ and complex $h$ in a neighborhood of the origin, and
satisfies
\[
\|\widehat\Psi_{h,k}-\operatorname{id}\|_\rho
\leq M|h|,
\qquad |h|<h_0,
\]
uniformly in $k$.

\item The normalized coordinate representation
\[
\widehat H_k^{\mathrm{ext}}(z)
:=
H^{\mathrm{ext}}\bigl(\chi_k^{-1}(z)\bigr)
-
H^{\mathrm{ext}}(x_k)
\]
extends holomorphically to $K_\rho$ and satisfies
$\|\widehat H_k^{\mathrm{ext}}\|_\rho\leq M$ uniformly in $k$.
\end{enumerate}
\end{definition}

The normalization in item~4 does not affect the Hamiltonian vector
field. It is required because the conjugate momentum $\nu_k$ need not
remain bounded even when the projected trajectory
$(t_k,g_k,\mu_k)$ stays in a controlled region.

For periodic time dependence, the time variable may instead be regarded
as an element of the circle $\mathbb T_T=\mathbb R/T\mathbb Z$. This
does not remove the need for uniform complex neighborhoods in the group
and momentum variables. In particular,
$T^*(\mathbb T_T\times G)$ is generally noncompact.

We shall use the following analytic interpolation result. It is the
local form of the standard interpolation theorem for analytic
near-identity maps used in backward error analysis
\cite{BenettinGiorgilli1994, hairer,hansen,Faltinsen2000}.

\begin{lemma}
\label{lem:analytic-interpolation}
Let $0<\rho_2<\rho_1<\rho$, and let
$\widehat\Psi_h$ be jointly holomorphic in $(z,h)$
for $z\in K_\rho$ and $|h|<h_0$, and real-valued
for real $z$ and $h$. Assume that
$\widehat\Psi_0=\operatorname{id}$ and
\[
\|\widehat\Psi_h-\operatorname{id}\|_\rho
\leq M|h|,
\qquad |h|<h_0.
\]

Let
$X_h=\sum_{r\geq0}h^rX_r$ be the formal
vector field defined by
$h^{-1}\log(\widehat\Psi_h^*)$.
Then there exist constants
$A,B,c_0,c,C,h_*>0$, depending only on
$M$, $h_0$, $\rho-\rho_1$,
$\rho_1-\rho_2$, and the dimension $d$,
such that
\begin{equation}
\label{eq:modified-vector-field-factorial-bound}
\|X_r\|_{\rho_1}
\leq AB^r r!,
\qquad r\geq0.
\end{equation}

The constant $c_0$ may be chosen so that
$Bc_0<1$. For $0<h<h_*$, set
$N(h):=\lfloor c_0/h\rfloor$ and
\[
X_h^{[N(h)]}
:=
\sum_{r=0}^{N(h)}h^rX_r.
\]
For sufficiently small $h$, the time-$s$ flow
$\varphi_s^{[N(h)]}$ exists for every initial
point $x\in K_{\rho_2}^{\mathbb R}$ and
$0\leq s\leq h$, and remains in
$K_{\rho_1}^{\mathbb R}$. Moreover,
\begin{equation}
\label{eq:analytic-interpolation-defect}
\sup_{x\in K_{\rho_2}^{\mathbb R}}
\left\|
\widehat\Psi_h(x)
-
\varphi_h^{[N(h)]}(x)
\right\|_\infty
\leq Ch\exp(-c/h).
\end{equation}
\end{lemma}

\begin{proof}
Write
\[
\widehat\Psi_h(z)
=
z+\sum_{j\geq1}h^jd_j(z).
\]
Joint analyticity and the assumed near-identity
bound imply, by Cauchy's estimates in $h$,
\[
\|d_j\|_\rho\leq M h_0^{1-j},
\qquad j\geq1.
\]
In particular, $f:=d_1$ is analytic and
uniformly bounded on $K_\rho$.

Choose $R>0$ such that
\[
4R<
\min\{\rho-\rho_1,\rho_1-\rho_2\}.
\]
Every point of $K_{\rho_1}$ has a complex
neighborhood of radius $2R$ contained in
$K_\rho$, while every point of
$K_{\rho_2}^{\mathbb R}$ has such a
neighborhood contained in $K_{\rho_1}$.

The coefficient estimates for modified
equations obtained by the Cauchy-majorant
argument of Theorem IX.7.5 in
\cite{hairer} therefore apply uniformly
on these neighborhoods. Since the formal
modified vector field is uniquely determined
by the logarithm of the pullback, they yield
\eqref{eq:modified-vector-field-factorial-bound},
with constants independent of the center.

Choose $c_0>0$ sufficiently small that
$Bc_0<1$. The factorial bound then gives
a uniform bound for $X_h^{[N(h)]}$
on $K_{\rho_1}$.
For sufficiently small $h$, its trajectories
starting in $K_{\rho_2}^{\mathbb R}$
exist throughout $[0,h]$ and remain in
$K_{\rho_1}^{\mathbb R}$.

Finally, the analytic interpolation argument
of Theorem IX.7.6 in \cite{hairer},
based on the coincidence of the Taylor
coefficients of the numerical map and the
truncated modified flow, applies uniformly
to the above neighborhoods.
The resulting constants depend only on the
stated analytic bounds, domain widths,
and dimension. This proves
\eqref{eq:analytic-interpolation-defect}.
\end{proof}

The factorial growth in
\eqref{eq:modified-vector-field-factorial-bound} results from the
successive differentiations appearing in the formal logarithm and from
Cauchy's estimates on the nested domains
$K_{\rho_1}\subset K_\rho$. Choosing
$N(h)$ proportional to $h^{-1}$ stops the asymptotic series near its
smallest term and gives the exponentially small remainder in
\eqref{eq:analytic-interpolation-defect}. The lemma is a local analytic
result and does not, by itself, guarantee that a long numerical
trajectory remains in the required domains.

We now apply the interpolation result to the exactly symplectic
extension constructed in
Proposition~\ref{prop:extended-map-Lie}. The next theorem states the
uniform assumptions needed to iterate the one-step estimate. These
assumptions are not automatic consequences of real analyticity at a
single point.

\begin{theorem}
\label{thm:BEA-Lie}
Let $\Psi_h^{\mathrm{ext}}$ be the exactly symplectic extended
variational map of Proposition~\ref{prop:extended-map-Lie}, and suppose
that it has order $p\geq1$ for the flow of
$H^{\mathrm{ext}}=\nu+H(t,g,\mu)$. Let
$x_k=(\Psi_h^{\mathrm{ext}})^k(x_0)$, and assume that the portion of the
numerical trajectory under consideration admits a uniformly controlled
analytic family of charts in the sense of
Definition~\ref{def:uniform-analytic-family}.

Let
\[
\overline H^{\mathrm{ext}}
=
H^{\mathrm{ext}}
+
\sum_{r\geq p}h^rH_r
\]
be the formal Hamiltonian furnished by
Theorem~\ref{thm:affine-modified-Hamiltonian}, with its additive
constants fixed consistently on every connected component of its
domain. Assume that, in each chart $\chi_k$, the functions
$H_r\circ\chi_k^{-1}$ extend holomorphically to $K_{\rho_1}$ and satisfy
the uniform bounds
\begin{equation}
\label{eq:uniform-global-modified-Hamiltonian-bounds}
\sup_k
\left\|
H_r\circ\chi_k^{-1}
\right\|_{\rho_1}
\leq
A_HB_H^r r!,
\qquad r\geq p,
\end{equation}
for constants $A_H,B_H>0$ independent of $k$ and $r$.

Then there exist constants $h_*,c_0,c,C>0$, independent of $k$ and
$h$, such that $B_Hc_0<1$ and, for $0<h<h_*$, the truncation index
$N(h):=\lfloor c_0/h\rfloor$ defines the truncated modified
Hamiltonian
\begin{equation}
\label{eq:truncated-modified-Hamiltonian}
\overline H_{N(h)}^{\mathrm{ext}}
:=
\nu+H(t,g,\mu)
+
\sum_{r=p}^{N(h)}h^rH_r(t,g,\mu).
\end{equation}

Assume, moreover, that for every step under
consideration, the time-$h$ trajectory of
$X_{\overline H_{N(h)}^{\mathrm{ext}}}$
starting at $x_k$ exists and remains in
$\chi_k^{-1}(K_{\rho_2}^{\mathbb R})$. Its time-$h$ flow satisfies the uniform one-step estimate
\begin{equation}
\label{eq:exponential-one-step-defect}
d_{\chi_k}\left(
\Psi_h^{\mathrm{ext}}(x_k),
\varphi_h^{\overline H_{N(h)}^{\mathrm{ext}}}(x_k)
\right)
\leq
Ch\exp(-c/h).
\end{equation}

If the uniform analytic and domain-containment assumptions remain valid
for $0\leq kh\leq\exp(c/(2h))$, then
\begin{equation}
\label{eq:modified-energy-exp-long}
\left|
\overline H_{N(h)}^{\mathrm{ext}}(x_k)
-
\overline H_{N(h)}^{\mathrm{ext}}(x_0)
\right|
\leq
C\exp\bigl(-c/(2h)\bigr).
\end{equation}
Over the same interval, the original extended Hamiltonian satisfies
\begin{equation}
\label{eq:original-H-ext-error}
\left|
H^{\mathrm{ext}}(x_k)
-
H^{\mathrm{ext}}(x_0)
\right|
\leq
Ch^p.
\end{equation}
\end{theorem}

\begin{proof}
By Theorem~\ref{thm:affine-modified-Hamiltonian},
the formal modified vector field associated with
$\Psi_h^{\mathrm{ext}}$ is the Hamiltonian vector
field of $\overline H^{\mathrm{ext}}$. Since the method has order $p$, the
coefficients of orders $1,\dots,p-1$ in
$\overline H^{\mathrm{ext}}-H^{\mathrm{ext}}$ vanish.

Apply Lemma~\ref{lem:analytic-interpolation} in the chart $\chi_k$.
The bounds in Definition~\ref{def:uniform-analytic-family} are
independent of $k$, so the constants in the lemma can also be chosen
independently of $k$. The truncation of the formal vector field is the
Hamiltonian vector field of
$\overline H_{N(h)}^{\mathrm{ext}}$. Therefore,
\eqref{eq:analytic-interpolation-defect} gives
\eqref{eq:exponential-one-step-defect}.

For each $j$, let
$y_{j+1}:=
\varphi_h^{\overline H_{N(h)}^{\mathrm{ext}}}(x_j)$.
The modified interpolating segment is assumed to remain in
$\chi_j^{-1}(K_{\rho_2}^{\mathbb R})$. Since
$\overline H_{N(h)}^{\mathrm{ext}}$ is constant along its own
Hamiltonian flow,
\[
\overline H_{N(h)}^{\mathrm{ext}}(y_{j+1})
=
\overline H_{N(h)}^{\mathrm{ext}}(x_j).
\]

We next establish a uniform derivative bound for the
truncated Hamiltonian. Set
$\rho_*:=(\rho_1+\rho_2)/2$.
By Cauchy's estimate between $K_{\rho_*}$
and $K_{\rho_1}$, assumption
\eqref{eq:uniform-global-modified-Hamiltonian-bounds}
gives a constant $C_1>0$ such that
\[
\left\|
d\bigl(H_r\circ\chi_j^{-1}\bigr)
\right\|_{\rho_*}
\leq C_1 A_H B_H^r r!,
\qquad r\geq p.
\]
Moreover, the normalized Hamiltonian bound in
Definition~\ref{def:uniform-analytic-family},
together with Cauchy's estimate between
$K_{\rho_*}$ and $K_\rho$, gives
\[
\sup_j
\left\|
d\bigl(H^{\mathrm{ext}}\circ\chi_j^{-1}\bigr)
\right\|_{\rho_*}
\leq C_0.
\]

Choose $c_0$ smaller if necessary so that
$B_Hc_0<q<1$. For sufficiently small $h$,
the ratio between consecutive terms in the
derivative estimate satisfies
$hB_H(r+1)\leq q$ for $p\leq r<N(h)$.
Hence the finite sum is uniformly bounded,
and
\[
\sup_j
\left\|
d\bigl(
\overline H_{N(h)}^{\mathrm{ext}}
\circ\chi_j^{-1}
\bigr)
\right\|_{\rho_*}
\leq C_2
\]
for a constant $C_2$ independent of $j$ and $h$.

Both $x_{j+1}$ and $y_{j+1}$ belong to
$\chi_j^{-1}(K_{\rho_2}^{\mathbb R})$.
By \eqref{eq:exponential-one-step-defect},
\[
\left\|
\chi_j(x_{j+1})-\chi_j(y_{j+1})
\right\|_\infty
\leq Ch\exp(-c/h).
\]
For sufficiently small $h$, the right-hand
side is smaller than $\rho_*-\rho_2$.
Consequently, the line segment joining
$\chi_j(x_{j+1})$ and $\chi_j(y_{j+1})$
lies in $K_{\rho_*}^{\mathbb R}$.
The mean-value estimate therefore gives
\[
\left|
\overline H_{N(h)}^{\mathrm{ext}}(x_{j+1})
-
\overline H_{N(h)}^{\mathrm{ext}}(x_j)
\right|
\leq
Ch\exp(-c/h).
\]

Summing this inequality from $j=0$ to $j=k-1$ yields
\[
\left|
\overline H_{N(h)}^{\mathrm{ext}}(x_k)
-
\overline H_{N(h)}^{\mathrm{ext}}(x_0)
\right|
\leq
Ckh\exp(-c/h).
\]
If $kh\leq\exp(c/(2h))$, this proves
\eqref{eq:modified-energy-exp-long}.

It remains to compare the truncated modified Hamiltonian with the
original one. From
\eqref{eq:uniform-global-modified-Hamiltonian-bounds},
\[
\left|
\overline H_{N(h)}^{\mathrm{ext}}(x)
-
H^{\mathrm{ext}}(x)
\right|
\leq
A_H
\sum_{r=p}^{N(h)}
h^rB_H^r r!.
\]
Because $hB_H(r+1)\leq q<1$ throughout the retained range, the sum is
bounded by a constant multiple of its first term. Consequently,
\[
\left|
\overline H_{N(h)}^{\mathrm{ext}}(x)
-
H^{\mathrm{ext}}(x)
\right|
\leq Ch^p
\]
uniformly on the controlled domains. Applying this estimate at $x_0$
and $x_k$, and using
\eqref{eq:modified-energy-exp-long}, gives
\[
\left|
H^{\mathrm{ext}}(x_k)
-
H^{\mathrm{ext}}(x_0)
\right|
\leq
2Ch^p+C\exp\bigl(-c/(2h)\bigr).
\]
For sufficiently small $h$, the exponential term is bounded by a
constant multiple of $h^p$, which proves
\eqref{eq:original-H-ext-error}.
\end{proof}

Estimate~\eqref{eq:original-H-ext-error} concerns the autonomous
extended Hamiltonian
$H^{\mathrm{ext}}=\nu+H(t,g,\mu)$, not the physical time-dependent
Hamiltonian $H$ alone. In a non-autonomous system, $H$ is not
generally conserved. Its change is balanced by the evolution of the
conjugate momentum $\nu$. The conclusion of
Theorem~\ref{thm:BEA-Lie} is therefore
\[
\nu_k+H(t_k,g_k,\mu_k)
=
\nu_0+H(t_0,g_0,\mu_0)+O(h^p)
\]
over the stated exponentially long interval. Obvioursly, it does not imply
$H(t_k,g_k,\mu_k)=H(t_0,g_0,\mu_0)+O(h^p)$.

The exponentially long interval in
Theorem~\ref{thm:BEA-Lie} is conditional on the numerical trajectory
and the modified interpolating segments remaining in domains satisfying
the same analytic bounds. Real analyticity of the Lagrangian at each
individual point does not imply this uniform property. In particular,
the theorem does not exclude escape from every compact region, approach
to the boundary of a retraction chart, or growth of the derivatives of
the Hamiltonian along the trajectory.

\begin{corollary}
\label{cor:BEA-plus-casimirs}
Assume the hypotheses of Theorem~\ref{thm:BEA-Lie}. Suppose, in
addition, that $G$ is connected, that the discrete Lagrangian is
invariant under the diagonal left action of $G$ on $G\times G$, and
that the reduced body momentum satisfies
\begin{equation}
\label{eq:coadjoint-update-section-four}
\mu_{k+1}
=
\operatorname{Ad}_{F_k^{-1}}^*\mu_k,
\qquad
F_k=g_k^{-1}g_{k+1}.
\end{equation}
Then every Casimir of the associated Lie-Poisson structure is preserved
exactly:
\begin{equation}
\label{eq:exact-discrete-Casimir}
C(\mu_k)=C(\mu_0),
\qquad k\geq0.
\end{equation}
Independently, the truncated modified extended Hamiltonian satisfies
the long-time estimate
\eqref{eq:modified-energy-exp-long}.
\end{corollary}

\begin{proof}
Let $C\in C^\infty(\mathfrak g^*)$ be a Casimir of the Lie-Poisson
bracket associated with the coadjoint convention fixed in
Section~\ref{Sec2-Lie}. A Casimir is constant on every connected
symplectic leaf of the Lie-Poisson manifold. Since $G$ is connected, each coadjoint orbit is connected and is
precisely a symplectic leaf of the Lie-Poisson structure on
$\mathfrak g^*$.

Equation~\eqref{eq:coadjoint-update-section-four} shows that
$\mu_{k+1}$ lies on the coadjoint orbit through $\mu_k$. Therefore,
$C(\mu_{k+1})=C(\mu_k)$. Iteration proves
\eqref{eq:exact-discrete-Casimir}. The estimate for the truncated
modified Hamiltonian follows separately from
Theorem~\ref{thm:BEA-Lie}.
\end{proof}

The two conclusions of Corollary~\ref{cor:BEA-plus-casimirs} have
complementary geometric origins. Exact Casimir preservation follows
from coadjoint transport, independently of backward error analysis,
whereas near-conservation of the modified extended Hamiltonian relies
on analytic interpolation and optimal truncation. In the presence of
left invariance, these structures can be related more closely: the
discrete Noether conservation law constrains not only the numerical
flow but also its modified vector field. We show next that the
Hamiltonian constructed in Theorem~\ref{thm:affine-modified-Hamiltonian}
can be chosen to descend to the Lie-Poisson variables. This is a
formal reduction statement and does not require the analyticity
assumptions of Theorem~\ref{thm:BEA-Lie}.

\begin{theorem}
\label{thm:reduced-modified-Hamiltonian}
Assume the hypotheses of
Theorem~\ref{thm:affine-modified-Hamiltonian} on an open set
$\mathcal U\subset\overline M$ whose nonempty fibers under
$\rho|_{\mathcal U}$ are connected. Suppose that $G$ is connected,
that $\mathcal U$ is invariant under the cotangent-lifted left action
of $G$, and that the continuous and discrete Lagrangians are
left-invariant. Then the modified Hamiltonian can be chosen as
\[
\overline H^{\mathrm{ext}}(t,\nu,g,\mu;h)
=\nu+H(t,\mu)+\sum_{r=p}^{\infty}h^rH_r(t,\mu).
\]
Consequently, its formal modified flow descends, on the corresponding
reduced domain, to the time-dependent Lie-Poisson equation
\[
\dot\mu
=\operatorname{ad}_{\delta\overline H_h/\delta\mu}^{*}\mu,
\qquad
\overline H_h(t,\mu)
:=H(t,\mu)+\sum_{r=p}^{\infty}h^rH_r(t,\mu).
\]
All these identities are understood coefficientwise in $h$.
\end{theorem}

\begin{proof}
For $\zeta\in\mathfrak g$, let
\[
J_\zeta(t,\nu,g,\mu)
:=\langle\operatorname{Ad}_g^*\mu,\zeta\rangle
\]
be the momentum-map component for the cotangent-lifted left action.
Left invariance and the discrete Noether theorem imply
$(\Psi_h^{\mathrm{ext}})^*J_\zeta=J_\zeta$; the update of $\nu$
does not affect this identity. Hence
$((\Psi_h^{\mathrm{ext}})^*-\operatorname{id})^mJ_\zeta=0$
for every $m\ge1$. The formal logarithm
\eqref{eq:formal-logarithm} therefore gives
$X_h[J_\zeta]=0$, where $X_h$ is the modified vector field.
By Theorem~\ref{thm:affine-modified-Hamiltonian},
$X_h=X_{\overline H^{\mathrm{ext}}}$, and our Poisson convention yields
\[
\{J_\zeta,\overline H^{\mathrm{ext}}\}_{\mathrm{ext}}=0.
\]
In particular, $\{J_\zeta,H_r\}_{\mathrm{ext}}=0$ for every
$r\ge p$. The Hamiltonian vector fields of the functions
$J_\zeta$ generate the cotangent-lifted left action, so these
identities imply infinitesimal invariance of every $H_r$.
Since $G$ is connected and $\mathcal U$ is $G$-invariant, each
coefficient is constant along the full $G$-orbits in $\mathcal U$.
In left-trivialized variables this action is
$(t,\nu,g,\mu)\mapsto(t,\nu,ag,\mu)$.

By \eqref{eq:modified-H-reduced}, each $H_r$ is also independent of
$\nu$; connectedness of the nonempty $\rho$-fibers ensures the
corresponding global descent. Thus $H_r$ depends only on $(t,\mu)$.
The continuous Hamiltonian $H$ has the same independence from $g$
because the continuous Lagrangian is left-invariant.
The left-trivialized Hamilton equations then give the stated
Lie-Poisson equation with the coadjoint convention fixed in
Section~\ref{Sec2-Lie}.
\end{proof}

The preceding theorem strengthens the geometric interpretation of the
modified dynamics: under symmetry, both the numerical momentum and its
formal modified evolution are governed by the same Lie-Poisson
geometry. It does not, however, imply conservation of the reduced
Hamiltonian, which remains explicitly time-dependent. To express the
long-time result in physically meaningful variables, we now return to
the general, not necessarily invariant, discrete Lagrangian and use
the conjugate-time update of Proposition~\ref{prop:extended-map-Lie}. The
extended-Hamiltonian bounds can then be rewritten as discrete energy
balances without the auxiliary variable $\nu$.

\begin{corollary}
\label{cor:long-time-discrete-energy-balance}
Assume the hypotheses of Theorem~\ref{thm:BEA-Lie}, and define
\[
\overline H_{N(h)}(t,g,\mu)
:=H(t,g,\mu)+\sum_{r=p}^{N(h)}h^rH_r(t,g,\mu).
\]
For every numerical trajectory and every index $k$ in the
exponentially long interval of that theorem,
\[
\begin{aligned}
\Bigg|&\overline H_{N(h)}(t_k,g_k,\mu_k)
-\overline H_{N(h)}(t_0,g_0,\mu_0)+\sum_{j=0}^{k-1}\partial_tL_d(t_j,g_j,g_{j+1};h)
\Bigg|\le C\exp\bigl(-c/(2h)\bigr),
\end{aligned}
\]
and
\[
\begin{aligned}
\Bigg|&H(t_k,g_k,\mu_k)-H(t_0,g_0,\mu_0)+\sum_{j=0}^{k-1}\partial_tL_d(t_j,g_j,g_{j+1};h)
\Bigg|\le Ch^p.
\end{aligned}
\]
Here $\partial_tL_d$ is taken with the two group endpoints and the
step size held fixed.
\end{corollary}

\begin{proof}
By the conjugate-momentum update of Proposition~\ref{prop:extended-map-Lie},
\[
\nu_k-\nu_0
=\sum_{j=0}^{k-1}\partial_tL_d(t_j,g_j,g_{j+1};h).
\]
Using
$\overline H_{N(h)}^{\mathrm{ext}}=\nu+\overline H_{N(h)}$,
the first assertion is exactly
\eqref{eq:modified-energy-exp-long}.
The identity $H^{\mathrm{ext}}=\nu+H$ similarly turns
\eqref{eq:original-H-ext-error} into the second assertion.
\end{proof}

The preceding corollary is a long-time analogue of the continuous
energy balance: the accumulated discrete time dependence compensates
the variation of the physical Hamiltonian, up to the stated defect.
Its validity depends essentially on the uniform control of the
modified Hamiltonian coefficients in
\eqref{eq:uniform-global-modified-Hamiltonian-bounds}. This
condition is not merely a technical restatement of the factorial
bounds for the modified vector field. The free-particle model used in the counterexample preceding
\eqref{eq:smooth-variational-contact} also provides a simple way
to distinguish the two obstructions.

\begin{remark}
The free-particle example of subsecction \ref{sec4.2} shows that a
pointwise $O(h^3)$ action defect need not yield a second-order
variational map: the perturbation
$\varepsilon h^3\sin((q_1-q_0)/h^2)$ has insufficiently controlled
endpoint derivatives. The obstruction addressed by
\eqref{eq:uniform-global-modified-Hamiltonian-bounds} is different.

Indeed, let
\[
L_{d,E}(q_0,q_1;h)=\frac{(q_1-q_0)^2}{2h},
\qquad
L_d(t,q_0,q_1;h)=L_{d,E}(q_0,q_1;h)+a\,t\,h^3,
\quad a\ne0.
\]
The projected method on $(q,\mu)$ is exact:
$q_{k+1}=q_k+h\mu_k$ and $\mu_{k+1}=\mu_k$. Nevertheless,
the extended update of Proposition~\ref{prop:extended-map-Lie} gives
$\nu_{k+1}=\nu_k+a h^3$. The corresponding extended map is the
exact time-$h$ flow of the modified Hamiltonian
\[
\overline H^{\mathrm{ext}}
=\nu+\frac12\mu^2-a h^2t.
\]
Thus the modified coefficient $H_2(t,\mu)=-at$ has uniformly bounded
derivatives but unbounded values on the full time axis. Moreover,
\[
H^{\mathrm{ext}}(x_k)-H^{\mathrm{ext}}(x_0)=a k h^3,
\]
which need not be $O(h^2)$ when $kh$ grows exponentially in $h^{-1}$.
The one-step map itself is analytic and near the identity; the
failure is the absence of a uniform bound of the form
\eqref{eq:uniform-global-modified-Hamiltonian-bounds}.

\end{remark}

\subsection{Application to the time-dependent rigid body}
\label{subsec:BEA-rigid-body}

The preceding results now admit a concrete interpretation for the
time-dependent rigid body. Its left-invariant discrete Lagrangian
gives exact coadjoint transport, as in
Corollary~\ref{cor:BEA-plus-casimirs}, whereas
Corollary~\ref{cor:long-time-discrete-energy-balance} concerns the
compensation between the physical Hamiltonian and the accumulated
discrete time dependence. These conclusions are logically distinct.
We first write the extended update explicitly and then verify the
order of the complete method, including the conjugate momentum $\nu$,
before invoking Theorem~\ref{thm:affine-modified-Hamiltonian} with
$p=2$ for the midpoint-Cayley scheme.

Let $G=SO(3)$. We identify $\mathfrak{so}(3)$ with $\mathbb R^3$ through
the hat map $x\mapsto\widehat x$. Using the Euclidean pairing, we also identify
$\mathfrak{so}(3)^*$ with $\mathbb R^3$. With the coadjoint convention
fixed in Section~\ref{Sec2-Lie}, one has
$\operatorname{ad}_x^*\mu=\mu\times x$. Let $\mathbb I:\mathbb R^3\to\mathbb R^3$ be symmetric and positive
definite. The continuous Lagrangian is
\[
L(t,R,\xi)
=
\frac12e^{-\kappa t}
\left\langle\xi,\mathbb I\xi\right\rangle.
\]
Its left-trivialized Legendre transformation is
$\mu=e^{-\kappa t}\mathbb I\xi$. Therefore,
$\xi=e^{\kappa t}\mathbb I^{-1}\mu$, and the corresponding Hamiltonian
is
\begin{equation*}
\label{eq:time-dependent-rigid-Hamiltonian}
H(t,R,\mu)
=
\frac12e^{\kappa t}
\left\langle\mu,\mathbb I^{-1}\mu\right\rangle.
\end{equation*}
The Hamiltonian is independent of $R$, but it depends explicitly on
time.

The associated extended Hamiltonian on
$T^*(\mathbb R\times SO(3))$ is
\begin{equation}
\label{eq:extended-rigid-Hamiltonian}
H^{\mathrm{ext}}(t,\nu,R,\mu)
=
\nu+
\frac12e^{\kappa t}
\left\langle\mu,\mathbb I^{-1}\mu\right\rangle.
\end{equation}
Using the convention
$\iota_{X_{H^{\mathrm{ext}}}}\omega_{\mathrm{ext}}
=dH^{\mathrm{ext}}$, its Hamilton equations are
\begin{equation*}
\label{eq:extended-rigid-body-equations}
\begin{aligned}
\dot t
&=1,
&
\dot\nu
&=
-\frac{\kappa}{2}e^{\kappa t}
\left\langle\mu,\mathbb I^{-1}\mu\right\rangle,
\\
\dot R
&=
R\,\widehat{e^{\kappa t}\mathbb I^{-1}\mu},
&
\dot\mu
&=
\mu\times
\bigl(e^{\kappa t}\mathbb I^{-1}\mu\bigr).
\end{aligned}
\end{equation*}
In particular, the signs in the equations for $\nu$ and $\mu$ follow,
respectively, from
$\dot\nu=-\partial H/\partial t$ and
$\dot\mu=\operatorname{ad}_{\delta H/\delta\mu}^*\mu$, since
$\mathbf{D}_RH=0$.

The extended Hamiltonian is conserved along the continuous extended
flow:
\begin{equation*}
\label{eq:rigid-extended-energy-balance}
\frac{d}{dt}
H^{\mathrm{ext}}(t,\nu,R,\mu)
=
0.
\end{equation*}
Equivalently, the change in the physical Hamiltonian is exactly
balanced by the conjugate momentum $\nu$. This conservation law is
different from the physical energy balance
\eqref{eq:continuous-energy-law} expressed in the velocity variables.

For the Cayley midpoint variational integrator, consider the discrete
Lagrangian
\begin{equation*}
\label{eq:rigid-midpoint-discrete-Lagrangian}
L_d(t_k,R_k,R_{k+1};h)
=
\frac h2e^{-\kappa t_{k+\frac12}}
\left\langle\xi_k,\mathbb I\xi_k\right\rangle,
\qquad
t_{k+\frac12}=t_k+\frac h2,
\end{equation*}
where
$\xi_k=h^{-1}\operatorname{cay}^{-1}
(R_k^{-1}R_{k+1})$.
This is the same discrete Lagrangian as
\eqref{eq:Ld-rigid}, now written with its initial-time argument
explicitly displayed.

When differentiating $L_d$ with respect to $t_k$, the endpoints
$R_k$ and $R_{k+1}$, and hence $\xi_k$, are held fixed. Consequently,
$\partial L_d/\partial t=-\kappa L_d$. The exactly symplectic extension
of Proposition~\ref{prop:extended-map-Lie} therefore gives
\begin{equation}
\label{eq:rigid-nu-update}
\nu_{k+1}
=
\nu_k-\kappa
L_d(t_k,R_k,R_{k+1};h).
\end{equation}
Together with \eqref{eq:rigid-complete-update}, the complete extended
update is thus obtained by solving
\eqref{eq:rigid-implicit-xi}, setting
$F_k=\operatorname{cay}(h\widehat\xi_k)$, and updating
\[
t_{k+1}=t_k+h,\qquad
R_{k+1}=R_kF_k,\qquad
\mu_{k+1}=F_k^{-1}\mu_k,
\]
followed by \eqref{eq:rigid-nu-update}. The equations determining
$(R_{k+1},\mu_{k+1})$ are independent of $\nu_k$, consistently with
the $\rho$-projectability of the extended map.

The second-order estimate in
Corollary~\ref{cor:midpoint-cayley-order} concerns the projected
variables $(g,\mu)$. To apply
Theorem~\ref{thm:affine-modified-Hamiltonian} with $p=2$, we must
also control the error in $\nu$.
For the midpoint-Cayley discrete Lagrangian, the expansions in the
proof of Corollary~\ref{cor:midpoint-cayley-order} have smooth,
uniformly bounded remainders including their first derivatives in
$t$ on the compact families considered. Therefore the pulled-back
action defect satisfies
\[
\Delta_d(t,h,g,\xi)=h^3\rho(t,h,g,\xi),
\qquad \partial_t\Delta_d=O(h^3).
\]
Let $g_1=\Gamma_{t,h}(g,\xi)$ and
$\delta L_h:=L_d-L_{d,E}$. Differentiation at fixed
$(g,\xi,h)$ gives
\[
\partial_t\Delta_d
=\partial_t\delta L_h(t,g,g_1)
+\bigl\langle
D_2\delta L_h(t,g,g_1),\,
\partial_t\Gamma_{t,h}(g,\xi)
\bigr\rangle.
\]
The endpoint derivative is $O(h^2)$ by the proof of
Theorem~\ref{thm:variational-order}. Since the first-order
reconstruction term is $hT_eL_g\xi$, with $(g,\xi)$ fixed,
$\partial_t\Gamma_{t,h}=O(h^2)$. Hence
\[
\partial_t\delta L_h
\bigl(t,g,\Gamma_{t,h}(g,\xi)\bigr)=O(h^3).
\]
For the exact discrete Lagrangian, differentiation with respect
to the initial time at fixed endpoints, followed by stationarity
of the exact action, gives
\[
\partial_tL_{d,E}(t,g_0,g_1;h)
=\int_0^h\partial_tL(t+s,g(s),\xi(s))\,ds.
\]
The integral is evaluated on the exact connecting trajectory,
with the endpoint values held fixed in the partial derivative.
When pulled back by the exact endpoint map, its derivative with
respect to the initial velocity is $O(h)$ uniformly on the
compact families under consideration.

For the same initial $(g,\mu)$, let $\xi$ parametrize the exact
endpoint and $\bar\xi$ the endpoint determined by the numerical
negative Legendre transformation. By
Theorem~\ref{thm:variational-order},
$\bar\xi-\xi=O(h^2)$. Thus
\[
\begin{aligned}
&\partial_tL_d\bigl(t,g,\Gamma_{t,h}(g,\bar\xi);h\bigr)
 -\partial_tL_{d,E}\bigl(t,g,\Gamma_{t,h}(g,\xi);h\bigr)\\
&\quad=\partial_t\delta L_h
\bigl(t,g,\Gamma_{t,h}(g,\bar\xi)\bigr)
+O(h)\,O(h^2)=O(h^3).
\end{aligned}
\]
The exact and discrete extended maps update $\nu$ by these
respective partial time derivatives. Their local error in the
conjugate momentum is therefore $O(h^3)$. Together with
Corollary~\ref{cor:midpoint-cayley-order}, this proves that the
complete midpoint-Cayley extension has order two.

Since the complete midpoint-Cayley extension
has order two, Theorem~\ref{thm:affine-modified-Hamiltonian}
applies with $p=2$. Moreover, the continuous and
discrete Lagrangians are left-invariant.
Consequently, on the invariant domains of
Theorem~\ref{thm:reduced-modified-Hamiltonian},
the modified Hamiltonian takes the stronger
reduced form
\[
\overline H^{\mathrm{ext}}
=
\nu+\frac12 e^{\kappa t}
\langle\mu,\mathbb I^{-1}\mu\rangle
+\sum_{r=2}^{\infty}h^rH_r(t,\mu).
\]
Thus every correction coefficient is independent
of both the attitude $R$ and the conjugate
momentum $\nu$.
Under the uniform analyticity,
domain-containment, and modified-Hamiltonian
coefficient bounds of
Theorem~\ref{thm:BEA-Lie}, its optimally
truncated series is conserved up to an
exponentially small variation over
exponentially long intervals. The original extended Hamiltonian
\eqref{eq:extended-rigid-Hamiltonian} consequently differs from its
initial value by $\mathcal O(h^2)$ over the same interval.

For the numerical solution, define
\[
\mathcal E_{\mathrm{ext}}^k
:=
H^{\mathrm{ext}}(t_k,\nu_k,R_k,\mu_k)
=
\nu_k+
\frac12e^{\kappa t_k}
\left\langle\mu_k,\mathbb I^{-1}\mu_k\right\rangle.
\]
Since the initial physical Hamiltonian
$H_0:=H(t_0,R_0,\mu_0)$ is strictly positive for $\mu_0\neq0$, we
measure the signed normalized variation of the extended Hamiltonian by
\begin{equation}
\label{eq:rigid-extended-Hamiltonian-error}
\varepsilon_{\mathrm{ext}}^k(h)
:=
\frac{
\mathcal E_{\mathrm{ext}}^k
-
\mathcal E_{\mathrm{ext}}^0
}{
H_0
},
\qquad
D_{\mathrm{ext}}(h;T)
:=
\max_{0\leq kh\leq T}
\left|
\varepsilon_{\mathrm{ext}}^k(h)
\right|.
\end{equation}
Thus $\varepsilon_{\mathrm{ext}}^k$ retains the sign of the numerical
variation, whereas $D_{\mathrm{ext}}(h;T)$ records its maximum absolute value over the integration interval.

By the conjugate-momentum update
\eqref{eq:rigid-nu-update}, we have
\[
\nu_k-\nu_0
=
\sum_{j=0}^{k-1}
\partial_t L_d(t_j,R_j,R_{j+1};h).
\]
Consequently, $D_{\mathrm{ext}}(h;T)$ also
measures the maximum normalized defect
in the discrete energy balance of
Corollary~\ref{cor:long-time-discrete-energy-balance}.

For every fixed $T$ contained in the interval covered by
Theorem~\ref{thm:BEA-Lie}, equation
\eqref{eq:original-H-ext-error} predicts
$D_{\mathrm{ext}}(h;T)=\mathcal O(h^2)$. For two consecutive step
sizes, the corresponding observed order is
\[
p_{\mathrm{ext}}(h_j)
:=
\frac{
\log\bigl(
D_{\mathrm{ext}}(h_j;T)/
D_{\mathrm{ext}}(h_{j+1};T)
\bigr)
}{
\log(h_j/h_{j+1})
}.
\]
An observed order approaching two is
consistent with the predicted second-order
behavior.

\paragraph{Fixed-time convergence of the conjugate momentum}

The second-order result established above concerns the entire extended
map, not only its projection onto $(R,\mu)$. We therefore assess the convergence of
the additional variable $\nu$ directly,
rather than inferring its accuracy from
the extended-Hamiltonian variation. For the continuous rigid body, define
\[
K(\mu):=\frac12\langle\mu,\mathbb I^{-1}\mu\rangle.
\]
The reduced equations imply
\[
\frac{d}{dt}K(\mu(t))
=e^{\kappa t}\left\langle
\mathbb I^{-1}\mu,\,
\mu\times\mathbb I^{-1}\mu\right\rangle=0.
\]
Integrating $\dot\nu=-\kappa e^{\kappa t}K(\mu)$ hence provides the
exact reference
\begin{equation}
\label{eq:rigid-exact-nu}
\nu(t)=\nu_0+K(\mu_0)
\bigl(e^{\kappa t_0}-e^{\kappa t}\bigr).
\end{equation}

We repeat the fixed-final-time convergence experiment of
Example~\ref{ex:numerical-variational-order}, with
$\kappa=-0.4$, $T=10$, and the same inertia, initial data, and step
sizes. Starting from $\nu_0=0$, we advance the conjugate momentum
using \eqref{eq:rigid-nu-update} at each step and define
\[
e_\nu(h):=|\nu_N-\nu(T)|,\qquad Nh=T.
\]
For the chosen initial data, $K(\mu_0)=1.45$ and
$\nu(10)=1.423442323611336$. The errors and observed orders are
reported in Table~\ref{tab:rigid-nu-convergence}.
\begin{table}[htbp]
\centering
\begin{tabular}{ccc}
\hline
$h$ & $e_\nu(h)$ & Observed order \\
\hline
$0.080$ & $4.4381\times10^{-3}$ & --- \\
$0.040$ & $1.1192\times10^{-3}$ & $1.9876$ \\
$0.020$ & $2.8040\times10^{-4}$ & $1.9969$ \\
$0.010$ & $7.0138\times10^{-5}$ & $1.9992$ \\
$0.005$ & $1.7537\times10^{-5}$ & $1.9998$ \\
\hline
\end{tabular}
\caption{Fixed-time convergence of the conjugate momentum $\nu$ for
the complete midpoint-Cayley extension at $T=10$. The exact
reference is given by \eqref{eq:rigid-exact-nu}.}
\label{tab:rigid-nu-convergence}
\end{table}

The observed order approaches two, consistently with the local
$O(h^3)$ estimate for $\nu$ derived above. Together with the configuration and momentum convergence
reported in Example~\ref{ex:numerical-variational-order}, this
provides a numerical check of the second-order behavior of every
component of the extended map over a fixed time interval.

\paragraph{Long-time extended-phase-space experiment}

We use the same physical parameters and initial data as in the long-time
experiment of Section~\ref{Sec:DiscreteNoether}:
$\kappa=-10^{-3}$, $\mathbb I=\operatorname{diag}(2,1,0.5)$,
$t_0=0$, $R_0=I$, and
$\xi(0)=(0.8,1.2,-0.6)^\top$. As explained in
Example~\ref{ex:numerical-variational-order}, the value
$\kappa=-10^{-3}$ keeps the time-dependent energy scale numerically
resolvable over the interval $[0,10^4]$ and is not a stability
restriction of the integrator. The initial body momentum is
$\mu_0=\mathbb I\xi(0)$, and we set $\nu_0=0$. The choice of $\nu_0$
only fixes the initial value of the extended Hamiltonian and does not
affect the projected evolution of $(R_k,\mu_k)$.

The first experiment uses $h=5\times10^{-2}$ and
$T=10^4$, corresponding to $2\times10^5$ time steps. At each
iteration, the conjugate momentum is advanced using
\eqref{eq:rigid-nu-update}, together with the configuration and
momentum updates of Section~3. To examine the step-size dependence of
the extended-Hamiltonian variation, the experiment is repeated over
the same interval with
$h\in\{2\times10^{-1},10^{-1},5\times10^{-2},
2.5\times10^{-2}\}$.

\begin{figure}[h!]
\centering
\begin{minipage}{0.48\textwidth}
    \centering
    \includegraphics[width=\textwidth]
    {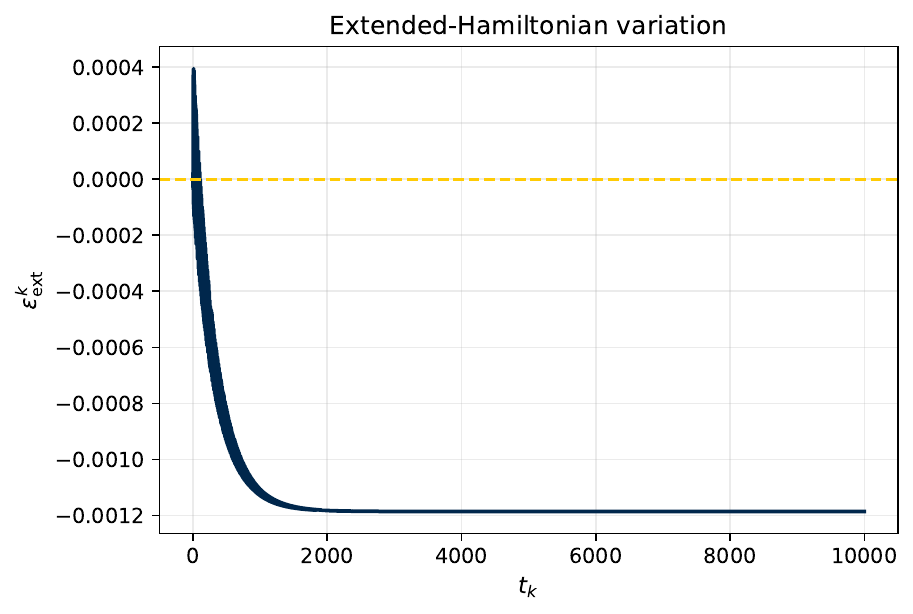}
    \small (a) Long-time extended-Hamiltonian variation
\end{minipage}\hfill
\begin{minipage}{0.48\textwidth}
    \centering
    \includegraphics[width=\textwidth]
    {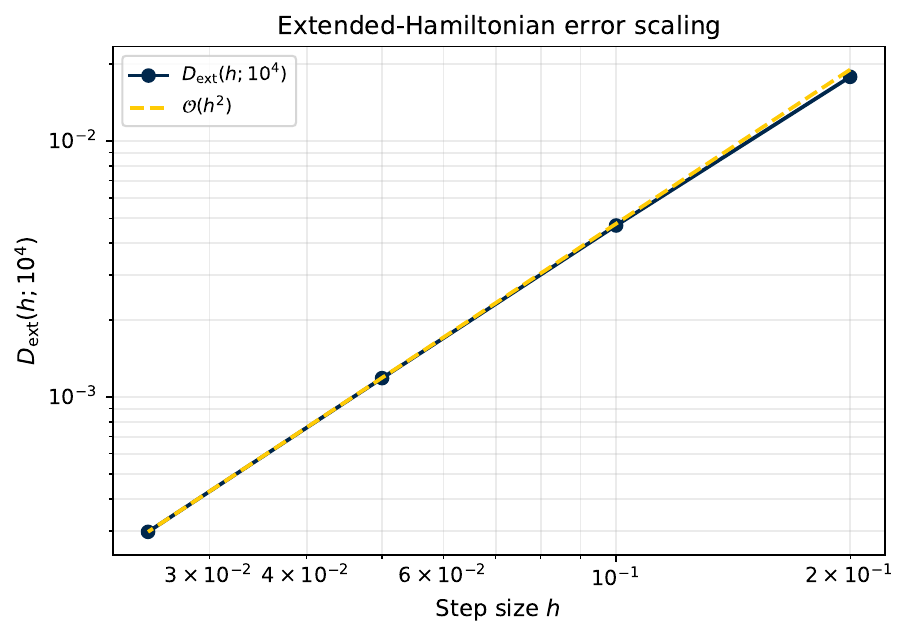}
    \small (b) Step-size dependence of the maximum variation
\end{minipage}
\caption{Extended-phase-space numerical results for the Cayley midpoint
variational integrator with $\kappa=-10^{-3}$ and $T=10^4$.
(a) Signed normalized variation
$\varepsilon_{\mathrm{ext}}^k$ for $h=5\times10^{-2}$; the dashed
yellow line represents exact conservation of $H^{\mathrm{ext}}$.
(b) Maximum normalized variation
$D_{\mathrm{ext}}(h;10^4)$ for the four tested step sizes, together
with a reference line proportional to $h^2$.}
\label{fig:rigid-extended-Hamiltonian}
\end{figure}

For $h=5\times10^{-2}$, the signed variation initially oscillates
about zero and subsequently approaches a bounded negative value. Over
the complete interval, it satisfies
$D_{\mathrm{ext}}(5\times10^{-2};10^4)
=1.186\times10^{-3}$, with no secular growth visible over the
$2\times10^5$ time steps.

The maximum normalized variations obtained in the step-size study are
\[
\begin{array}{c|c|c}
h & D_{\mathrm{ext}}(h;10^4) & p_{\mathrm{ext}}\\
\hline
0.200 & 1.784\times10^{-2} & \text{---}\\
0.100 & 4.682\times10^{-3} & 1.9299\\
0.050 & 1.186\times10^{-3} & 1.9813\\
0.025 & 2.974\times10^{-4} & 1.9953
\end{array}
\]
The observed orders approach two as the step size decreases,
consistent with the $\mathcal O(h^2)$ estimate for the original
extended Hamiltonian in
\eqref{eq:original-H-ext-error}. The maximum implicit residual over
all four experiments is $4.628\times10^{-15}$.

The behavior of the trace in
Figure~\ref{fig:rigid-extended-Hamiltonian}(a) also reflects the
time scale of the underlying non-autonomous dynamics. Indeed, the
reduced equation can be written as
\[
\dot\mu=e^{\kappa t}\mu\times\mathbb I^{-1}\mu.
\]
With $t_0=0$, the change of time variable
\begin{equation}
\label{eq:rigid-effective-time}
s(t)=\int_0^t e^{\kappa u}\,du
=\frac{1-e^{\kappa t}}{-\kappa}
\end{equation}
transforms it into the autonomous free rigid-body equation in the
effective time $s$. For $\kappa<0$, one has
$s(t)\to 1/|\kappa|$ as $t\to\infty$; the continuous momentum and
attitude therefore evolve progressively more slowly. The exact
conjugate-momentum update also decays in magnitude, since
$\dot\nu=-\kappa e^{\kappa t}K(\mu_0)$.
Consequently, the flattening of the numerical variation is partly
associated with the saturation of the underlying dynamics, rather
than constituting a direct verification of exponentially accurate
Hamiltonian interpolation. The refinement study in
Figure~\ref{fig:rigid-extended-Hamiltonian}(b) instead gives numerical
evidence consistent with the $O(h^2)$ estimate
\eqref{eq:original-H-ext-error} for the tested interval.

This experiment evaluates the original extended Hamiltonian
\eqref{eq:extended-rigid-Hamiltonian}; it does not numerically evaluate
the optimally truncated modified Hamiltonian
\eqref{eq:truncated-modified-Hamiltonian}. Such an evaluation would
require the explicit computation of the correction coefficients
$H_2,H_3,\ldots$, which are not used by the present implementation.
Consequently, the exponentially small estimate
\eqref{eq:modified-energy-exp-long} remains a theoretical result,
whereas Figure~\ref{fig:rigid-extended-Hamiltonian}
provides numerical evidence for the $\mathcal O(h^2)$ behavior of
the original extended Hamiltonian. The experiment does not, by itself,
verify the uniform analytic-domain hypotheses or the exponential-time
assertion of Theorem~\ref{thm:BEA-Lie}.

\paragraph{Periodic modulation without asymptotic slowing}

To distinguish bounded long-time energy variation from the slowing
induced by the factor $e^{\kappa t}$, we consider a second
left-invariant rigid-body Lagrangian with periodic scalar modulation,
\begin{equation}
\label{eq:periodic-rigid-body-L}
L^{\mathrm{per}}(t,R,\xi)
=\frac{1}{2a(t)}\langle\xi,\mathbb I\xi\rangle,
\qquad
a(t)=1+\varepsilon\sin(\omega t),
\quad 0<\varepsilon<1.
\end{equation}
This model retains the same $SO(3)$ geometry and discrete Noether
symmetry, but the time dependence does not decay. Its Hamiltonian is
$H^{\mathrm{per}}(t,\mu)=a(t)K(\mu)$, and the continuous reduced
equations read
\[
\dot\mu=a(t)\mu\times\mathbb I^{-1}\mu,
\qquad
\dot\nu=-a'(t)K(\mu).
\]
The effective time is now
\[
s(t)=\int_0^t a(u)\,du
=t+\frac{\varepsilon}{\omega}
\bigl(1-\cos(\omega t)\bigr),
\]
which grows without bound. Moreover, $K(\mu(t))=K(\mu_0)$ and
$\nu(t)=\nu_0+K(\mu_0)[a(0)-a(t)]$.

We use the same midpoint-Cayley construction, now with the
midpoint factor $e^{-\kappa t_{k+1/2}}$ in
\eqref{eq:Ld-rigid} replaced by $a(t_{k+1/2})^{-1}$.
In particular, the implicit body-velocity equation is obtained from
\eqref{eq:rigid-implicit-xi} by the same substitution, while the
exactly symplectic extension has the update
\begin{equation}
\label{eq:periodic-rigid-nu-update}
\nu_{k+1}=\nu_k-
\frac{a'(t_{k+1/2})}{a(t_{k+1/2})}
L_d^{\mathrm{per}}(t_k,R_k,R_{k+1};h).
\end{equation}
The configuration and momentum updates remain
$R_{k+1}=R_kF_k$ and $\mu_{k+1}=F_k^{-1}\mu_k$.

\begin{figure}[htbp]
\centering
\begin{minipage}{0.48\textwidth}
\centering
\includegraphics[width=\textwidth]{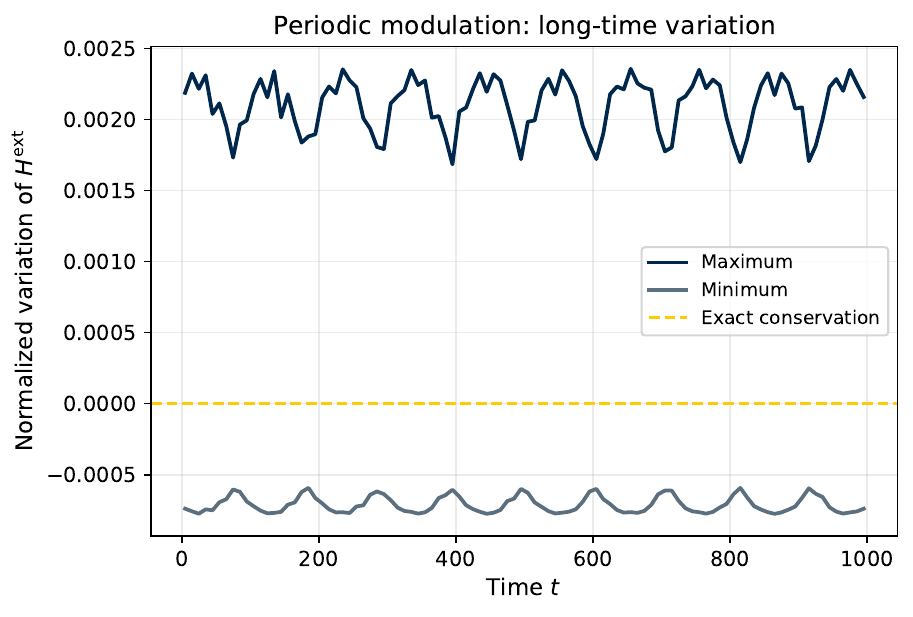}
\small (a) Envelope over $T=1000$ for $h=0.05$
\end{minipage}\hfill
\begin{minipage}{0.48\textwidth}
\centering
\includegraphics[width=\textwidth]{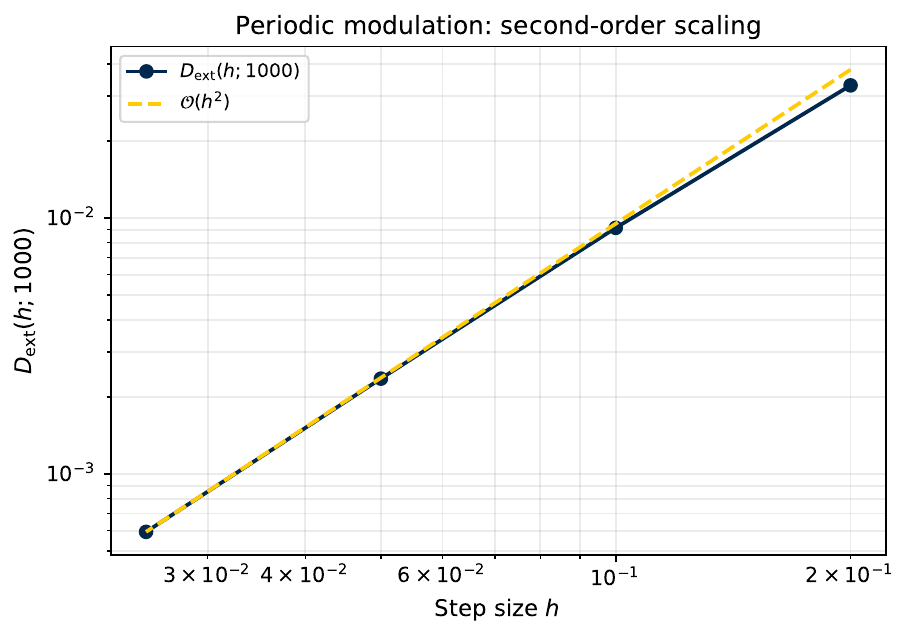}
\small (b) Maximum variation versus step size
\end{minipage}
\caption{Extended-Hamiltonian numerical results for the periodically
modulated rigid body \eqref{eq:periodic-rigid-body-L}, with
$\varepsilon=0.3$, $\omega=1$, and $T=1000$.
(a) Maximum and minimum values of the signed
normalized extended-Hamiltonian variation,
computed over successive time intervals
of length $10$.
(b) Maximum normalized variation for four step sizes and a reference
line proportional to $h^2$.}
\label{fig:periodic-rigid-extended-Hamiltonian}
\end{figure}

For the experiment we take $\varepsilon=0.3$, $\omega=1$,
$T=1000$, $R_0=I$, $\xi(0)=(0.8,1.2,-0.6)^\top$,
$\mathbb I=\operatorname{diag}(2,1,0.5)$, $t_0=0$,
$\nu_0=0$, and
$h\in\{0.2,0.1,0.05,0.025\}$. The initial extended Hamiltonian
is $H_0^{\mathrm{ext}}=K(\mu_0)=1.45$.
We use the same normalized variation as in
\eqref{eq:rigid-extended-Hamiltonian-error}, with
$H^{\mathrm{ext}}$ replaced by
$H^{\mathrm{ext,per}}=\nu+a(t)K(\mu)$.
The results are displayed in
Figure~\ref{fig:periodic-rigid-extended-Hamiltonian}.

The maximum normalized variations are
\[
\begin{array}{c|c|c}
h & D_{\mathrm{ext}}(h;1000) & p_{\mathrm{ext}}\\
\hline
0.200 & 3.300\times10^{-2} & \text{---}\\
0.100 & 9.157\times10^{-3} & 1.8494\\
0.050 & 2.358\times10^{-3} & 1.9571\\
0.025 & 5.941\times10^{-4} & 1.9889
\end{array}
\]
The maximum implicit residual in these four runs was below
$2.5\times10^{-13}$. In contrast with the first long-time experiment,
the effective time grows linearly, and the numerical variation
continues to oscillate throughout the integration interval. Its
maximum decreases with an observed order approaching two.
This complementary experiment therefore supports the $O(h^2)$
extended-energy scaling without relying on asymptotic slowing; it
still does not constitute a numerical verification of the
exponentially small modified-Hamiltonian defect
\eqref{eq:modified-energy-exp-long}.

The Casimir and spatial-momentum errors reported in
Figure~\ref{fig:so3_diagnostics} remain at the level of floating-point
roundoff. Those exact-arithmetic preservation properties arise from
coadjoint transport and discrete Noether symmetry, independently of
backward error analysis. Taken together,
Figures~\ref{fig:so3_diagnostics},
\ref{fig:rigid-extended-Hamiltonian}, and
\ref{fig:periodic-rigid-extended-Hamiltonian} illustrate the
complementary structural and asymptotic mechanisms analyzed in this
paper: exact symmetry-induced conservation laws and bounded
near-conservation of the extended Hamiltonian under two distinct
forms of explicit time dependence.

\section{Application: time-dependent mechanical consensus on
\texorpdfstring{$SO(3)$}{SO(3)}}
\label{sec:app-consensus}

Consensus and synchronization on Lie groups have been studied from several
complementary perspectives. Potential-based coordination and symmetry
reduction for networks of rigid bodies were developed in
\cite{HanssmannLeonardSmith2006,NairLeonard2007}, while the stabilization
of synchronized mechanical motions by controlled-Lagrangian and energy
methods was investigated in \cite{NairLeonard2008}. Consensus optimization
on Riemannian manifolds was studied in
\cite{SarletteSepulchre2009}, including the geometric and topological
limitations that distinguish nonlinear configuration spaces from Euclidean
ones. More recently, second-order Laplacian flows and a unified framework
connecting Euclidean consensus, synchronization on the circle, and
consensus on Lie groups admitting bi-invariant metrics were developed in
\cite{SeshanBanavarMahindrakar2022,
ChandrasekaranBanavarMahindrakar2025,
ChandrasekaranBanavarMahindrakarBook2025}.

The continuous geometry of relative-configuration potentials and the
stability of autonomous rigid body networks are therefore not new. The
purpose of this section is different. We consider a time-dependent mechanical network, establish a local convergence result
under explicit assumptions on its time-varying coupling weights, and apply
the discrete variational and backward error constructions of
Sections~\ref{Sec:DiscreteNoether} and~\ref{sec:BEA}. The resulting method
preserves the momentum associated with the diagonal symmetry exactly and
admits an affine modified Hamiltonian on the extended phase space.

The exponential time dependence introduced below damps the angular
velocities themselves. Consequently, the limiting motion considered here is
static attitude consensus, rather than synchronization along a common
nonzero rotational motion of the type studied in
\cite{NairLeonard2007,NairLeonard2008}. We therefore use the term
{consensus} for convergence to a common attitude with vanishing
angular velocities.

\subsection{Continuous time-dependent attitude-consensus dynamics}
\label{subsec:time-dependent-SO3-network}

Let $\mathcal G=(\mathcal V,\mathcal E)$ be a fixed, connected and
undirected graph, with vertex set $\mathcal V=\{1,\ldots,N\}$. Its edges, $\mathcal E,$
represent the pairs of agents that exchange relative-attitude information.
For every $(i,j)\in\mathcal E$, let
$w_{ij}:[0,\infty)\to(0,\infty)$ be a continuously differentiable effective
coupling weight satisfying $w_{ij}=w_{ji}$. 

The configuration space is the product Lie group $SO(3)^N$. We write
$\mathbf R=(R_1,\ldots,R_N)$ and define the body angular velocities
$\Omega_i\in\mathbb R^3$ by
\begin{equation}
\label{eq:network-reconstruction}
\dot R_i=R_i\widehat\Omega_i,
\qquad
i=1,\ldots,N.
\end{equation}
For each edge, the relative attitude is $R_j^\top R_i$. We use the function $V:SO(3)\to\mathbb R$ given by
$V(R)=\operatorname{tr}(I-R)$. If $R$ is a rotation through the angle
$\vartheta\in[0,\pi]$, then
$V(R)=2(1-\cos\vartheta)$. Hence $V\geq0$, $V(R)=V(R^{-1})$, and
$V(R)=0$ if and only if $R=I$.

The time-dependent network potential is
\begin{equation}
\label{eq:SO3-network-potential}
U(t,\mathbf R)
:=
\sum_{(i,j)\in\mathcal E}
w_{ij}(t)
V(R_j ^T R_i)
=\sum_{(i,j)\in\mathcal E}
w_{ij}(t)
\operatorname{tr}(I-R_j^\top R_i),
\end{equation}
where each undirected edge is included once. Its zero set is the consensus
manifold
$\mathcal C=\{(R,\ldots,R)\in SO(3)^N:R\in SO(3)\}$. Indeed, connectivity
of $\mathcal G$ implies that $U(t,\mathbf R)=0$ if and only if
$R_i=R_j$ along every edge, and hence if and only if
$\mathbf R\in\mathcal C$.

Let $\mathbb I_i:\mathbb R^3\to\mathbb R^3$ be a symmetric and positive
definite transformation. For a constant $\kappa<0$, consider
\begin{equation}
\label{eq:network-Lagrangian}
L(t,\mathbf R,\boldsymbol\Omega)
=
e^{-\kappa t}
\left[
\frac12
\sum_{i=1}^N
\Omega_i^\top\mathbb I_i\Omega_i
-
U(t,\mathbf R)
\right].
\end{equation}
This is a time-dependent mechanical Lagrangian of Caldirola-Kanai type  \cite{Caldirola1941,Kanai1948}.
It is not the uncontrolled Lagrangian of a physical vehicle network. Its
Euler-Lagrange equations describe a mechanical attitude network with
isotropic viscous damping and time-dependent internal coupling.

For later use, define
\begin{equation}
\label{eq:SO3-network-gradient}
\mathcal T_i(t,\mathbf R)
:=
\mathbf{D}_{R_i}U(t,\mathbf R)
=
\sum_{j:\,(i,j)\in\mathcal E}
w_{ij}(t)
\left(
R_j^\top R_i-R_i^\top R_j
\right)^\vee .
\end{equation}
Here $\mathbf{D}_{R_i}U\in\mathfrak{so}(3)^*\simeq\mathbb R^3$ is characterized
by differentiating along $R_i\exp(s\widehat\eta_i)$. Formula
\eqref{eq:SO3-network-gradient} follows from
$-\operatorname{tr}(R\widehat\eta)
=\langle(R-R^\top)^\vee,\eta\rangle$.

\begin{proposition}\label{prop:SO3-network-EL}Given the Lagrangian defined in \eqref{eq:network-Lagrangian}, its Euler-Lagrange equations are defined by \eqref{eq:network-reconstruction} and

\begin{equation}
\label{eq:SO3-network-continuous}
\mathbb I_i\dot\Omega_i
+
\Omega_i\times(\mathbb I_i\Omega_i)
-
\kappa\mathbb I_i\Omega_i
+
\mathcal T_i(t,\mathbf R)
=
0,
\qquad
i=1,\ldots,N.
\end{equation}

\end{proposition}

\begin{proof}
The fiber and left-trivialized configuration derivatives of
\eqref{eq:network-Lagrangian} are
$\partial L/\partial\Omega_i
=e^{-\kappa t}\mathbb I_i\Omega_i$ and
$\mathbf{D}_{R_i}L=-e^{-\kappa t}\mathcal T_i$. Substitution into the left-trivialized Euler-Lagrange equations
\eqref{eq:EL-Lie} gives
$\mathbb I_i\dot\Omega_i
-\operatorname{ad}_{\Omega_i}^*(\mathbb I_i\Omega_i)
-\kappa\mathbb I_i\Omega_i+\mathcal T_i=0$.
Under the convention fixed in Section~\ref{Sec2-Lie},
$\operatorname{ad}_{\Omega_i}^*\mu_i=\mu_i\times\Omega_i$.
Consequently,
$-\operatorname{ad}_{\Omega_i}^*(\mathbb I_i\Omega_i)
=\Omega_i\times(\mathbb I_i\Omega_i)$, which yields
\eqref{eq:SO3-network-continuous}.
\end{proof}

In the Abelian case, the corresponding construction reduces to a
second-order damped Laplacian system with time-dependent edge weights. On
$S^1$, the same relative potential produces sinusoidal coupling. These are
the Euclidean and circular models unified with Lie group consensus in
\cite{ChandrasekaranBanavarMahindrakar2025,
ChandrasekaranBanavarMahindrakarBook2025}. Equation
\eqref{eq:SO3-network-continuous} is their mechanical attitude counterpart,
with the additional non-autonomous dependence carried by the weights
$w_{ij}(t),$ for $(i,j)\in\mathcal E$.

We next turn to the stability of the continuous network.
Pointwise connectivity of a time-dependent graph is not, by itself,
sufficient to establish stability of a second-order mechanical network.
Rapid variation of a positive stiffness can produce parametric excitation.
We therefore state explicit conditions under which the time dependence
cannot inject mechanical energy.

\begin{assumption}
\label{ass:monotone-network-weights}
For every $(i,j)\in\mathcal E$, there exist constants
$w_{ij}^-,w_{ij}^+>0$ such that
\begin{equation}
\label{eq:weight-assumptions}
0<w_{ij}^-\leq w_{ij}(t)\leq w_{ij}^+,
\qquad
\dot w_{ij}(t)\leq0,
\qquad
t\geq0.
\end{equation}
\end{assumption}

The lower bound preserves the coercivity of the connected network, while
monotonicity prevents the time-dependent coupling from increasing the
potential energy. Each weight consequently has a positive limit
$w_{ij}^\infty:=\displaystyle{\lim_{t\to\infty}}w_{ij}(t)$.

Define the mechanical energy, without the Caldirola-Kanai multiplier, by
\begin{equation}
\label{eq:network-mechanical-energy}
\mathbb E(t,\mathbf R,\boldsymbol\Omega)
=
\frac12
\sum_{i=1}^N
\Omega_i^\top\mathbb I_i\Omega_i
+
U(t,\mathbf R).
\end{equation}

\begin{lemma}
\label{lem:network-energy-dissipation} Under Assumption~\ref{ass:monotone-network-weights},  every solution satisfying
\eqref{eq:network-reconstruction} and \eqref{eq:SO3-network-continuous} verifies that 
\begin{equation}
\label{eq:network-energy-dissipation}
\frac{d}{dt}\mathbb E
=
\kappa
\sum_{i=1}^N
\Omega_i^\top\mathbb I_i\Omega_i
+
\sum_{(i,j)\in\mathcal E}
\dot w_{ij}(t)
\operatorname{tr}(I-R_j^\top R_i) \leq0.
\end{equation}

\end{lemma}

\begin{proof}
Pairing \eqref{eq:SO3-network-continuous} with $\Omega_i$ and using
$\Omega_i^\top[\Omega_i\times(\mathbb I_i\Omega_i)]=0$ gives the
derivative of the kinetic energy. The configuration-dependent part of
$dU/dt$ cancels the sum of the terms
$\langle\mathcal T_i,\Omega_i\rangle$. The remaining explicit time
derivative is the second term on the right-hand side of
\eqref{eq:network-energy-dissipation}. Both terms are non-positive because
$\kappa<0$, $\dot w_{ij}\leq0$, and
$\operatorname{tr}(I-R_j^\top R_i)\geq0$.
\end{proof}

Let
$\mathcal A:=\mathcal C\times\{0\}\subset
SO(3)^N\times(\mathbb R^3)^N$ denote the static consensus set. Stability
of $\mathcal A$ is understood as stability of a set and therefore modulo
the free diagonal left action of $SO(3)$.

\begin{theorem}
\label{thm:time-dependent-local-consensus}
Suppose that $\mathcal G$ is connected,
Assumption~\ref{ass:monotone-network-weights} holds, and $\kappa<0$.
Then $\mathcal A$ is locally asymptotically stable for
\eqref{eq:network-reconstruction} and \eqref{eq:SO3-network-continuous}.
More precisely, there exists a diagonal-action-invariant neighborhood
$\mathcal N$ of $\mathcal A$ such that every solution starting
sufficiently close to $\mathcal A$ remains in $\mathcal N$ and satisfies
\begin{equation}
\label{eq:network-consensus-limit}
\operatorname{dist}(\mathbf R(t),\mathcal C)\rightarrow0,
\qquad
\Omega_i(t)\rightarrow0,
\qquad
i=1,\ldots,N.
\end{equation}
\end{theorem}

\begin{proof}
Choose local coordinates transverse to $\mathcal C$ by writing
$R_i=R\exp(\widehat\eta_i)$, where the common infinitesimal rotation is
removed, for instance, by imposing $\sum_i\eta_i=0$. Since
$\operatorname{tr}(I-\exp(\widehat x))
=\|x\|^2+\mathcal O(\|x\|^4)$, the quadratic part of the potential is
$\displaystyle{\sum_{(i,j)\in\mathcal E}}w_{ij}(t)\|\eta_i-\eta_j\|^2$.
Its kernel consists precisely of diagonal variations because
$\mathcal G$ is connected. The lower bounds in
\eqref{eq:weight-assumptions} imply uniform positive definiteness on the
transverse subspace. Thus, on a sufficiently small invariant neighborhood
of $\mathcal C$, the potential is uniformly equivalent to the squared
distance from $\mathcal C$. Together with positive definiteness of the
inertia operators, this shows that \eqref{eq:network-mechanical-energy} is
locally and uniformly positive definite with respect to $\mathcal A$.

By Lemma~\ref{lem:network-energy-dissipation},
the mechanical energy is non-increasing along
solutions. Together with the uniform local
coercivity established above, this implies
that solutions with sufficiently small initial
energy remain in a normal neighborhood of
$\mathcal C$ in which the consensus manifold
is the only critical manifold of the potential.
This proves local stability of $\mathcal A$.

Let $c=-\kappa>0$. Integration of
\eqref{eq:network-energy-dissipation} gives
$\int_0^\infty\sum_i
\Omega_i^\top\mathbb I_i\Omega_i\,dt<\infty$.
The solution remains in a compact subset of the chosen neighborhood, and
the right-hand side of \eqref{eq:SO3-network-continuous} is bounded.
Hence every $\dot\Omega_i$ is bounded and
$\sum_i\|\Omega_i\|^2$ is uniformly continuous. Barbalat's lemma then
implies $\Omega_i(t)\to0$ for every $i$.

It remains to identify the limiting configurations. By monotonicity and
the positive lower bounds, $w_{ij}(t)\to w_{ij}^\infty>0$. The vector
field therefore converges uniformly on compact subsets to the autonomous
vector field obtained by replacing $w_{ij}(t)$ with
$w_{ij}^\infty$. Let a point of the $\omega$-limit set be obtained along
a sequence $t_n\to\infty$. The shifted trajectories based at $t_n$ admit
a subsequence converging on compact time intervals to a trajectory of the
limiting autonomous system. Since $\Omega_i(t)\to0$, this limiting
trajectory has identically vanishing angular velocities. Its invariance
then requires $\mathbf{D}_{R_i}U_\infty(\mathbf R)=0$ for every $i$, where
$U_\infty$ is the potential with weights $w_{ij}^\infty$. Within the
chosen normal neighborhood, the only such critical configurations belong
to $\mathcal C$. Consequently, the $\omega$-limit set is contained in
$\mathcal A$, which proves \eqref{eq:network-consensus-limit}.
\end{proof}

The monotonicity condition is sufficient rather than necessary. It selects
a nontrivial class of time-dependent networks for which the
continuous convergence argument remains rigorous. General oscillatory or
switching weights require a different non-autonomous stability analysis
and are not covered by Theorem~\ref{thm:time-dependent-local-consensus}.

We now complement this convergence result with the continuous conservation law associated with the diagonal symmetry. The potential \eqref{eq:SO3-network-potential} depends only on relative
attitudes and is therefore invariant under the diagonal left action
$P\cdot(R_1,\ldots,R_N)=(PR_1,\ldots,PR_N)$. This invariance is unaffected
by the explicit time dependence of the weights.

The momentum conjugate to $\Omega_i$ is
$\mu_i=e^{-\kappa t}\mathbb I_i\Omega_i$. Applying
Theorem~\ref{thm:noether-Lie} to the diagonal action gives conservation of
the total spatial momentum
\begin{equation}
\label{eq:network-continuous-momentum}
\Pi(t)
=
\sum_{i=1}^N
R_i(t)\mu_i(t)
=
e^{-\kappa t}
\sum_{i=1}^N
R_i(t)\mathbb I_i\Omega_i(t).
\end{equation}
Thus the conserved Noether quantity is the exponentially weighted spatial
momentum, not the unweighted physical angular momentum. This distinction
is a consequence of the time-dependent Lagrangian
\eqref{eq:network-Lagrangian}.

Conservation of \eqref{eq:network-continuous-momentum} is compatible with
Theorem~\ref{thm:time-dependent-local-consensus} since
$e^{-\kappa t}$ grows when $\kappa<0$, a nonzero constant value of
$\Pi(t)$ can coexist with exponentially decaying physical angular
velocities.

\subsection{Variational integrator and geometric properties}\label{sec6.2}

For each agent and time step, define
$F_{i,k}=R_{i,k}^\top R_{i,k+1}
=\operatorname{cay}(h\widehat\Omega_{i,k})$ and
$R_{i,k+\frac12}
=R_{i,k}\operatorname{cay}(\frac h2\widehat\Omega_{i,k})$.
Let $\mathbf F_k=(F_{1,k},\ldots,F_{N,k})$ and
$\mathbf R_{k+\frac12}
=(R_{1,k+\frac12},\ldots,R_{N,k+\frac12})$.

The midpoint discrete Lagrangian is
\begin{equation}
\label{eq:SO3-network-discrete-L}
L_d^k(\mathbf R_k,\mathbf R_{k+1})
=
h e^{-\kappa t_{k+\frac12}}
\left[
\frac12
\sum_{i=1}^N
\Omega_{i,k}^\top\mathbb I_i\Omega_{i,k}
-
U(t_{k+\frac12},\mathbf R_{k+\frac12})
\right].
\end{equation}
Its mixed representation is
$\widehat L_d^k(\mathbf R_k,\mathbf F_k)
=L_d^k(\mathbf R_k,\mathbf R_k\mathbf F_k)$, and its Lie algebra
representation is
$\ell_d^k(\mathbf R_k,\boldsymbol\Omega_k)
=\widehat L_d^k(\mathbf R_k,
\operatorname{cay}(h\widehat\Omega_{1,k}),\ldots,
\operatorname{cay}(h\widehat\Omega_{N,k}))$.

Define the constitutive post-momentum
\begin{equation*}
\label{eq:network-post-momentum}
\mathcal{M}_{i,k}(R_k,\Omega_k)
:=
\frac{1}{h}
\left(
\mathrm d^L\operatorname{cay}^{-1}_{F_{i,k}}
\right)^*
\frac{\partial\ell_d^k}{\partial\Omega_i}(R_k,\Omega_k).
\end{equation*}  
For numerical implementation, the derivatives
entering the implicit equations can be expressed
in terms of the midpoint interaction torques.
Set
\[
C_{i,k}
:=
\operatorname{cay}
\left(\frac h2\widehat\Omega_{i,k}\right),
\qquad
\mathcal T_{i,k}^{m}
:=
\mathcal T_i
(t_{k+\frac12},\mathbf R_{k+\frac12}),
\]
and introduce the linear map
\[
Q_{i,k}(a)
:=
I_3+\frac a2\widehat\Omega_{i,k}
+\frac{a^2}{4}
\Omega_{i,k}\Omega_{i,k}^{\top}.
\]

Using the left-trivialized differential of
the inverse Cayley map in~\eqref{eq:dlog-cayley},
the constitutive post-momentum becomes
\begin{equation}
\label{eq:network-explicit-post-momentum}
\mathcal M_{i,k}
=
e^{-\kappa t_{k+\frac12}}
Q_{i,k}(h)^\top
\left[
\mathbb I_i\Omega_{i,k}
-\frac h2
Q_{i,k}(h/2)^{-\top}
\mathcal T_{i,k}^{m}
\right].
\end{equation}
Moreover, differentiation of the mixed
discrete Lagrangian while holding
$\mathbf F_k$ fixed gives
\begin{equation}
\label{eq:network-explicit-configuration-derivative}
\mathbf D_{R_i}\widehat L_d^k
=
-h e^{-\kappa t_{k+\frac12}}
C_{i,k}\mathcal T_{i,k}^{m}.
\end{equation}

Indeed, a variation of $\Omega_{i,k}$
induces the left-trivialized variation
\[
\delta R_{i,k+\frac12}
\longleftrightarrow
\frac h2
Q_{i,k}(h/2)^{-1}
\delta\Omega_{i,k},
\]
whereas a variation
$\delta R_{i,k}=R_{i,k}\widehat\eta_i$
with $\mathbf F_k$ fixed induces
$\delta R_{i,k+\frac12}
\leftrightarrow
C_{i,k}^{\top}\eta_i$.
These identities yield the two expressions
above and provide an explicit evaluation
of the coupled discrete momentum equations.

\begin{theorem}
\label{thm:network-integrator}
Assume that the discrete Legendre transformation associated with
\eqref{eq:SO3-network-discrete-L} is locally invertible. Given
$(\mathbf {R}_k,\mu_k)$, determine
$\boldsymbol\Omega_k$ from the coupled implicit equations
\begin{equation}
\label{eq:network-implicit-equations}
\mu_{i,k}
=
F_{i,k}
\mathcal M_{i,k}(\mathbf R_k,\boldsymbol\Omega_k)
-
\mathbf{D}_{R_i}\widehat L_d^k(\mathbf R_k,\mathbf F_k),
\qquad
i=1,\ldots,N.
\end{equation}
The phase-space update is then
\begin{equation}
\label{eq:SO3-network-update}
R_{i,k+1}=R_{i,k}F_{i,k},
\qquad
\mu_{i,k+1}
=
\mathcal M_{i,k}(\mathbf R_k,\boldsymbol\Omega_k).
\end{equation}
This is the left-trivialized phase-space map generated by the discrete
variational principle on $SO(3)^N$.
\end{theorem}

\begin{proof}
Apply Corollary~\ref{prop:cayley-momentum-update}
to the product group $SO(3)^N$. For each component, that corollary gives 
$\mu_{i,k+1}=\mathcal M_{i,k}$ and
$\mu_{i,k}
=\operatorname{Ad}_{F_{i,k}}^*\mu_{i,k+1}
-\mathbf{D}_{R_i}\widehat L_d^k$. Under the Euclidean identification
$\mathfrak{so}(3)^*\simeq\mathbb R^3$,
$\operatorname{Ad}_{F_{i,k}}^*\mu=F_{i,k}\mu$. Substitution of the
constitutive post-momentum therefore gives
\eqref{eq:network-implicit-equations}. Once these equations have been
solved, reconstruction and the post-momentum relation give
\eqref{eq:SO3-network-update}. Local invertibility of the discrete
Legendre transformation guarantees that the coupled implicit equations
determine a local phase-space map.
\end{proof}

The dependence of the potential on all midpoint attitudes makes
\eqref{eq:network-implicit-equations} a coupled nonlinear system. It
cannot, in general, be replaced by $N$ independent copies of the
single rigid body equation.

Having constructed the discrete phase-space map, we next
examine the preservation of the diagonal symmetry at the discrete level.
For every $P\in SO(3)$, the increments and midpoint relative attitudes
satisfy $(PR_{i,k})^\top(PR_{i,k+1})=F_{i,k}$ and
$(PR_{j,k+\frac12})^\top(PR_{i,k+\frac12})
=R_{j,k+\frac12}^\top R_{i,k+\frac12}$. Hence
\eqref{eq:SO3-network-discrete-L} is invariant under the diagonal left
action at every time step.

\begin{theorem}
\label{thm:network-noether}
Along the discrete flow
\eqref{eq:network-implicit-equations} and \eqref{eq:SO3-network-update},
the total spatial momentum
\begin{equation}
\label{eq:SO3-network-total-momentum}
\Pi_k
=
\sum_{i=1}^N R_{i,k}\mu_{i,k}
\end{equation}
is preserved exactly in exact arithmetic
$\Pi_{k+1}=\Pi_k$.
\end{theorem}

\begin{proof}
For $\zeta\in\mathfrak{so}(3)$, the infinitesimal generator of the
diagonal left action is
$\zeta_{SO(3)^N}(\mathbf R)
=(T_e\mathcal R_{R_1}\zeta,\ldots,
T_e\mathcal R_{R_N}\zeta)$.
Its body representation at $R_i$ is
$\operatorname{Ad}_{R_i^{-1}}\zeta$. Therefore, the component of the
discrete momentum map in the direction $\zeta$ is
$\sum_i\langle\mu_{i,k},
\operatorname{Ad}_{R_{i,k}^{-1}}\zeta\rangle
=\langle\sum_i\operatorname{Ad}_{R_{i,k}}^*\mu_{i,k},\zeta\rangle$.
Theorem~\ref{thm:discrete-noether} implies that this quantity is
independent of $k$. Since
$\operatorname{Ad}_{R_{i,k}}^*\mu_{i,k}=R_{i,k}\mu_{i,k}$ under the
Euclidean identification, the conserved momentum is precisely
\eqref{eq:SO3-network-total-momentum}.
\end{proof}

The individual body momenta are not separately conserved because the interaction
potential produces the force terms in
\eqref{eq:network-implicit-equations}. In particular, unlike the free rigid body,
the individual quadratic quantities $\|\mu_{i,k}\|^2$ are not invariants in
general. Indeed, \eqref{eq:network-implicit-equations} and
\eqref{eq:SO3-network-update} give
$\mu_{i,k+1}=F_{i,k}^{-1}
\bigl(\mu_{i,k}+\mathbf{D}_{R_i}\widehat L_d^k(\mathbf R_k,\mathbf F_k)\bigr)$,
so the orthogonality of $F_{i,k}$ does not imply
$\|\mu_{i,k+1}\|=\|\mu_{i,k}\|$ when the discrete interaction force is
nonzero. The discrete Lagrangian is invariant under the diagonal
$SO(3)$ action, rather than under independent left actions on every factor.
Consequently, only the total spatial momentum
$\Pi_k$ in \eqref{eq:SO3-network-total-momentum} is preserved. This conservation
law holds for arbitrary time dependence of the weights;
Assumption~\ref{ass:monotone-network-weights} is required for the continuous
convergence theorem, but not for the discrete Noether theorem.

We finally place the discrete network map in the extended
Hamiltonian framework of Section~\ref{sec:BEA} to describe its long-time behavior. The Hamiltonian associated with \eqref{eq:network-Lagrangian} is
\begin{equation*}
\label{eq:network-Hamiltonian}
H(t,\mathbf R,\boldsymbol\mu)
=
\frac12 e^{\kappa t}
\sum_{i=1}^N
\mu_i^\top\mathbb I_i^{-1}\mu_i
+
e^{-\kappa t}U(t,\mathbf R),
\end{equation*}
where $\mu_i=e^{-\kappa t}\mathbb I_i\Omega_i$. Its autonomous extension
is
\begin{equation*}
\label{eq:network-extended-Hamiltonian}
H^{\mathrm{ext}}(t,\nu,\mathbf R,\boldsymbol\mu)
=
\nu+H(t,\mathbf R,\boldsymbol\mu).
\end{equation*}

The exactly symplectic extension of the discrete map advances the
conjugate momentum according to
\begin{equation}
\label{eq:network-nu-update}
\nu_{k+1}
=
\nu_k+
\frac{\partial L_d}{\partial t}
(t_k,\mathbf R_k,\mathbf R_{k+1};h),
\end{equation}
where the endpoint attitudes are held fixed in the partial derivative.
For \eqref{eq:SO3-network-discrete-L}, the discrete velocities and
midpoint attitudes are then independent of $t_k$, and
\begin{equation*}
\label{eq:network-time-derivative-Ld}
\frac{\partial L_d}{\partial t}
=
h e^{-\kappa t_{k+\frac12}}
\left[
-\kappa
\left(
\frac12\sum_{i=1}^N
\Omega_{i,k}^\top\mathbb I_i\Omega_{i,k}
-
U(t_{k+\frac12},\mathbf R_{k+\frac12})
\right)
-
\frac{\partial U}{\partial t}
(t_{k+\frac12},\mathbf R_{k+\frac12})
\right].
\end{equation*}

The diagonal symmetry has a further consequence
for the modified dynamics. Unlike the
left-invariant rigid body of Subsection~5.5,
the interaction potential is not invariant
under independent left translations of the
$N$ attitudes. Thus the full reduction of
Theorem~\ref{thm:reduced-modified-Hamiltonian}
does not apply directly.

Nevertheless, the discrete Noether conservation
law constrains the formal modified Hamiltonian.
In particular, its dependence on the attitudes
can be expressed entirely through relative
configurations.

\begin{proposition}
\label{prop:network-diagonal-modified-Hamiltonian}
Assume the hypotheses of
Theorem~\ref{thm:affine-modified-Hamiltonian}
for the extended network map on an open set
$\mathcal U\subset\overline M$ invariant
under the cotangent-lifted diagonal action
of $SO(3)$, with connected nonempty fibers
of $\rho|_{\mathcal U}$.

Define the relative attitudes
\[
S_i:=R_1^\top R_i,
\qquad i=2,\ldots,N,
\]
and set
$\mathbf S:=(S_2,\ldots,S_N)$.

Then the formal modified Hamiltonian can be
chosen in the form
\begin{equation}
\label{eq:network-diagonal-modified-Hamiltonian}
\overline H^{\mathrm{ext}}
=
\nu+H(t,\mathbf S,\boldsymbol\mu)
+
\sum_{r=p}^{\infty}
h^r H_r(t,\mathbf S,\boldsymbol\mu).
\end{equation}
Here the functions on the right denote
their representations on the diagonal
quotient, identified by the section
$(I,S_2,\ldots,S_N)$.

In particular, all modified Hamiltonian
coefficients are invariant under a common
spatial rotation, although they may depend
on the relative attitudes. The associated
formal modified vector field descends to
the corresponding Poisson quotient.

This result is formal and does not require
the analytic assumptions of
Theorem~\ref{thm:BEA-Lie}.
\end{proposition}

\begin{proof}
For each $\zeta\in\mathfrak{so}(3)$,
consider the momentum-map component
\[
J_\zeta
=
\left\langle
\sum_{i=1}^{N}R_i\mu_i,\zeta
\right\rangle.
\]
By Theorem~\ref{thm:network-noether},
the extended variational map preserves
$J_\zeta$ exactly:
\[
(\Psi_h^{\mathrm{ext}})^*J_\zeta
=
J_\zeta.
\]
Consequently, the formal logarithm of
the pullback operator annihilates
$J_\zeta$. The formal modified vector
field therefore satisfies
\[
X_h[J_\zeta]=0.
\]

By Theorem~\ref{thm:affine-modified-Hamiltonian},
$X_h$ is Hamiltonian, with formal Hamiltonian
$\overline H^{\mathrm{ext}}$. Hence
\[
\{J_\zeta,\overline H^{\mathrm{ext}}\}_{\mathrm{ext}}
=0.
\]
Comparing coefficients of equal powers
of $h$ yields
\[
\{J_\zeta,H_r\}_{\mathrm{ext}}=0,
\qquad r\geq p.
\]

The Hamiltonian vector fields of the
components $J_\zeta$ generate the
cotangent-lifted diagonal action.
Since $SO(3)$ is connected and
$\mathcal U$ is invariant under this action,
each $H_r$ is constant along its orbits.

Under left trivialization, the diagonal
action is
\[
P\cdot
(t,\nu,\mathbf R,\boldsymbol\mu)
=
(t,\nu,P\mathbf R,\boldsymbol\mu),
\qquad P\in SO(3).
\]
Its configuration orbits are characterized
by the relative attitudes
$(R_1^\top R_2,\ldots,R_1^\top R_N)$,
while the body momenta remain unchanged.

Thus every $H_r$ depends on the
configuration only through $\mathbf S$.
Its independence of $\nu$ follows from
Theorem~\ref{thm:affine-modified-Hamiltonian}.
The continuous Hamiltonian has the same
diagonal invariance. This proves
\eqref{eq:network-diagonal-modified-Hamiltonian}.

Finally, the diagonal action is symplectic,
so its invariant functions define the
Poisson algebra of the quotient.
The invariant formal Hamiltonian therefore
generates a well-defined formal Hamiltonian
vector field on that quotient.
\end{proof}

\begin{corollary}\label{cor:network-BEA}
Assume that the midpoint-Cayley network
integrator is defined on common local
regularity domains satisfying the smoothness
hypotheses of Corollary~\ref{cor:midpoint-cayley-order}. Assume also the
invariant-domain hypotheses of
Proposition~\ref{prop:network-diagonal-modified-Hamiltonian}.
Then the complete extended variational map
has order two, and its formal modified
Hamiltonian has the reduced form
\eqref{eq:network-diagonal-modified-Hamiltonian}
with $p=2$.

Suppose, in addition, that the extended
numerical trajectory satisfies all the
hypotheses of Theorem~\ref{thm:BEA-Lie},
including the uniformly controlled analytic
family of charts, the coefficient bounds
\eqref{eq:uniform-global-modified-Hamiltonian-bounds},
and the required containment of the modified
interpolating trajectories. Then the optimally truncated modified
extended Hamiltonian satisfies
\eqref{eq:modified-energy-exp-long}, while
the original extended Hamiltonian satisfies
\eqref{eq:original-H-ext-error} with $p=2$,
over the exponentially long interval
specified in that theorem.
\end{corollary}

\begin{proof}
The midpoint-Cayley construction applies
to the product Lie group $SO(3)^N$.
By Corollary~\ref{cor:midpoint-cayley-order}, its projected phase-space
map has local error $O(h^3)$.

We must also consider the conjugate momentum
$\nu$. The argument used in Subsection~5.5
to establish the order of the complete
extended map applies on the present product
group. Indeed, smooth variational contact
and smooth dependence on the initial time give
\[
\Delta_d=h^3\rho,
\qquad
\partial_t\Delta_d=O(h^3).
\]
The endpoint estimates in the proof of
Theorem~4.5 yield
$D_2\delta L_h=O(h^2)$, while
$\partial_t\Gamma_{t,h}=O(h^2)$
at fixed initial configuration and velocity.
The chain rule therefore gives
\[
\partial_t\delta L_h
\bigl(t,\mathbf R,
\Gamma_{t,h}(\mathbf R,\boldsymbol\Omega)\bigr)
=O(h^3).
\]

Moreover, the difference between the
initial velocities associated with the
exact and numerical endpoints is $O(h^2)$.
The derivative with respect to the initial
velocity of the exact discrete
time derivative, pulled back by the exact
endpoint map, is $O(h)$.
Consequently, the local error in the
conjugate-momentum update is also $O(h^3)$.

Thus the complete extended map has order
two. Proposition~5 gives exact symplecticity
and $\rho$-projectability, while
Proposition~\ref{prop:network-diagonal-modified-Hamiltonian}
establishes the reduced form of its formal
modified Hamiltonian, with the expansion
starting at order $h^2$.

Finally, under the additional uniform
analyticity, coefficient, and
domain-containment hypotheses, the
conclusions follow directly from
Theorem~\ref{thm:BEA-Lie}.
\end{proof}

The exact conservation of \eqref{eq:SO3-network-total-momentum} and the
near-conservation of the modified extended Hamiltonian are logically
distinct. The former follows from diagonal symmetry and holds without
analyticity assumptions. The latter follows from symplectic backward
error analysis and requires the hypotheses stated in
Corollary~\ref{cor:network-BEA}.

\subsection{Numerical experiment}
\label{subsec:network-long-time-experiment}

The results of Subsection~\ref{sec6.2} establish
complementary properties of the discrete network
dynamics: exact preservation of the total spatial
momentum, second-order accuracy of the complete
extended variational map, and conditional long-time
estimates for the extended Hamiltonian.
We now examine these properties numerically
for a mechanical attitude-consensus system
on $SO(3)^N$.

We first assess the convergence of the complete
extended map $(\mathbf R,\boldsymbol\mu,\nu)$
against a high-accuracy reference solution on
a fixed time interval. We then investigate
attitude consensus, angular-velocity decay,
preservation of the total spatial momentum,
and retention of the $SO(3)^N$ constraints
over a large number of time steps.
Finally, we study the long-time variation
of the original extended Hamiltonian and
its dependence on the step size.

The numerical experiments illustrate these
properties without asserting that the uniform
analyticity, coefficient bounds, and
domain-containment assumptions required for
the exponential-time estimates of
Corollary~\ref{cor:network-BEA}
have been independently verified.

We consider four agents connected through the cycle graph with edge set
$\mathcal E=\{(1,2),(2,3),(3,4),(4,1)\}$. The inertia tensors are
$\mathbb I_1=\operatorname{diag}(2,1,0.5)$,
$\mathbb I_2=\operatorname{diag}(1.8,1.1,0.6)$,
$\mathbb I_3=\operatorname{diag}(2.2,0.9,0.55)$, and
$\mathbb I_4=\operatorname{diag}(1.9,1.2,0.45)$. For each $(i,j)\in\mathcal E$, we choose the time-dependent coupling
weight
$w_{ij}(t)=w_\infty+(w_{ij}^0-w_\infty)e^{-\beta_{ij}t}$, where
$w_\infty=0.5$,
$(w_{12}^0,w_{23}^0,w_{34}^0,w_{41}^0)=(1,1.2,0.9,1.1)$, and
$(\beta_{12},\beta_{23},\beta_{34},\beta_{41})
=10^{-3}(0.5,0.8,1,0.6)$. Since $w_{ij}^0>w_\infty>0$ and
$\beta_{ij}>0$ for every edge, these weights are positive,
monotonically decreasing, and satisfy
$\lim_{t\to\infty}w_{ij}(t)=w_\infty$. Consequently, they satisfy
Assumption~\ref{ass:monotone-network-weights}. We take
$\kappa=-10^{-3}$.

The initial attitudes are given by $R_i(0)=\exp(\widehat r_i)$, where
$r_1=(0.30,-0.14,0.08)^\top$,
$r_2=(-0.20,0.24,-0.10)^\top$,
$r_3=(0.12,0.17,-0.23)^\top$, and
$r_4=(-0.16,-0.20,0.18)^\top$. Their norms are smaller than
$0.4\,\mathrm{rad}$, so the initial attitudes belong to a common normal
neighborhood of the identity. The initial body angular velocities are
$\Omega_1(0)=(0.16,-0.07,0.05)^\top$,
$\Omega_2(0)=(0.09,0.12,-0.04)^\top$,
$\Omega_3(0)=(-0.08,0.06,0.10)^\top$, and
$\Omega_4(0)=(0.11,-0.05,0.07)^\top$. The initial momenta are obtained
from the continuous Legendre transformation,
$\mu_i(0)=\mathbb I_i\Omega_i(0)$.

Before considering long times, we examine the
second-order accuracy of the complete extended
network map, including the conjugate momentum
$\nu$. This provides a numerical assessment of
the order established in
Corollary~\ref{cor:network-BEA}, independently
of the long-time variation of the extended
Hamiltonian considered below.

We use the initial data and parameters
specified above, with $T=10$ and step sizes
$h\in\{0.08,0.04,0.02,0.01,0.005\}$.
The conjugate momentum is initialized by
$\nu_0=-H(0,\mathbf R_0,\boldsymbol\mu_0)$,
as in the long-time experiment below.
An independent reference solution of the
continuous extended Hamilton equations is
computed with the eighth-order adaptive
DOP853 method, using relative and absolute
tolerances $3\times10^{-13}$ and
$10^{-14}$, respectively.

We measure the final-time configuration,
body-momentum, and conjugate-momentum errors by
\begin{equation}
\label{eq:network-extended-convergence-errors}
\begin{aligned}
e_R(h)
&:=
\left(
\sum_{i=1}^{N}
\left\|
\left(
\log\bigl(
R_{i,\mathrm{ref}}(T)^\top R_{i,h}(T)
\bigr)
\right)^\vee
\right\|^2
\right)^{1/2},\\
e_\mu(h)
&:=
\left(
\sum_{i=1}^{N}
\|\mu_{i,h}(T)-\mu_{i,\mathrm{ref}}(T)\|^2
\right)^{1/2},\\
e_\nu(h)
&:=
|\nu_h(T)-\nu_{\mathrm{ref}}(T)|.
\end{aligned}
\end{equation}
Here $\log$ denotes the principal matrix
logarithm on $SO(3)$ near the identity.
For $\alpha\in\{R,\mu,\nu\}$, the observed
order is
\[
p_\alpha(h)
:=
\log_2\left(
\frac{e_\alpha(2h)}{e_\alpha(h)}
\right).
\]

\begin{table}[htbp]
\centering
\small
\begin{tabular}{c|cc|cc|cc}
\hline
$h$ & $e_R(h)$ & $p_R$
& $e_\mu(h)$ & $p_\mu$
& $e_\nu(h)$ & $p_\nu$\\
\hline
$0.080$ & $5.4795\times10^{-2}$ & ---
& $7.0630\times10^{-2}$ & ---
& $3.3795\times10^{-5}$ & ---\\
$0.040$ & $1.3515\times10^{-2}$ & 2.0195
& $1.8911\times10^{-2}$ & 1.9011
& $8.3869\times10^{-6}$ & 2.0106\\
$0.020$ & $3.3622\times10^{-3}$ & 2.0071
& $4.8091\times10^{-3}$ & 1.9754
& $2.0905\times10^{-6}$ & 2.0043\\
$0.010$ & $8.3944\times10^{-4}$ & 2.0019
& $1.2074\times10^{-3}$ & 1.9939
& $5.2221\times10^{-7}$ & 2.0012\\
$0.005$ & $2.0979\times10^{-4}$ & 2.0005
& $3.0216\times10^{-4}$ & 1.9985
& $1.3052\times10^{-7}$ & 2.0003\\
\hline
\end{tabular}
\caption{Fixed-time global errors and
observed convergence orders of the complete
extended network integrator at $T=10$,
using an independent DOP853 reference.}
\label{tab:network-extended-state-convergence}
\end{table}

\begin{figure}[htbp]
\centering
\includegraphics[width=0.78\linewidth]
{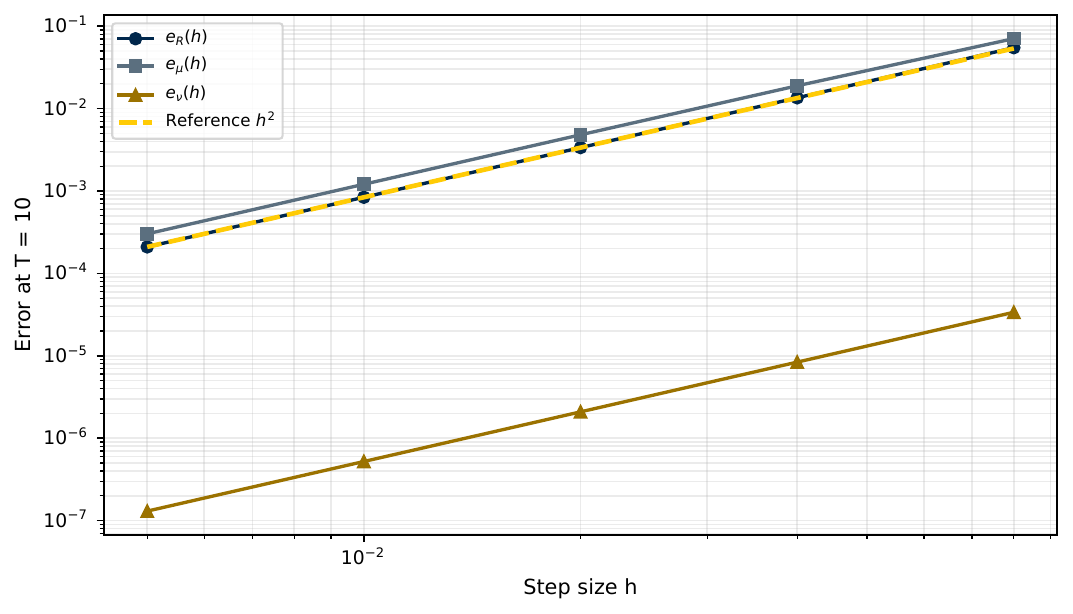}
\caption{Fixed-time convergence of the
configuration, body momenta, and
conjugate-time momentum of the
midpoint-Cayley network integrator
at $T=10$. The dashed line indicates
a reference slope proportional to $h^2$.}
\label{fig:network-extended-state-convergence}
\end{figure}

Table~\ref{tab:network-extended-state-convergence}
and Figure~\ref{fig:network-extended-state-convergence}
show observed orders approaching two in
all three components of the extended
phase-space map. These results are
consistent with the second-order accuracy
established in Corollary~\ref{cor:network-BEA}.
Unlike the long-time experiments below,
this test directly measures the global
solution error on a fixed time interval.

For the primary experiment, the step size is $h=5\times10^{-2}$ and
the final time is $T=10^4$, corresponding to $2\times10^5$ time steps.
No projection or re-orthogonalization is applied. The implicit
discrete equations are solved with a relative residual tolerance of
$2\times10^{-14}$.

Let $\boldsymbol\Omega_k=(\Omega_{1,k},\ldots,\Omega_{N,k})$ and let
$\displaystyle{\Pi_k=\sum_{i=1}^N R_{i,k}\mu_{i,k}}$ denote the total spatial
momentum. 
We monitor the following quantities:
\begin{equation*}
\label{eq:network-numerical-diagnostics}
\begin{aligned}
\varepsilon_{\mathcal C}^k
&=
\sum_{(i,j)\in\mathcal E}
\left(3-\operatorname{tr}(R_{j,k}^{\top}R_{i,k})\right),
&
\varepsilon_\Omega^k
&=
\left(
\sum_{i=1}^N\|\Omega_{i,k}\|^2
\right)^{1/2},
\\
\varepsilon_\Pi^k
&=
\frac{\|\Pi_k-\Pi_0\|}{\|\Pi_0\|},
&
\varepsilon_R^k
&=
\max_{1\leq i\leq N}
\|R_{i,k}^{\top}R_{i,k}-I\|_{\mathrm F}.
\end{aligned}
\end{equation*}
Here $\varepsilon_{\mathcal C}^k$ is the attitude-consensus error,
$\varepsilon_\Omega^k$ is the norm of the collection of body angular
velocities, $\varepsilon_\Pi^k$ is the relative total
spatial-momentum error, and $\varepsilon_R^k$ measures the maximum
orthogonality error among the agents.

\begin{figure}[h!]
\centering
\includegraphics[width=\linewidth]{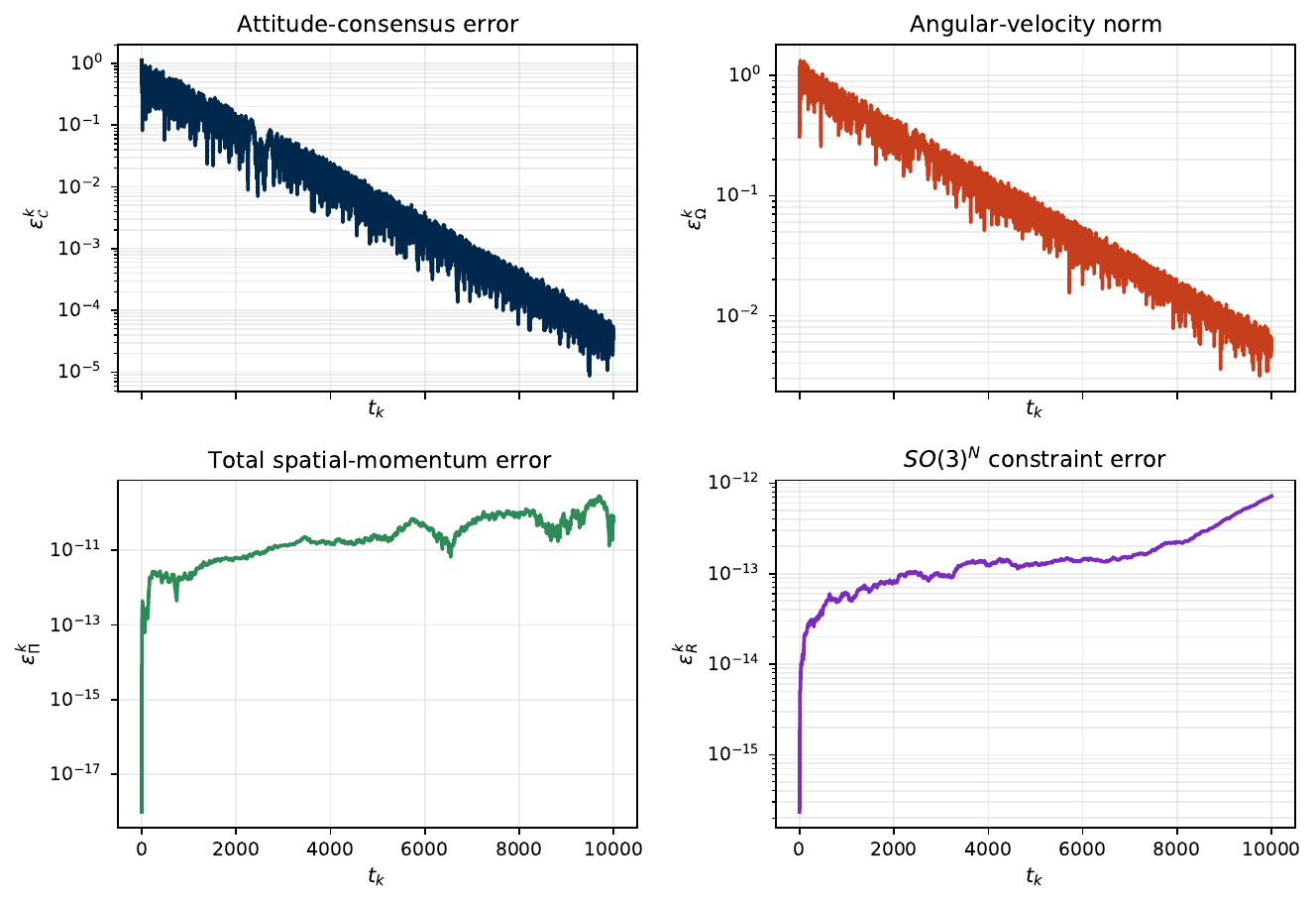}
\caption{Long-time numerical results for the time-dependent
mechanical attitude-consensus system on $SO(3)^4$, computed with
$h=5\times10^{-2}$ and $T=10^4$. Top left: attitude-consensus error
$\varepsilon_{\mathcal C}^k$. Top right: angular-velocity
$\varepsilon_\Omega^k$. Bottom left: relative total spatial-momentum
error $\varepsilon_\Pi^k$. Bottom right: maximum $SO(3)^4$ constraint
error $\varepsilon_R^k$.}
\label{fig:network-long-time-diagnostics}
\end{figure}

Figure~\ref{fig:network-long-time-diagnostics} shows convergence
towards the consensus manifold over the complete integration interval.

The consensus error decreases from 1.125 to
$3.736\times10^{-5}$, while the final
angular-velocity norm is $5.210\times10^{-3}$.
 The maximum relative spatial-momentum error is
$3.071\times10^{-10}$, and the maximum orthogonality error is
$7.184\times10^{-13}$. The maximum relative residual of the implicit
discrete equations is $2.000\times10^{-14}$. These values confirm that
the momentum and group-preservation properties remain numerically
effective over $2\times10^5$ consecutive updates. The maximum relative
attitude angle along the four edges at the final time is
$3.750\times10^{-3}\,\mathrm{rad}$, approximately $0.215^\circ$.

To study the extended Hamiltonian, let
$H^{\mathrm{ext}}=\nu+H(t,\boldsymbol R,\boldsymbol\mu)$ and initialize
$\nu_0=-H(0,\boldsymbol R_0,\boldsymbol\mu_0)$, so that
$H^{\mathrm{ext}}_0=0$. We monitor its normalized variation through
\begin{equation}
\label{eq:network-extended-Hamiltonian-diagnostic}
\varepsilon_{\mathrm{ext}}^k
=
\frac{
H^{\mathrm{ext}}(t_k,\nu_k,\boldsymbol R_k,\boldsymbol\mu_k)
-
H^{\mathrm{ext}}(t_0,\nu_0,\boldsymbol R_0,\boldsymbol\mu_0)
}{
1+
\left|
H^{\mathrm{ext}}(t_0,\nu_0,\boldsymbol R_0,\boldsymbol\mu_0)
\right|
}.
\end{equation}

With the chosen initialization, the denominator
in \eqref{eq:network-extended-Hamiltonian-diagnostic}
equals one. Moreover, summing the conjugate-momentum
update \eqref{eq:network-nu-update} gives
\[
\begin{aligned}
\varepsilon_{\mathrm{ext}}^k
={}&H(t_k,\mathbf R_k,\boldsymbol\mu_k)
-H(t_0,\mathbf R_0,\boldsymbol\mu_0)+\sum_{j=0}^{k-1}
\partial_t L_d(t_j,\mathbf R_j,\mathbf R_{j+1};h).
\end{aligned}
\]
Thus, the measured variation is precisely the
discrete energy-balance defect appearing in
Corollary~\ref{cor:long-time-discrete-energy-balance}.
It concerns the original extended Hamiltonian,
rather than the formal modified Hamiltonian
constructed in Subsection~5.3.

\begin{figure}[h!]
\centering
\includegraphics[width=0.82\linewidth]
{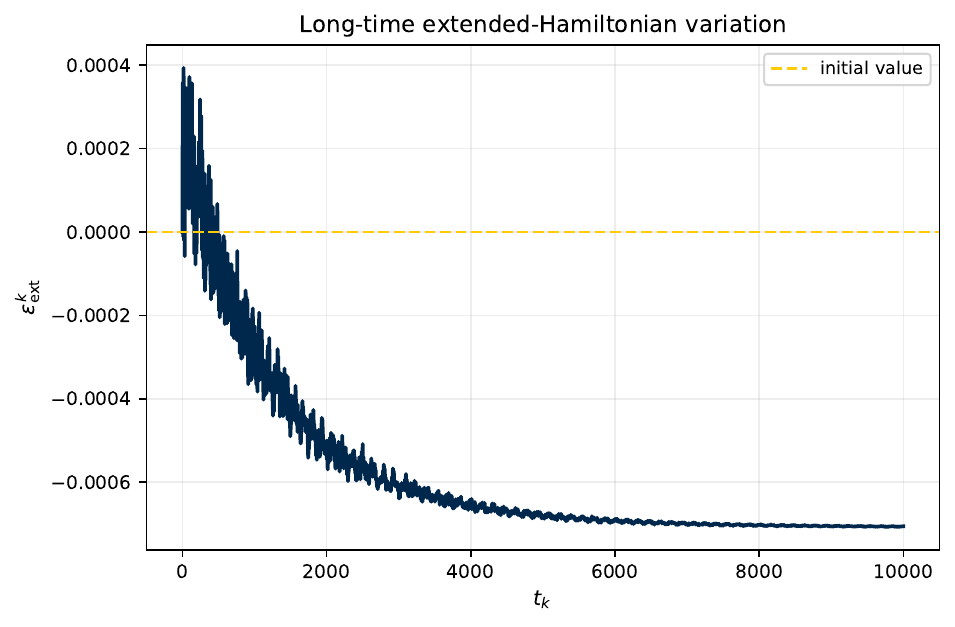}
\caption{Normalized variation
$\varepsilon_{\mathrm{ext}}^k$ of the original extended Hamiltonian
over the primary long-time experiment with
$h=5\times10^{-2}$ and $T=10^4$. The variation remains bounded and
approaches a value of approximately $-7.06\times10^{-4}$ without
subsequent secular growth.}
\label{fig:network-extended-Hamiltonian-long-time}
\end{figure}

\begin{figure}[h!]
\centering
\includegraphics[width=0.82\linewidth]
{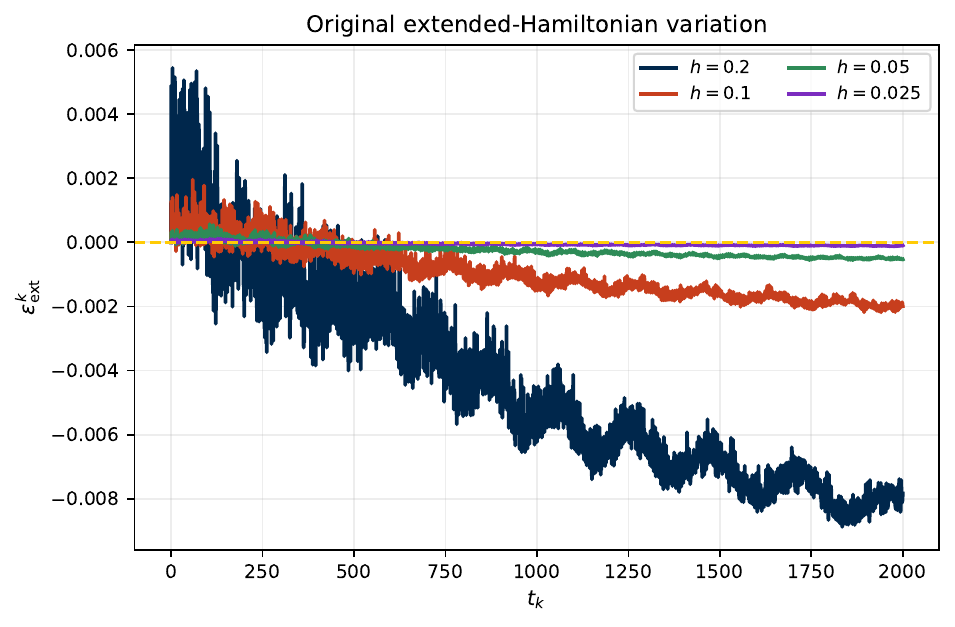}
\caption{Normalized variation of the original extended Hamiltonian for
$h=0.2$, $0.1$, $0.05$, and $0.025$ over the fixed interval
$0\leq t\leq2000$. The magnitude of the variation decreases
systematically as the step size is reduced.}
\label{fig:network-extended-Hamiltonian-step-study}
\end{figure}

For the step-size study, the system is integrated up to $T=2000$ with
$h\in\{0.2,0.1,0.05,0.025\}$. Figure
\ref{fig:network-extended-Hamiltonian-step-study} shows the normalized
extended-Hamiltonian variation for the four step sizes.

Define
$D_{\mathrm{ext}}(h)
=\max_k|\varepsilon_{\mathrm{ext}}^k|$ and the observed order by
$p_{\mathrm{ext}}(h_j)
=\log(D_{\mathrm{ext}}(h_{j-1})/D_{\mathrm{ext}}(h_j))/
\log(h_{j-1}/h_j)$. The computed values are reported in
Table~\ref{tab:network-extended-Hamiltonian-order}.

\begin{table}[h!]
\centering
\begin{tabular}{ccc}
\hline
$h$ & $D_{\mathrm{ext}}(h)$ & $p_{\mathrm{ext}}$\\
\hline
$0.200$ & $8.876\times10^{-3}$ & ---\\
$0.100$ & $2.193\times10^{-3}$ & $2.0171$\\
$0.050$ & $5.662\times10^{-4}$ & $1.9535$\\
$0.025$ & $1.413\times10^{-4}$ & $2.0023$\\
\hline
\end{tabular}
\caption{Maximum normalized variation of the original extended
Hamiltonian and observed order over $0\leq t\leq2000$.}
\label{tab:network-extended-Hamiltonian-order}
\end{table}

\begin{figure}[h!]
\centering
\includegraphics[width=0.72\linewidth]
{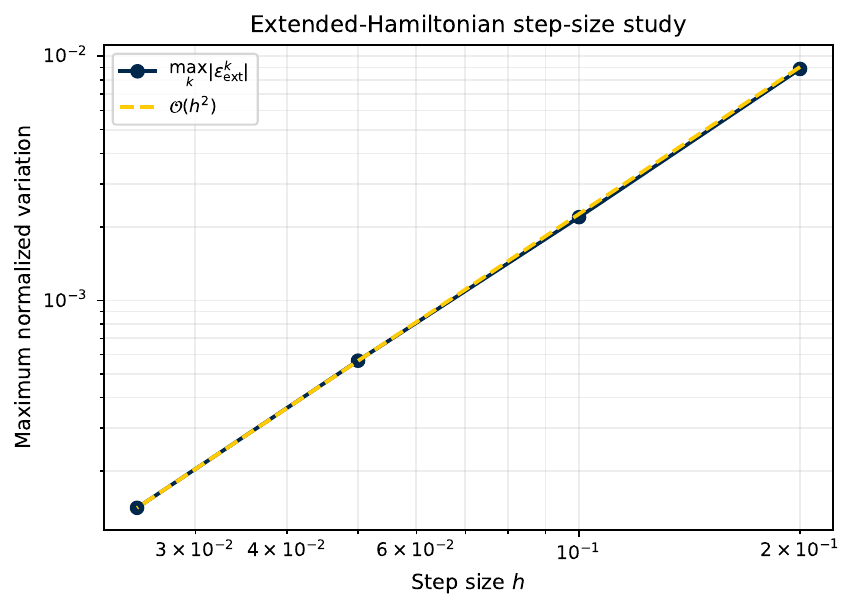}
\caption{Log-log representation of
$D_{\mathrm{ext}}(h)$ against the step size. The numerical values
follow the reference slope $\mathcal O(h^2)$, in agreement with the
second-order accuracy of the midpoint-Cayley variational integrator.}
\label{fig:network-extended-Hamiltonian-order}
\end{figure}

The observed orders in
Table~\ref{tab:network-extended-Hamiltonian-order}
are close to two, and
Figure~\ref{fig:network-extended-Hamiltonian-order}
shows an approximately quadratic dependence
of $D_{\mathrm{ext}}(h)$ on the step size
over the tested range.

This behavior is consistent with the
$O(h^2)$ long-time bound
\eqref{eq:original-H-ext-error}
under the hypotheses of
Corollary~\ref{cor:network-BEA}.
However, the experiment measures the
variation of the original extended
Hamiltonian, not that of its optimally
truncated modified counterpart.
Consequently, these numerical results do
not independently establish the uniform
analytic assumptions of the corollary
or verify the exponentially small
interpolation defect.

A video of the complete long-time evolution is available at
\url{https://youtu.be/xaXW6O1S-gw}. 
The animation displays the computed body
attitudes of the four agents, together with
the evolution of the consensus error,
angular-velocity norm, and original
extended-Hamiltonian variation.

The body frames are shown at fixed display
positions, which do not represent
translational degrees of freedom.
The numerical states are sampled
nonuniformly in time to resolve both the
initial transient and the long-time
approach to attitude consensus.

\section{Conclusions and future work}
\label{sec: con}

We have developed a variational framework for time-dependent mechanical
systems on Lie groups, from left-trivialized Euler--Lagrange equations,
energy balance, and Noether laws to discrete geometric integration.
Explicit time dependence changes the energy evolution but does not
destroy momentum conservation associated with time-independent
symmetries of the Lagrangian.

Time-indexed discrete Lagrangians admit variational and lift-based
characterizations, while retractions provide intrinsic, computable
group increments. The underlying momentum update and Noether
preservation do not depend on the particular retraction. On their regularity domains, the discrete maps on $T^*G$ are
symplectic. For left-invariant discrete Lagrangians, the reduced
maps preserve the Lie-Poisson bracket and the coadjoint orbits.
When $G$ is connected, they also preserve every Casimir function. 
Smooth variational contact gives the required control of discrete
Legendre transformations and establishes second-order accuracy of the
midpoint-Cayley scheme. Rigid-body experiments confirm this order and
illustrate long-time preservation of the group constraint, Casimir,
and spatial momentum. Exponential and periodic time modulations provide
complementary tests of the extended energy balance.

For backward error analysis, we constructed an exactly symplectic,
$\rho$-projectable map on $T^*(\mathbb R\times G)$. Its formal modified
Hamiltonian is affine in the momentum conjugate to time and, under
symmetry, admits reduction. Under uniform analyticity, coefficient
bounds, and domain-containment hypotheses, optimal truncation gives
exponentially accurate interpolation and near-conservation of the
modified extended Hamiltonian over exponentially long intervals.
The associated discrete energy-balance estimate concerns the original
extended Hamiltonian, not conservation of the physical Hamiltonian.
These conclusions are conditional on uniform control along the
trajectory and are not global statements for arbitrary Lie groups.

For the mechanical attitude network on $SO(3)^N$, we prove local
asymptotic stability of the static consensus set for the continuous
dynamics when $\kappa<0$, the graph is connected, and the coupling
weights are non-increasing and uniformly bounded away from zero. The diagonal $SO(3)$ symmetry preserves total spatial
momentum in the continuous and discrete models, but the interactions
generally prevent conservation of individual body momenta or their
quadratic functions. The same symmetry restricts the formal modified
Hamiltonian to relative attitudes and body momenta. Numerical tests
confirm second-order convergence of the complete extended map,
including the momentum conjugate to time, and illustrate attitude
synchronization, angular-velocity decay, Noether conservation,
group-constraint preservation, and bounded variation of the original
extended Hamiltonian.

Future directions include oscillatory, switching, and state-dependent
couplings, for which persistent excitation, joint connectivity, or
switched-system methods would replace the present monotonicity
argument. Global consensus is not claimed: the topology of $SO(3)^N$
and the additional critical configurations of the potential motivate
almost-global results, alternative navigation functions, and estimates
of the local region of attraction. Further numerical developments
include higher-order discrete Lagrangians and retractions,
variable-time variational integrators, scalable implicit solvers,
roundoff-aware implementations, and explicit low-order modified
Hamiltonians. Extension to noncompact groups such as $SE(3)$ requires
new coercivity and domain-control arguments and mechanically
appropriate Lagrangians and interactions; it does not follow from a
formal change of configuration group.

\section*{Acknowledgments}

The authors acknowledge financial support from the Spanish
Ministry of Science and Innovation under grant
PID2022-137909NB-C21, PCI2024-155047-2 and from the Severo Ochoa Programme
for Centres of Excellence in R\&D (CEX2023-001347-S).

\medskip
\noindent\textbf{Use of artificial intelligence.}
Generative AI tools were used to review and improve the
writing of the manuscript and to assist with code for the
numerical experiments. The authors take full responsibility
for the content of the manuscript, including all
mathematical statements, proofs, and numerical results.

\end{document}